\documentclass[11pt]{amsart}
\usepackage{amscd,amssymb,amsthm,amsmath,amssymb,textcomp,supertabular,longtable,enumerate,rotating,mathrsfs,mathtools,hyperref,latexsym,amscd,amsbsy,mathrsfs}
\usepackage{pdflscape}
\usepackage{graphicx}
\usepackage{array}
\usepackage[matrix,arrow,curve]{xy}
\usepackage{xcolor}
\usepackage{longtable}
\usepackage{enumitem}
\usepackage{MnSymbol}
\usepackage{multicol}
\usepackage{tabulary}
\usepackage{colortbl}
\usepackage{geometry}

\usepackage[
    maxbibnames=999,    
    maxcitenames=1,     
    mincitenames=1,     
    style=apa,    
    url=false,
    doi=true,
]{biblatex}
\newtheorem{lemma}[equation]{Lemma}

\newtheorem*{calabiproblem*}{Calabi Problem}

\theoremstyle{definition}

\makeatletter\@addtoreset{equation}{section} \makeatother

\swapnumbers

\newtheoremstyle{dotless}{}{}{\rm}{}{\sc}{}{ }{}

\theoremstyle{dotless}

\newcommand{\DR}{\mathbb{R}} 

\newcommand{\DP}{\mathbb{P}}
\newcommand{\DA}{\mathbb{A}}

\newcommand{\MA}{\mathcal{A}}

\newtheorem*{theorem*}{Theorem}
\newtheorem*{maintheorem*}{Main Theorem.}

\author{Elena Denisova}

\title{$\delta$-invariants of Du Val del Pezzo surfaces of degree $\ge 4$}

\address{\emph{Elena Denisova}
\newline
\textnormal{School of Mathematics, The University of Edinburgh, Edinburgh EH9 3JZ, UK.}
\newline
\textnormal{\texttt{e.denisova@sms.ed.ac.uk}}}

\begin{document}

\maketitle
\begin{abstract}
In this article, we compute $\delta$-invariants of Du Val del Pezzo surfaces of degrees $4$, $5$, $6$, $7$ and $8$.
\end{abstract}
\section{Introduction}
\subsection{History and Results.} 
It is known that a smooth Fano variety admits a Kähler--Einstein metric if and only if it is $K$-polystable.
In the case of two-dimensional Fano varieties (del Pezzo surfaces), Tian and Yau showed that a smooth del Pezzo surface is $K$-polystable if and only if it is not the blow-up of $\mathbb{P}^2$ at one or two points (see \cite{TianYau1987,Ti90}).
Significant progress has been made for smooth Fano threefolds (see 
 \cite{AbbanZhuangSeshadri, Fano21, Liu23, CheltsovFujitaKishimotoPark23, LiuZhao24, GuerreiroGiovenzanaViswanathan23, Malbon24, CheltsovPark22, BelousovLoginov24, BelousovLoginov23, CheltsovFujitaKishimotoOkada23, Denisova24, Denisova23, CheltsovDenisovaFujita24}).
However, for Fano varieties in higher dimensions, many questions remain open. In several cases, the problem reduces to computing the $\delta$-invariant of (possibly singular) del Pezzo surfaces (see \cite{Fano21,CheltsovDenisovaFujita24,CheltsovFujitaKishimotoOkada23}, etc.). For smooth del Pezzo surfaces  $\delta$-invariants were computed in \cite{Fano21}:
 \begin{table}[h!]
 \renewcommand{\arraystretch}{1.6}
\centering
\begin{tabular}{ | c | c | c | c | c | c | c | c | c | c | c | c |c | c |}
 \hline 
 $X$ & $\DP^2$ & $\DP^1\times \DP^1$ & $S_8$  & $S_7$ & $S_6$ & $S_5$ & $S_4$ & $S_3^1$ & $S_3^2$ & $S_2^1$  & $S_2^2$ & $S_1^1$ & $S_1^2$ \\\hline
$\delta(X)$ & $1$ & $1$ & $\frac{6}{7}$ & $\frac{21}{25}$ & $1$ & $\frac{15}{13}$ & $\frac{4}{3}$ & $\frac{3}{2}$ & $\frac{27}{17}$ & $\frac{9}{5}$ & $\frac{15}{8}$ & $\frac{15}{7}$ & $\frac{12}{5}$\\\hline
    \end{tabular}
\end{table}
 \\where $S_d$ is a blowup of $\DP^2$ at $9-d$ points in general position; $S_3^1$ is $S_3$ with an Eckardt point, $S_3^2$ is $S_3$ without an Eckardt point; $S_2^1$ is $S_2$ such that
the linear system $|-K_{S_2}|$ contains a tacnodal curve;  $S_2^2$ is $S_2$ such that
the linear system $|-K_{S_2}|$ does not contain a tacnodal curve; $S_1^1$ is $S_1$ such that
the linear system $|-K_{S_1}|$ contains a cuspidal curve;  $S_1^2$ is $S_1$ such that
the linear system $|-K_{S_1}|$ does not contain a cuspidal curve. In this thesis, we compute  $\delta$-invariants  of singular Du Val del Pezzo surfaces. In this article, we compute $\delta$-invariants  of Du Val del Pezzo surfaces of degree $d$ for $d\in\{4,5,6,7,8\}$.
\begin{maintheorem*}
Let $X$ be a singular del Pezzo surface of degree $d$ where $d\in\{4,5,6,7,8\}$. Then the  $\delta$-invariant of $X$ is uniquely determined by the degree of $X$, the number of lines on $X$, and the type of singularities on $X$ which is given in the following table:\\
\begin{minipage}{5.7cm}
 \renewcommand{\arraystretch}{1.1}
  \begin{longtable}{ | c | c | c | c | }
   \hline
   $K_X^2$ & $\#$ lines & $\mathrm{Sing}(X)$ & $\delta(X)$\\
  \hline\hline
\endhead 
 $8$ & $0$ & $\DA_1$ & $\frac{3}{4}$\\
\hline
  $7$ & $2$ & $\DA_1$ & $\frac{21}{31}$\\
\hline
 $6$ & $3$ & $\DA_1$ & $\frac{3}{4}$\\
\hline
 $6$ & $4$ & $\DA_1$ & $\frac{9}{11}$\\
\hline
 $6$ & $2$ & $2\DA_1$ & $\frac{9}{14}$\\
\hline
 $6$ & $2$ & $\DA_2$ & $\frac{3}{5}$\\
\hline
 $6$ & $1$ & $\DA_2+\DA_1$ & $\frac{1}{2}$\\
\hline
 $5$ & $7$ & $\DA_1$ & $\frac{15}{17}$\\
\hline
 $5$ & $5$ & $2\DA_1$ & $\frac{15}{19}$\\
\hline 
  \end{longtable}
  \end{minipage}
  \begin{minipage}{6.7cm}
 \renewcommand{\arraystretch}{1.0}
    \begin{longtable}{ | c | c | c | c | }
   \hline
   $K_X^2$ & $\#$ lines & $\mathrm{Sing}(X)$ & $\delta(X)$\\
  \hline\hline
\endhead 
 $5$ & $4$ & $\DA_2$ & $\frac{5}{7}$\\
\hline
$5$ & $3$ & $\DA_2+\DA_1$ & $\frac{15}{23}$\\
\hline
$5$ & $2$ & $\DA_3$ & $\frac{5}{9}$\\
\hline
 $5$ & $1$ & $\DA_4$ & $\frac{3}{7}$\\
\hline
  $4$ & $12$ & $\DA_1$ & $1$\\
\hline
 $4$ & $9$ & $2\DA_1$ & $1$\\
\hline
 $4$ & $8$ & $2\DA_1$ & $1$\\
\hline
 $4$ & $6$ & $3\DA_1$ & $1$\\
\hline
 $4$ & $4$ & $4\DA_1$ & $1$\\
\hline
 $4$ & $8$ & $\DA_2$ & $\frac{6}{7}$\\
\hline

  \end{longtable}
  \end{minipage}
    \begin{minipage}{5.5cm}
 \renewcommand{\arraystretch}{1.1}
    \begin{longtable}{ | c | c | c | c | }
   \hline
   $K_X^2$ & $\#$ lines & $\mathrm{Sing}(X)$ & $\delta(X)$\\
  \hline\hline
\endhead 
 $4$ & $6$ & $\DA_2+\DA_1$ & $\frac{6}{7}$\\
\hline
 $4$ & $4$ & $\DA_2+2\DA_1$ & $\frac{6}{7}$\\
\hline
 $4$ & $5$ & $\DA_3$ & $\frac{2}{3}$\\
\hline
 $4$ & $4$ & $\DA_3$ & $\frac{3}{4}$\\
\hline
 $4$ & $3$ & $\DA_3+\DA_1$ & $\frac{3}{4}$\\
\hline
 $4$ & $2$ & $\DA_3+2\DA_1$ & $\frac{3}{4}$\\
\hline
 $4$ & $3$ & $\DA_4$ & $\frac{6}{11}$\\
\hline
  $4$ & $2$ & $\mathbb{D}_4$ & $\frac{1}{2}$\\
\hline
  $4$ & $1$ & $\mathbb{D}_5$ & $\frac{3}{8}$\\
\hline
  \end{longtable}
  \end{minipage}\\
  \end{maintheorem*}

\noindent {\bf Acknowledgments:} I am grateful to my supervisor Professor Ivan Cheltsov for the introduction to the topic and continuous support.

\section{Proof of Main Theorem via  Kento Fujita’s formulas}
\noindent In this work in order to find  $\delta$-invariants of Du Val del Pezzo surfaces we apply Abban--Zhuang theory and use Kento Fujita’s formulas. Let $X$ be a Du Val del Pezzo surface, and let $S$ be a minimal resolution of $X$. For a birational morphism  $f\colon\widetilde{X}\to X$ and $E$ be a prime divisor in $\widetilde{X}$ we say that $E$ is a prime divisor over $X$.
If~$E$~is $f$-exceptional, we say that $E$ is an~exceptional  prime divisor over~$X$.
We will denote the~subvariety $f(E)$ by $C_X(E)$. 
Let \index{$S_X(E)$}
$$
S_X(E)=\frac{1}{(-K_X)^2}\int_{0}^{\tau}\mathrm{vol}(f^*(-K_X)-vE)dv\text{ and }A_X (E) = 1 + \mathrm{ord}_E(K_{\widetilde{X}} - f^*(K_X)),
$$
where $\tau=\tau(E)$ is the~pseudo-effective threshold of $E$ with respect to $-K_X$.
Let $Q$ be a point in $X$. We can define a local $\delta$-invariant and a global $\delta$-invariant of $X$ as
$$
\delta_Q(X)=\inf_{\substack{E/X\\ Q\in C_X(E)}}\frac{A_X(E)}{S_X(E)}\text{ and }\delta(X)=\inf_{Q\in X}\delta_Q(X)
$$
where the~infimum runs over all prime divisors $E$ over the surface $X$ such that $Q\in C_X(E)$. Similarly, for the  surface $S$ and a point $P$ on $S$ we define local $\delta$-invariant and a global $\delta$-invariant of $S$ as
$$
\delta_P(S)=\inf_{\substack{F/S\\ P\in C_S(F)}}\frac{A_S(F)}{S_S(F)}
\text{
and }\delta(S)=\inf_{P\in S}\delta_P(S)$$
where $S_S(F)$ and $A_S(F)$ are defined as $S_X(E)$ and $A_X(E)$ above. It is clear that
$$\delta(X)=\delta(S)\text{ and }\delta_Q(X)=\inf_{P: Q=f(P)}\delta_P(S)$$
We now fix a point $P$ on $S$ and choose a smooth curve $A$ on $S$ containing $P$. 
Set
$$
\tau(A)=\mathrm{sup}\Big\{v\in\mathbb{R}_{\geqslant 0}\ \big\vert\ \text{the divisor  $-K_S-vA$ is pseudo-effective}\Big\}.
$$
For~$v\in[0,\tau]$, let $P(v)$ be the~positive part of the~Zariski decomposition of the~divisor $-K_S-vA$,
and let $N(v)$ be its negative part. 
Then we set $$
S\big(W^A_{\bullet,\bullet};P\big)=\frac{2}{K_S^2}\int_0^{\tau(A)} h(v) dv,
\text{ where }
h(v)=\big(P(v)\cdot A\big)\times\big(N(v)\cdot A\big)_P+\frac{\big(P(v)\cdot A\big)^2}{2}.
$$
It follows from {\cite[Theorem 1.7.1]{Fano21}} that:
\begin{equation}\label{estimation1}
    \delta_P(S)\geqslant\mathrm{min}\Bigg\{\frac{1}{S_S(A)},\frac{1}{S(W_{\bullet,\bullet}^A,P)}\Bigg\}.
\end{equation}
Unfortunately, this approach does not always give us a good estimation. If this is the case, we apply the generalization of this method. Let $\sigma: \widehat{S}\to S$ be a weighted blowup of the point $P$ on $S$. Suppose, in addition, that $\widehat{S}$ is a Mori Dream Space Then
\begin{itemize}
\item the~$\sigma$-exceptional curve $E_P\cong \DP^1$ such that $\sigma(E_P)=P$,
\item the~log pair $(\widehat{S},E_P)$ has purely log terminal singularities.
\end{itemize}
We write
$$
K_{E_P}+\Delta_{E_P}=\big(K_{\widehat{S}}+E_P\big)\big\vert_{E_P},
$$
where $\Delta_{E_P}$ is an~effective $\mathbb{Q}$-divisor on $E_P$ known as the~different of the~log pair $(\widehat{S},E_P)$.
Note that the~log pair $(E_P,\Delta_{E_P})$ has at most Kawamata log terminal singularities, and the~divisor $-(K_{E_P}+\Delta_{E_P})$ is $\sigma\vert_{E_P}$-ample.
\\Let $O$ be a point on $E_P$. 
Set
$$
\tau(E_P)=\mathrm{sup}\Big\{v\in\mathbb{R}_{\geqslant 0}\ \big\vert\ \text{the divisor  $\sigma^*(-K_S)-vE_P$ is pseudo-effective}\Big\}.
$$
For~$v\in[0,\tau]$, let $\widehat{P}(v)$ be the~positive part of the~Zariski decomposition of the~divisor $\sigma^*(-K_S)-vE_P$,
and let $\widehat{N}(v)$ be its negative part. 
Then we set $$
S\big(W^{E_P}_{\bullet,\bullet};O\big)=\frac{2}{K_{\widehat{S}}^2}\int_0^{\tau(E_P)} \widehat{h}(v) dv,
\text{ where }
\widehat{h}(v)=\big(\widehat{P}(v)\cdot E_P\big)\times\big(\widehat{N}(v)\cdot E_P\big)_O+\frac{\big(\widehat{P}(v)\cdot E_P\big)^2}{2}.
$$
Let
$A_{E_P,\Delta_{E_P}}(O)=1-\mathrm{ord}_{\Delta_{E_P}}(O)$.
It follows from {\cite[Theorem 1.7.9]{Fano21}} and {\cite[Corollary 1.7.12]{Fano21}} that
\begin{equation}
\label{estimation2}
\delta_P(S)\geqslant\mathrm{min}\Bigg\{\frac{A_S(E_P)}{S_S(E_P)},\inf_{O\in E_P}\frac{A_{E_P,\Delta_{E_P}}(O)}{S\big(W^{E_P}_{\bullet,\bullet};O\big)}\Bigg\},
\end{equation}
where the~infimum is taken over all points $O\in E_P$. Now for all the points $P$ on $S$ we know either values of local $\delta$-invariants or estimations of them. Taking the minimum we compute $\delta(S)$ --- the  global  $\delta$-invariant of $S$ and thus, $\delta(X)=\delta(S)$ --- the  global  $\delta$-invariant of $X$. We apply this method to minimal resolutions of Du Val del Pezzo surfaces to prove the Main Theorem. Throughout this work small circles correspond to a $(-1)$-curves and large circles correspond to a $(-2)$-curves on dual graphs.
\section{Du Val del Pezzo surfaces of degrees $7$ and $8$}
\noindent It was mentioned in \cite[Table 2.1]{Fano21} $\delta$-invariants of smooth del Pezzo surfaces of degrees $7$ and $8$ are given in the following table:
\begin{table}[h!]
 \renewcommand{\arraystretch}{1.4}
\centering
\begin{tabular}{ | c | c | c | c |}
 \hline 
 $X$  & $\DP^1\times \DP^1$ & blowup of $\DP^2$ in one point & blowup of $\DP^2$ in two points\\\hline
$\delta(X)$ & $1$ & $\frac{6}{7}$ & $\frac{21}{25}$\\\hline
    \end{tabular}
\end{table}
\\ In this chapter, we compute  $\delta$-invariants  of singular  Du Val del Pezzo surfaces of degrees $7$ and $8$. More precisely we prove that if $X$ is a del Pezzo surface with $\DA_1$-singularity of degree $8$ then $\delta(X)=\frac{4}{3}$ and if $X$ is a del Pezzo surface with $\DA_1$-singularity of degree $7$ then $\delta(X)=\frac{21}{31}$.

\subsection{Del Pezzo surface of degree $8$ with $\DA_1$ singularity}
\begin{lemma}
Let $X$ be a singular del Pezzo surface of degree $8$ with one $\DA_1$ singularity.
Then $X$ contains $0$ lines. Let $E$ be the unique $(-2)$-curve. One has $\delta_P(X)=\frac{3}{4}\text{ for }P\in X $. Thus $\delta(X)=\frac{3}{4}$.
\end{lemma}
\begin{proof}
 Let $F$ be a fiber of natural projection on $\DP^1$. Then $\tau(F)=4$ and the Zariski Decomposition of the divisor $-K_X-vF$ is given by:
$$P(v)=-K_X-vF-\frac{v}{2}E\text{ and }N(v)=\frac{v}{2}E\text{ if }v\in[0,4].$$
Moreover, 
$$(P(v))^2=\frac{(4-v)^2}{2}\text{ and }P(v)\cdot F=2-\frac{v}{2}\text{ if }v\in[0,4].$$
We have $S_X(F)=\frac{4}{3}$.
Thus, $\delta_P(S)\le \frac{3}{4}$. Note that 
$h(v)\le \frac{(4-v)(4+v)}{8}\text{ if }v\in[0,4]$.
So $S(W_{\bullet,\bullet}^F;P)\le\frac{4}{3}$.
Thus, $\delta_P(X)=\frac{3}{4}$ for $P\in X$ and  $\delta(X)=\frac{3}{4}$.
\end{proof}

\subsection{Del Pezzo surface of degree $7$ with $\DA_1$ singularity}
\begin{lemma}
    Let $X$ be a singular del Pezzo surface with one $\DA_1$  singularity, $S$ be a minimal resolution of $X$. Then $X$ contains $2$ lines, and $S$ can be obtained by blowing up $\DP^2$ at point $P_1$; and a point $P_2$ on the exceptional divisor corresponding to $P_1$. Let $E_1$, $E_2$ be the exceptional divisors corresponding to $P_1$, $P_2$; $L_{12}$ be a $(-1)$-curve which is a strict transform of the line passing through $P_1$.  The dual graph of $(-1)$ and $(-2)$ curves is given in the following picture:
\begin{figure}[h!]
    \begin{center}
      \includegraphics[width=4.7cm]{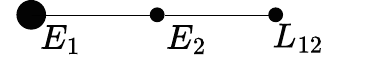}
   \end{center}
   \caption{Dual graph: $(-K_S)^2=7$, singularity $\DA_1$}
\end{figure}
\par One has
     $$\delta_P(S)=\begin{cases}
      \frac{21}{31}\text{ if }P\in E_2,\\
      \frac{7}{9}\text{ if }P\in E_1\backslash E_2,\\
      \frac{21}{25}\text{ if }P\in L_{12}\backslash E_2.\\
      \end{cases}\text{ and }
\delta_P(S) \ge\frac{21}{23}\text{ otherwise. }$$
Thus $\delta_P(X)=\frac{21}{31}.$
\end{lemma}
\begin{proof}
{\bf Step 1.} Suppose $P\in L_{12}$. Then $\tau(L_{12})=3$ and the Zariski Decomposition of the divisor $-K_S-vL_{12}$ is given by:
\begin{align*}
    &&P(v)=\begin{cases}-K_S-vL_{12}\text{ if }v\in[0,1],\\
-K_S-vL_{12}-(v-1)(E_1+2E_2)\text{ if }v\in[1,3].
\end{cases}\\
&&N(v)=\begin{cases}0\text{ if }v\in[0,1],\\
(v-1)(E_1+2E_2)\text{ if }v\in[1,3].
\end{cases}
\end{align*}
Moreover, 
$$(P(v))^2=\begin{cases}7-2v-v^2\text{ if }v\in[0,1],\\
(3-v)^2\text{ if }v\in[1,3].
\end{cases}P(v)\cdot L_{12}=\begin{cases}1+v\text{ if }v\in[0,1],\\
3-v\text{ if }v\in[1,3].
\end{cases}$$
We have  $S_S(L_{12})=\frac{25}{21}$.
Thus, $\delta_P(S)\le \frac{21}{25}$ for $P\in L_{12}$. Note that for $P\in L_{12}\backslash E_2$ we have:
$$h(v)\le 
\begin{cases} 
\frac{(1+v)^2}{2} \text{ if }v\in[0,1],\\
 \frac{(3-v)^2}{2} \text{ if }v\in[1,3].
\end{cases}$$
So 
$S(W_{\bullet,\bullet}^{L_{12}};P)\le\frac{5}{7}<\frac{25}{21}$.
Thus, $\delta_P(S)=\frac{21}{25}$ if $P\in L_{12}\backslash E_2$. 
\\{\bf Step 2.} Suppose $P\in E_{1}$ . Then $\tau(E_1)=2$ and the Zariski Decomposition of the divisor $-K_S-vE_{1}$ is given by:
$$P(v)=\begin{cases}
-K_S-vE_{1}\text{ if }v\in[0,1],\\
-K_S-vE_{1}-(v-1)E_2\text{ if }v\in[1,2].
\end{cases}N(v)=\begin{cases}0\text{ if }v\in[0,1],\\
(v-1)E_2\text{ if }v\in[1,2].
\end{cases}$$
Moreover, 
$$(P(v))^2=\begin{cases}7-2v^2\text{ if }v\in[0,1],\\
(2-v)(4+v)\text{ if }v\in[1,2].
\end{cases}P(v)\cdot E_{1}=\begin{cases}2v\text{ if }v\in[0,1],\\
v+1\text{ if }v\in[1,2].
\end{cases}$$
We have
$S_S(E_{1})=\frac{9}{7}$.
Thus, $\delta_P(S)\le \frac{7}{9}$ for $P\in E_{1}$.  Note that for $P\in E_1\backslash E_2$ 
$$h(v)\le \begin{cases}2v^2\text{ if }v\in[0,1],\\
\frac{(v+1)^2}{2}\text{ if }v\in[1,2].
\end{cases}$$
 So
$S(W_{\bullet,\bullet}^{E_1};P)\le\frac{23}{21}< \frac{9}{7}$.
Thus, $\delta_P(S)=\frac{7}{9}$ if $P\in E_{1}\backslash E_{2}$.\\
{\bf Step 3.} Suppose $P\in E_{2}$ . Then $\tau(E_2)=4$ and the Zariski Decomposition of the divisor $-K_S-vE_{2}$ is given by:
$$P(v)=\begin{cases}
-K_S-vE_{2}-\frac{v}{2}E_1\text{ if }v\in[0,1],\\
-K_S-vE_{2}-\frac{v}{2}E_1-(v-1)L_{12}\text{ if }v\in[1,4].
\end{cases}N(v)=\begin{cases}\frac{v}{2}E_1\text{ if }v\in[0,1],\\
\frac{v}{2}E_1+(v-1)L_{12}\text{ if }v\in[1,4].
\end{cases}$$
Moreover, 
$$(P(v))^2=\begin{cases}7-2v-\frac{v^2}{2}\text{ if }v\in[0,1],\\
\frac{(4-v)^2}{2}\text{ if }v\in[1,4].
\end{cases}P(v)\cdot E_{2}=\begin{cases}1+\frac{v}{2}\text{ if }v\in[0,1],\\
2-\frac{v}{2}\text{ if }v\in[1,4].
\end{cases}$$
We have
$S_S(E_{2})=\frac{31}{21}$.
Thus, $\delta_P(S)\le \frac{21}{31}$ for $P\in E_{2}$.  Note that if $P\in E_2\backslash L_{12}$ or if $P\in E_2\cap L_{12}$ then
$$h(v)\le\begin{cases}
 \frac{(v + 2) (3 v + 2)}{8}\text{ if }v\in[0,1],\\
\frac{ (4 - v) (v + 4)}{8}\text{ if }v\in[1,4].
\end{cases}
\text{ or }
h(v)\le\begin{cases}
\frac{(2+v)^2}{8}\text{ if }v\in[0,1],\\
\frac{3 (4 - v) v}{8}\text{ if }v\in[1,4].
\end{cases}$$
So
$S(W_{\bullet,\bullet}^{E_2};P)\le\frac{17}{21} < \frac{31}{21}$
or  $S(W_{\bullet,\bullet}^{E_2};P)\le\frac{13}{21} \le \frac{31}{21}$.
Thus, $\delta_P(S)=\frac{21}{31}$ if $P\in E_{2}$.
\\{\bf Step 4.} Suppose $P\not\in  E_1\cup E_2\cup L_2$. Consider a blowup $\pi:\widetilde{S}\to S$ at point $P$ with the exceptional divisor $E_P$. Suppose $\widetilde{E}_1$ is a strict transform of $E_1$  and $L_{1P}$ is a strict transform of the line through $P_1$ and a projection of $P$ on $\DP^2$. Then  $\tau(E_P)=3$ and the Zariski Decomposition of the divisor $\sigma^*(-K_S)-vE_{P}$ is given by:
\begin{align*}
&&P(v)=\begin{cases}\sigma^*(-K_S)-vE_{P}\text{ if }v\in[0,2],\\
\sigma^*(-K_S)-vE_{P}-(v-2)(2L_{1P}+\widetilde{E}_1)\text{ if }v\in[2,3].
\end{cases}\\
&&N(v)=\begin{cases}0\text{ if }v\in[0,2],\\
(v-2)(2L_{1P}+\widetilde{E}_1)\text{ if }v\in[2,3].
\end{cases}
\end{align*}
Moreover, 
$$(P(v))^2=\begin{cases}7-v^2\text{ if }v\in[0,2],\\
(3-v)(5-v)\text{ if }v\in[2,3].
\end{cases}P(v)\cdot E_{P}=\begin{cases}v\text{ if }v\in[0,2],\\
4-v\text{ if }v\in[2,3].
\end{cases}$$
We have $S_S(E_{P})=\frac{38}{21}$.
Thus, $\delta_P(S)\le \frac{2}{38/21}=\frac{21}{19}$.  Note that for $O\in E_P$ we have: 
$$h(v)=\begin{cases}\frac{v^2}{2}\text{ if }v\in[0,2],\\
    \frac{(4 - v) (3 v - 4)}{2}\text{ if }v\in[2,3].
\end{cases}$$
So
$S(W_{\bullet,\bullet}^{E_P};O)\le\frac{23}{21}$.
Thus, $\delta_P(S)\le \frac{21}{23}$ if $P\not\in E_{1}\cup E_2\cup L_{12}$.
\end{proof}

\section{Du Val Del Pezzo surfaces of degree 6}
\noindent It was mentioned in \cite[Table 2.1]{Fano21} that $\delta(X)=1$ when $X$ is a smooth del Pezzo surface of degree $6$. In this chapter, we compute  $\delta$-invariants  of singular  Du Val del Pezzo surfaces of degree $6$. 
\begin{maintheorem*}
Let $X$ be a singular del Pezzo surface of degree $6$. Then the  $\delta$-invariant of $X$ is uniquely determined by the degree of $X$, the number of lines on $X$, and the type of singularities on $X$ which is given in the following table:
 {
  \begin{longtable}{ | c | c | c | c | }
   \hline
   $K_X^2$ & $\#$ lines & $\mathrm{Sing}(X)$ & $\delta(X)$\\
  \hline\hline
\endhead 
 $6$ & $3$ & $\DA_1$ & $\frac{3}{4}$\\
\hline
 $6$ & $4$ & $\DA_1$ & $\frac{9}{11}$\\
\hline
$6$ & $2$ & $2\DA_1$ & $\frac{9}{14}$\\
\hline
$6$ & $2$ & $\DA_2$ & $\frac{3}{5}$\\
\hline
 $6$ & $1$ & $\DA_2+\DA_1$ & $\frac{1}{2}$\\
 \hline
  \end{longtable}}
  \end{maintheorem*}
\subsection{General results for degree $6$}
Let $X$ be a del Pezzo surface of degree $6$ with at most Du Val singularities,  $S$ be a minimal resolution of $X$ and $P$ be a point on $S$.   Then: 
\begin{lemma}\label{deg6-generalpoint}
Let $P$ be a general point on $S$. Consider the blowup $\sigma:\widetilde{S}\to S$ of $S$ at $P$ with the exceptional divisor $A$ there exist $(-1)$-curves and $(-2)$-curves  which form  one of the following dual graphs:
\begin{figure}[h]
    \centering
    \includegraphics[width=12cm]{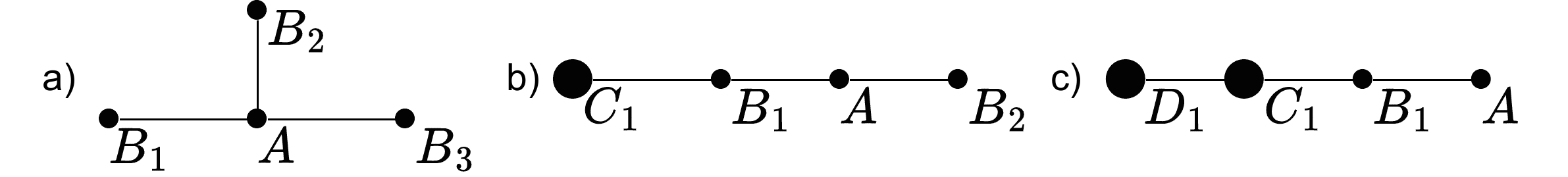}
    \caption{Dual graph: $(-K_S)^2=6$ for a general point}
\end{figure}
\par Then $\tau(A)=3$ and the Zariski Decomposition of the divisor $\sigma^*(-K_S)-vA$ is given by:
{{\allowdisplaybreaks \begin{align*}
&{\text{\bf a). }}&P(v)=\begin{cases}
\sigma^*(-K_{S})-vA\text{ if }v\in[0,2],\\
\sigma^*(-K_{S})-vA-(v-2)(B_1+B_2+B_3)\text{ if }v\in[2,3].
\end{cases}\\&&
N(v)=\begin{cases}
0\text{ if }v\in[0,2],\\
(v-2)(B_1+B_2+B_3)\text{ if }v\in[2,3].
\end{cases}\\
&{\text{\bf b). }}&P(v)=\begin{cases}
\sigma^*(-K_{S})-vA\text{ if }v\in[0,2],\\
\sigma^*(-K_{S})-vA-(v-2)(2B_1+C_1+B_2)\text{ if }v\in[2,3].
\end{cases}\\&&
N(v)=\begin{cases}
0\text{ if }v\in[0,2],\\
(v-2)(2B_1+C_1+B_2)\text{ if }v\in[2,3].
\end{cases}\\
&{\text{\bf c). }}&P(v)=\begin{cases}
\sigma^*(-K_{S})-vA\text{ if }v\in[0,2],\\
\sigma^*(-K_{S})-vA-(v-2)(3B_1+2C_1+D_1)\text{ if }v\in[2,3].
\end{cases}\\&&
N(v)=\begin{cases}
0\text{ if }v\in[0,2],\\
(v-2)(3B_1+2C_1+D_1)\text{ if }v\in[2,3].
\end{cases}
\end{align*}}}
Moreover, 
$$(P(v))^2=
\begin{cases}
6-v^2\text{ if }v\in[0,2],\\
2(3-v)^2\text{ if }v\in[2,3].
\end{cases} 
P(v)\cdot A=
\begin{cases}
v\text{ if }v\in[0,2],\\
2(3-v)\text{ if }v\in[2,3].
\end{cases}$$
In this case  $\delta_P(S)\ge 1$.
\end{lemma}
\begin{proof}
 The Zariski Decomposition in part a). follows from 
 $\sigma^*(-K_S)-vA\sim (3-v)A+B_1+B_2+B_3$. A similar statement holds in other parts. We have
$S_{S}(A)=\frac{5}{3}$.
Thus, $\delta_P(S)\le \frac{2}{5/3}=\frac{6}{5}$. Moreover,
$$h(v)=\begin{cases}
    \frac{v^2}{2}\text{ if }v\in[0,2],\\
    2 (3 - v) (2 v - 3)\text{ if }v\in[2,3].
    \end{cases}$$
So
$S(W_{\bullet,\bullet}^{A};O)\le 1$.
We get that $\delta_P(S)\ge 1$.
\end{proof}
\begin{lemma}\label{deg6-1_1_2_point}
Suppose $P$ belongs to a $(-1)$-curve $A$ and there exist $(-1)$-curves and $(-2)$-curves   which form the following dual graph:
\begin{figure}[h]
    \centering
   \includegraphics[width=11cm]{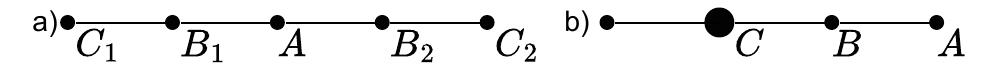}
    \caption{Dual graph: $(-K_S)^2=6$ and $\delta_P(S)=1$}
\end{figure}
   Then $\tau(A)=2$ and the Zariski Decomposition of the divisor $-K_S-vA$ is given by:
{\allowdisplaybreaks \begin{align*}
&{\text{\bf a). }}&P(v)=\begin{cases}
-K_{S}-vA\text{ if }v\in[0,1],\\
-K_{S}-vA-(v-1)(B_1+B_2)\text{ if }v\in[1,2].
\end{cases}
\\&&N(v)=\begin{cases}
0\text{ if }v\in[0,1],\\
(v-1)(B_1+B_2)\text{ if }v\in[1,2].
\end{cases}\\
&{\text{\bf b). }}&P(v)=\begin{cases}
-K_{S}-vA\text{ if }v\in[0,1],\\
-K_{S}-vA-(v-1)(2B+C)\text{ if }v\in[1,2].
\end{cases}
\\&&N(v)=\begin{cases}
0\text{ if }v\in[0,1],\\
(v-1)(2B+C)\text{ if }v\in[1,2].
\end{cases}
\end{align*}}
Moreover, 
$$(P(v))^2=\begin{cases}
6-2v-v^2\text{ if }v\in[0,1],\\
(2-v)(4-v)\text{ if }v\in[1,2].
\end{cases} 
P(v)\cdot A=\begin{cases}
v+1\text{ if }v\in[0,1],\\
3-v\text{ if }v\in[1,2].
\end{cases}$$
In this case  $\delta_P(S)=1\text{ if }P\in A\backslash B$.
\end{lemma}

\begin{proof}
The Zariski Decomposition in part a). follows from
$-K_S-vA\sim_{\DR} \big(2-v\big)A+\frac{1}{2}\big(3B_1+C_1+3B_2+C_2\big)$.
A similar statement holds in other parts. We have
$S_{S}(A)=1$.
Thus, $\delta_P(S)\le 1$ for $P\in A$. Moreover, for $P\in A\backslash B$:
$$h(v)\le \begin{cases}
\frac{(v+1)^2}{2}\text{ if }v\in[0,1],\\
\frac{(3 - v) (v + 1)}{2}\text{ if }v\in[1,2].
\end{cases}$$
So
$S(W_{\bullet,\bullet}^{A};P)\le 1$.
We get that $\delta_P(S)=1$ for $P\in A\backslash B$.
\end{proof}
\begin{lemma}\label{deg6-910_2_points}
Suppose $P$ belongs to a $(-1)$-curve $A$ and there exist $(-1)$-curves and $(-2)$-curves   which form the following dual graph:
\begin{figure}[h]
    \centering
   \includegraphics[width=4cm]{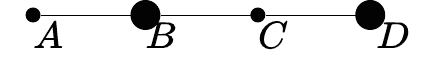}
    \caption{Dual graph: $(-K_S)^2=6$ and $\delta_P(S)=\frac{9}{10}$}
\end{figure}
 \par Then $\tau(A)=2$ and the Zariski Decomposition of the divisor $-K_S-vA$ is given by:
$$P(v)=-K_S-vA-\frac{v}{2}B \text{ and } N(v)=\frac{v}{2}B\text{ if }v\in[0,2].
$$
Moreover, 
$$(P(v))^2=\frac{(2-v)(6+v)}{2} \text{ and } P(v)\cdot A=1+\frac{v}{2}\text{ if }v\in[0,2].$$
In this case  $\delta_P(S)=\frac{9}{10}\text{ if }P\in A\backslash B$.
\end{lemma}
\begin{proof}
The Zariski Decomposition follows from
$-K_S-vA\sim_{\DR} \big(2-v\big)A+3B+4C+2D$. We have
$S_S(A)=\frac{10}{9}$.
Thus, $\delta_P(S)\le \frac{9}{10}$ for $P\in A$.  Note that if $P\in A\backslash B$ we have
$h(v)=\frac{(2+v)^2}{8}\text{ if }v\in[0,2]$.
So $S(W_{\bullet,\bullet}^{A};P)\le\frac{7}{9}\le \frac{10}{9}$.
Thus, $\delta_P(S)=\frac{9}{10}$ if $P\in A\backslash B$.
\end{proof}
\begin{lemma}\label{deg6-911_1_2_points}
Suppose $P$ belongs to a $(-2)$-curve $A$ and there exist $(-1)$-curves and $(-2)$-curves   which form the following dual graph and no other  $(-1)$-curves and $(-2)$-curves intersect $A$:
\begin{figure}[h]
    \centering
   \includegraphics[width=9cm]{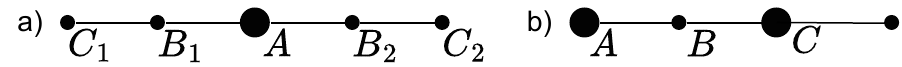}
    \caption{Dual graph: $(-K_S)^2=6$ and $\delta_P(S)=\frac{9}{11}$ with $\tau(A)=2$ }
\end{figure}
 \par  Then $\tau(A)=2$ and the Zariski Decomposition of the divisor $-K_S-vA$ is given by:
{\allowdisplaybreaks\begin{align*}
&{\text{\bf a). }}&P(v)=\begin{cases}
-K_S-vA\text{ if }v\in[0,1],\\
-K_S-vA-(v-1)(B_1+B_2)\text{ if }v\in[1,2].
\end{cases}
\\&&N(v)=\begin{cases}
0\text{ if }v\in[0,1],\\
(v-1)(B_1+B_2)\text{ if }v\in[1,2].
\end{cases}\\
&{\text{\bf b). }}&P(v)=\begin{cases}
-K_S-vA\text{ if }v\in[0,1],\\
-K_S-vA-(v-1)(2B_1+C)\text{ if }v\in[1,2].
\end{cases}
\\&&N(v)=\begin{cases}
0\text{ if }v\in[0,1],\\
(v-1)(2B_1+C)\text{ if }v\in[1,2].
\end{cases}
\end{align*}}
Moreover, 
$$(P(v))^2=\begin{cases}
6-2v^2\text{ if }v\in[0,1],\\
8-4v\text{ if }v\in[1,2].
\end{cases} 
P(v)\cdot A=\begin{cases}
2v\text{ if }v\in[0,1],\\
2\text{ if }v\in[1,2].
\end{cases}$$
In this case  $\delta_P(S)=\frac{9}{11}\text{ if }P\in A\backslash B$.
\end{lemma}
\begin{proof}
The Zariski Decomposition in part a). follows from
$-K_S-vA\sim_{\DR} \big(2-v\big)A+2B_1+C_1+2B_2+C_2$.
A similar statement holds in other parts. We have
$S_{S}(A)=\frac{11}{9}$.
Thus, $\delta_P(S)\le \frac{9}{11}$ for $P\in A\backslash B$. Moreover,
$$h(v)=\begin{cases}
2v^2\text{ if }v\in[0,1],\\
2v\text{ if }v\in[1,2].
\end{cases}$$
So
$S(W_{\bullet,\bullet}^{A};P)\le \frac{11}{9}$.
We get that $\delta_P(S)=\frac{9}{11}$ for $P\in A\backslash B$. 
\end{proof}
\begin{lemma}\label{deg6-911_1_2_3_points}
Suppose $P$ belongs to a $(-1)$-curve $A$ and there exist $(-1)$-curves and $(-2)$-curves   which form the following dual graph:
\begin{figure}[h]
    \centering
   \includegraphics[width=4.5cm]{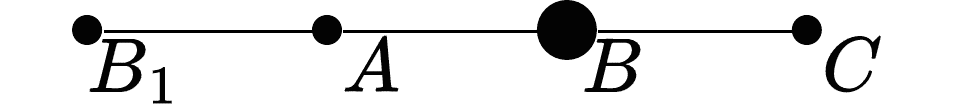}
    \caption{Dual graph: $(-K_S)^2=6$ and $\delta_P(S)=\frac{9}{11}$ with $\tau(A)=3$}
\end{figure}
\par   Then $\tau(A)=3$ and the Zariski Decomposition of the divisor $-K_S-vA$ is given by:
{\allowdisplaybreaks\begin{align*}
&&P(v)=\begin{cases}
-K_S-vA-\frac{v}{2}B\text{ if }v\in[0,1],\\
-K_S-vA-\frac{v}{2}B-(v-1)B_1\text{ if }v\in[1,2],\\
-K_S-vA-(v-1)(B+B_1)-(v-2)C\text{ if }v\in[2,3],
\end{cases}\\&&N(v)=\begin{cases}
\frac{v}{2}B\text{ if }v\in[0,1],\\
\frac{v}{2}B+(v-1)B_1\text{ if }v\in[1,2],\\
(v-1)(B+B_1)+(v-2)C\text{ if }v\in[2,3].
\end{cases}
\end{align*}}
Moreover, 
$$(P(v))^2=\begin{cases}
\frac{(2-v)(6+v)}{2}\text{ if }v\in[0,1],\\
\frac{v^2}{2}-4v+7\text{ if }v\in[1,2],\\
(3-v)^2\text{ if }v\in[2,3].\\
\end{cases} 
P(v)\cdot A=\begin{cases}
1+\frac{v}{2}\text{ if }v\in[0,1],\\
2-\frac{v}{2}\text{ if }v\in[1,2],\\
3-v\text{ if }v\in[2,3].
\end{cases}$$
In this case  $\delta_P(S)=\frac{9}{11}\text{ if }P\in A\backslash B$.
\end{lemma}
\begin{proof}
The Zariski Decomposition follows from
$-K_S-vA\sim_{\DR} \big(3-v\big)A+2B_1+2B+C$. We have
$S_{S}(A)=\frac{11}{9}$.
Thus, $\delta_P(S)\le \frac{9}{11}$ for $P\in A$. Moreover, if $P\in A\backslash B$:
$$h(v)\le \begin{cases}
 \frac{(v + 2)^2}{8} \text{ if }v\in[0,1],\\
\frac{3 (4 - v) v}{8}\text{ if }v\in[1,2],\\
 \frac{ (3 - v) (v + 1)}{2}\text{ if }v\in[2,3].
\end{cases}$$
So
$S(W_{\bullet,\bullet}^{A};P)\le 1\le \frac{11}{9}$.
We get that $\delta_P(S)=\frac{9}{11}$ for $P\in A\backslash B$. 
\end{proof}
\begin{lemma}\label{deg6-45_32_3_points}
Suppose $P$ belongs to a $(-1)$-curve $A$ and there exist $(-1)$-curves and $(-2)$-curves   which form the following dual graph:
\begin{figure}[h]
    \centering
   \includegraphics[width=3.7cm]{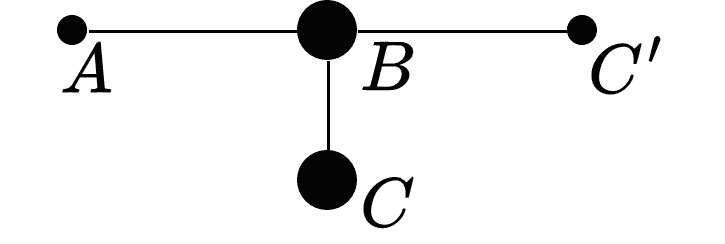}
    \caption{Dual graph: $(-K_S)^2=6$ and $\delta_P(S)=\frac{4}{5}$ }
\end{figure}
 \par     Then $\tau(A)=3$ and the Zariski Decomposition of the divisor $-K_S-vA$ is given by:
\begin{align*}
&&P(v)=\begin{cases}
-K_S-vA-\frac{v}{3}(2B+C)\text{ if }v\in\big[0,\frac{3}{2}\big],\\
-K_S-vA-(v-1)(2B+C)-(2v-3)C'\text{ if }v\in\big[\frac{3}{2},3\big].
\end{cases}\\&& N(v)=\begin{cases}
\frac{v}{3}(2B+C)\text{ if }v\in\big[0,\frac{3}{2}\big],\\
(v-1)(2B+C)+(2v-3)C'\text{ if }v\in\big[\frac{3}{2},3\big].
\end{cases}
\end{align*}
Moreover, 
$$(P(v))^2=\begin{cases}
6-2v-\frac{v^2}{3}\text{ if }v\in\big[0,\frac{3}{2}\big],\\
(3-v)^2\text{ if }v\in\big[\frac{3}{2},3\big].
\end{cases} 
P(v)\cdot A=\begin{cases}
1+\frac{v}{3}\text{ if }v\in\big[0,\frac{3}{2}\big],\\
3-v\text{ if }v\in\big[\frac{3}{2},3\big].
\end{cases}$$
In this case  $\delta_P(S)=\frac{4}{5}\text{ if }P\in A\backslash B$.
\end{lemma}
\begin{proof}
The Zariski Decomposition follows from
$-K_S-vA\sim_{\DR} \big(3-v\big)A+4B+2C+3C'$.
We have
$S_{S}(A)=\frac{5}{4}$.
Thus, $\delta_P(S)\le \frac{4}{5}$ for $P\in A$. Moreover, if $P\in A\backslash B$:
$$h(v)=\begin{cases}
\frac{(3+v)^2}{18}\text{ if }v\in\big[0,\frac{3}{2}\big],\\
\frac{(3-v)^2}{2}\text{ if }v\in\big[\frac{3}{2},3\big].
\end{cases}$$
So $S(W_{\bullet,\bullet}^{A};P)\le\frac{7}{12}\le \frac{5}{4}$. We get that $\delta_P(S)=\frac{4}{5}$ for $P\in A\backslash B$.
\end{proof}
\begin{lemma}\label{deg6-34_1_3_points}
Suppose $P$ belongs to a $(-2)$-curve $A$ and there exist $(-1)$-curves and $(-2)$-curves   which form the following dual graph:
\begin{figure}[h]
    \centering
   \includegraphics[width=12cm]{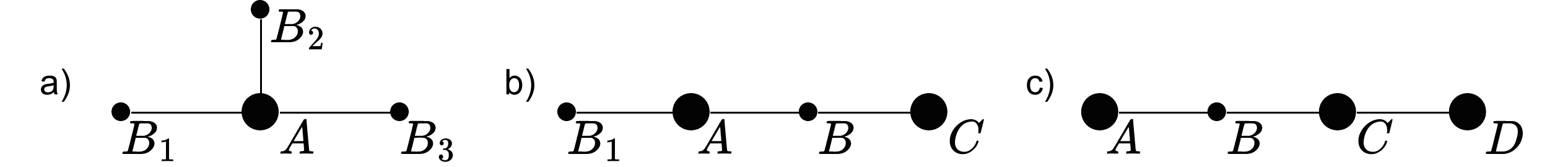}
    \caption{Dual graph: $(-K_S)^2=6$ and $\delta_P(S)=\frac{3}{4}$ with $\tau(A)=3$}
\end{figure}
  \par Then $\tau(A)=3$ and the Zariski Decomposition of the divisor $-K_S-vA$ is given by:
{\allowdisplaybreaks\begin{align*}
&{\text{\bf a). }}&P(v)=\begin{cases}
-K_S-vA\text{ if }v\in[0,1],\\
-K_S-vA-(v-1)(B_1+B_2+B_3)\text{ if }v\in[1,3].
\end{cases}
\\&&N(v)=\begin{cases}0\text{ if }v\in[0,1],\\
(v-1)(B_1+B_2+B_3)\text{ if }v\in[1,3].
\end{cases}\\
&{\text{\bf b). }}&P(v)=\begin{cases}
-K_S-vA\text{ if }v\in[0,1],\\
-K_S-vA-(v-1)(B_1+2B+C)\text{ if }v\in[1,3].
\end{cases}
\\&&N(v)=\begin{cases}0\text{ if }v\in[0,1],\\
(v-1)(B_1+2B+C)\text{ if }v\in[1,3].
\end{cases}\\
&{\text{\bf c). }}&P(v)=\begin{cases}
-K_S-vA\text{ if }v\in[0,1],\\
-K_S-vA-(v-1)(3B+2C+D)\text{ if }v\in[1,3].
\end{cases}
\\&&N(v)=\begin{cases}0\text{ if }v\in[0,1],\\
(v-1)(3B+2C+D)\text{ if }v\in[1,3].
\end{cases}
\end{align*}}
Moreover, 
$$(P(v))^2=\begin{cases}2(3-v^2)\text{ if }v\in[0,1],\\
(3-v)^2\text{ if }v\in[1,3].
\end{cases} P(v)\cdot A=\begin{cases}2v\text{ if }v\in[0,1],\\
3-v\text{ if }v\in[1,3].
\end{cases}$$
In this case  $\delta_P(S)=\frac{3}{4}\text{ if }P\in A\backslash B$.
\end{lemma}
\begin{proof}
The Zariski Decomposition in part a). follows from
$-K_S-vA\sim_{\DR} \big(3-v\big)A+2B_1+2B_2+2B_3$.
A similar statement holds in other parts. We have
$S_S(A)=\frac{4}{3}$.
Thus, $\delta_P(S)\le \frac{3}{4}$ for $P\in A$. Note that  if $P\in A\backslash B$ we have:
$$h(v)\le \begin{cases} 
2v^2\text{ if }v\in[0,1],\\
\frac{ (3 - v) (v + 1)}{2}\text{ if }v\in[1,3].
\end{cases}$$
So 
$S(W_{\bullet,\bullet}^{A};P)\le\frac{10}{9}\le \frac{4}{3}$.
Thus, $\delta_P(S)=\frac{3}{4}$ if $P\in A\backslash B$.
\end{proof}
 \begin{lemma}\label{deg6-34_2_points}
Suppose $P$ belongs to a $(-1)$-curve $A$ and there exist $(-1)$-curves and $(-2)$-curves   which form the following dual graph:
\begin{figure}[h]
    \centering
   \includegraphics[width=7.2cm]{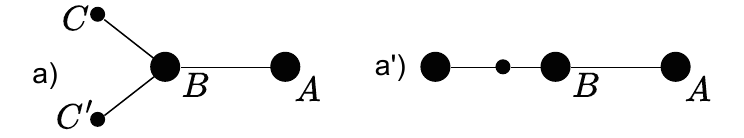}
    \caption{Dual graph: $(-K_S)^2=6$ and $\delta_P(S)=\frac{3}{4}$ with $\tau(A)=2$}
\end{figure}
 \par   Then $\tau(A)=2$ and the Zariski Decomposition of the divisor $-K_S-vA$ is given by:
$$P(v)=-K_S-vA-\frac{v}{2}B \text{ and } N(v)=\frac{v}{2}B\text{ if }v\in[0,2].
$$
Moreover, 
$$(P(v))^2=\frac{3(2-v)(6+v)}{2} \text{ and } P(v)\cdot A= \frac{3v}{2} \text{ if }v\in[0,2].$$
In this case  $\delta_P(S)=\frac{3}{4}\text{ if }P\in A\backslash B$.
\end{lemma}
\begin{proof}
The Zariski Decomposition in part a). follows from
$-K_S-vA\sim_{\DR} \big(2-v\big)A+4B+3C+3C'$.
We have $S_S(A)=\frac{4}{3}$.
Thus, $\delta_P(S)\le \frac{3}{4}$ for $P\in A$.  Note that if $P\in A\backslash B$ we have
$h(v)=\frac{9v^2}{8}\text{ if }v\in[0,2]$.
 So $S(W_{\bullet,\bullet}^{A};P)\le\frac{4}{3}$.
Thus, $\delta_P(S)=\frac{3}{4}$ if $P\in  A\backslash B$.
\end{proof}
\begin{lemma}\label{deg6-914_2_4_points}
Suppose $P$ belongs to a $(-1)$-curve $A$ and there exist $(-1)$-curves and $(-2)$-curves   which form the following dual graph:
\begin{figure}[h]
    \centering
   \includegraphics[width=4.5cm]{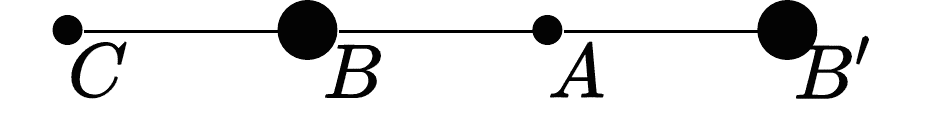}
    \caption{Dual graph: $(-K_S)^2=6$ and $\delta_P(S)=\frac{9}{14}$ }
\end{figure}
\par   Then $\tau(A)=4$ and the Zariski Decomposition of the divisor $-K_S-vA$ is given by:
{\allowdisplaybreaks\begin{align*}
&&P(v)=\begin{cases}
-K_S-vA-\frac{v}{2}(B+B')\text{ if }v\in[0,2],\\
-K_S-vA-\frac{v}{2}B'-(v-1)B-(v-2)C\text{ if }v\in[2,4].
\end{cases}
\\&&N(v)=\begin{cases}
\frac{v}{2}(B+B')\text{ if }v\in[0,2],\\
\frac{v}{2}B'+(v-1)B+(v-2)C\text{ if }v\in[2,4].
\end{cases}
\end{align*}}
Moreover, 
$$(P(v))^2=\begin{cases}
2(3-v)\text{ if }v\in[0,2],\\
\frac{(v-4)^2}{2}\text{ if }v\in[2,4].
\end{cases} 
P(v)\cdot A=\begin{cases}
1\text{ if }v\in[0,2],\\
2-\frac{v}{2}\text{ if }v\in[2,4].
\end{cases}$$
In this case  $\delta_P(S)=\frac{9}{14}\text{ if }P\in A$.
\end{lemma}
\begin{proof}
The Zariski Decomposition follows from
$-K_S-vA\sim_{\DR} \big(4-v\big)A+3B+2C+2B'$. We have
$S_{S}(A)=\frac{14}{9}$.
Thus, $\delta_P(S)\le \frac{9}{14}$ for $P\in A$. 
Moreover if $P\in A\cap B'$ or $P\in A\backslash B'$:
$$h(v)=\begin{cases}
 \frac{v+1}{2} \text{ if }v\in[0,2],\\
 \frac{(4-v)(4+v)}{8} \text{ if }v\in[2,4].
\end{cases}\text{ or }h(v)\le\begin{cases}
 \frac{v+1}{2} \text{ if }v\in[0,2],\\
 \frac{3 (4 - v) v}{8} \text{ if }v\in[2,4].
\end{cases}$$
So
$S(W_{\bullet,\bullet}^{A};P)\le \frac{11}{9}\le \frac{14}{9}$
or 
$S(W_{\bullet,\bullet}^{A};P)\le \frac{4}{3}\le \frac{14}{9}$.
We get that $\delta_P(S)=\frac{9}{14}$ for $P\in A$.
\end{proof}
\begin{lemma}\label{deg6-35_1_4_points}
Suppose $P$ belongs to a $(-2)$-curve $A$ and there exist $(-1)$-curves and $(-2)$-curves   which form the following dual graph:
\begin{figure}[h]
    \centering
   \includegraphics[width=9.2cm]{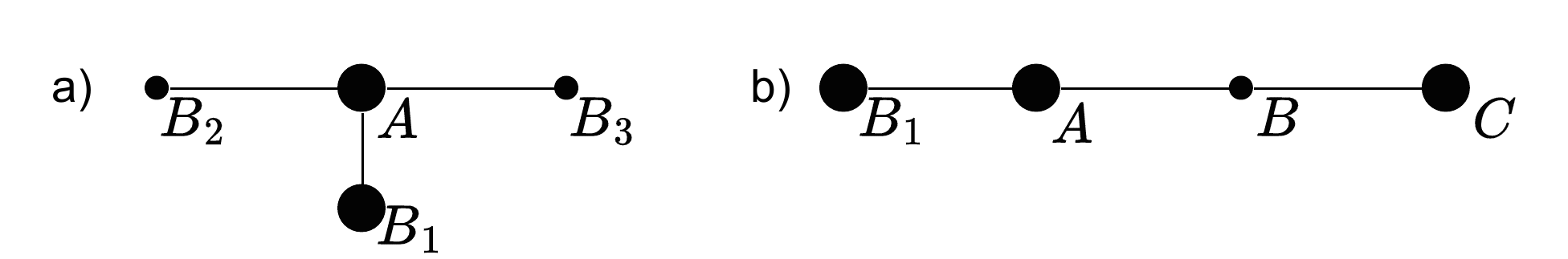}
    \caption{Dual graph: $(-K_S)^2=6$ and $\delta_P(S)=\frac{3}{5}$ }
\end{figure}
 \par  Then $\tau(A)=4$ and the Zariski Decomposition of the divisor $-K_S-vA$ is given by:
{\allowdisplaybreaks\begin{align*}
&{\text{\bf a). }}&P(v)=\begin{cases}
-K_S-vA-\frac{v}{2}B_1\text{ if }v\in[0,1],\\
-K_S-vA-\frac{v}{2}B_1-(v-1)(B_2+B_3)\text{ if }v\in[1,4].
\end{cases} \\&&N(v)=\begin{cases}
\frac{v}{2}B_1\text{ if }v\in[0,1],\\
\frac{v}{2}B_1-(v-1)(B_2+B_3)\text{ if }v\in[1,4].
\end{cases}\\
&{\text{\bf b). }}&P(v)=\begin{cases}
-K_S-vA-\frac{v}{2}B_1\text{ if }v\in[0,1],\\
-K_S-vA-\frac{v}{2}B_1-(v-1)(2B+C)\text{ if }v\in[1,4].
\end{cases} \\&&N(v)=\begin{cases}
\frac{v}{2}B_1\text{ if }v\in[0,1],\\
\frac{v}{2}B_1-(v-1)(2B+C)\text{ if }v\in[1,4].
\end{cases}
\end{align*}}
Moreover, 
$$(P(v))^2=\begin{cases} 
\frac{3(2-v)(2+v)}{2} \text{ if }v\in[0,1],\\
\frac{(4-v)^2}{2} \text{ if }v\in[1,4].
\end{cases} P(v)\cdot A=\begin{cases} \frac{3v}{2} \text{ if }v\in[0,1],\\
2-\frac{v}{2}\text{ if }v\in[1,4].
\end{cases}$$
In this case  $\delta_P(S)=\frac{3}{5}\text{ if }P\in A\backslash B$.
\end{lemma}
\begin{proof}
The Zariski Decomposition follows from
$-K_S-vA\sim_{\DR} \big(4-v\big)A+2B_1+3B_2+3B_3$. A similar statement holds in other parts.  We have
$S_S(A)=\frac{5}{3}$.
Thus, $\delta_P(S)\le \frac{3}{5}$ for $P\in E_2$.  
Note that if $P\in A\backslash B_1$ or if $P\in A\cap B_1$: then:
$$h(v)\le\begin{cases}
\frac{9v^2}{8}\text{ if }v\in[0,1],\\
\frac{3 (4 - v) v}{8} \text{ if }v\in[1,4].
    \end{cases}\text{ or }h(v)\le\begin{cases}
     \frac{15v^2}{8} \text{ if }v\in[0,1],\\
     \frac{(4-v)(4+v)}{8} \text{ if }v\in[1,4].
    \end{cases}$$
So 
$S(W_{\bullet,\bullet}^{A};P)\le\frac{5}{4}\le \frac{5}{3}$ or 
$S(W_{\bullet,\bullet}^{A};P)\le\frac{4}{3}\le \frac{5}{3}$. Thus, $\delta_P(S)=\frac{3}{5}$ if $P\in A\backslash B$.
\end{proof}
\begin{lemma}\label{deg6-12_6_points}
Suppose $P$ belongs to a $(-1)$-curve $A$ and there exist $(-1)$-curves and $(-2)$-curves   which form the following dual graph:
\begin{figure}[h]
    \centering
   \includegraphics[width=5cm]{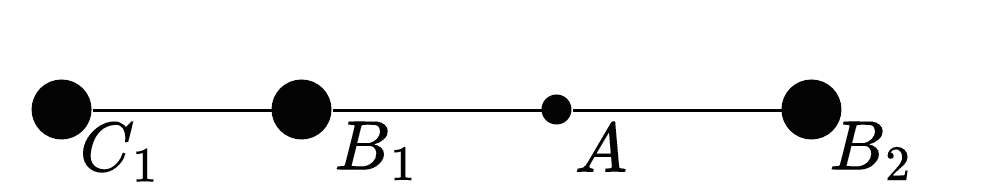}
    \caption{Dual graph: $(-K_S)^2=6$ and $\delta_P(S)=\frac{1}{2}$ }
\end{figure}
 \par    Then $\tau(A)=6$ and the Zariski Decomposition of the divisor $-K_S-vA$ is given by:
$$P(v)=-K_S-vA-\frac{v}{3}(2B_1+C_1)-\frac{v}{2}B_2\text{ and } N(v)=\frac{v}{3}(2B_1+C_1)+\frac{v}{2}B_2\text{ if }v\in[0,6].
$$
Moreover, 
$$(P(v))^2=\frac{(6-v)^2}{2} P(v)\cdot A=1-\frac{v}{6}\text{ if }v\in[0,6].$$
In this case  $\delta_P(S)=\frac{1}{2}\text{ if }P\in A$.
\end{lemma}
\begin{proof}
The Zariski Decomposition follows from
$-K_S-vA\sim_{\DR} \big(6-v\big)A+4B_1+2C_1+3B_2$. We have $S_S(A)=2$.
Thus, $\delta_P(S)\le \frac{1}{2}$ for $P\in A$.  Note that 
$h(v)\le \frac{(6 - v) (7 v + 6)}{72}\text{ if }v\in[0,6].$
So
$S(W_{\bullet,\bullet}^{A};P)\le\frac{5}{3}\le 2 $.
Thus, $\delta_P(S)=\frac{1}{2}$ if $P\in  A$.
\end{proof}


\subsection{Finding $\delta$-invariants for degree $6$}

Let $X$ be a singular del Pezzo surface of degree $6$ with and $S$ be a minimal resolution of $X$. Then there are several possible cases:
\begin{itemize}
    \item[I.] $X$ has an $\DA_1$ singularity and contains $3$ lines. In this case, we let  $E$ be the exceptional divisor,  $L_{i}$ for $i\in\{1,2,3\}$ are the lines on $S$,  
    \item[II.] $X$ has an $\DA_1$ singularity and contains $4$ lines. In this case, we let $E$ be the exceptional divisor, 
 $L_{i}$ and $L_i'$ for $i\in\{1,2\}$ be the lines on $S$,  
    \item[III.] $X$ has  two $\DA_1$ singularities and contains $2$ lines. In this case, we let $E_i$ for $i\in\{1,2\}$ be the exceptional divisors,  $L_{1}$ and $L_{12}$ are the lines on $S$,  
    \item[IV.] $X$ has an $\DA_2$ singularity and contains $2$ lines. In this case, we let $E_i$ for $i\in\{1,2\}$ be the exceptional divisors,  $L_{1}$ and $L_1'$ be the lines on $S$,  
    \item[V.] $X$ has $\DA_2$ and $\DA_1$ singularities and contains $1$ line. In this case, we let $E_i$ for $i\in\{1,2,3\}$ be the exceptional divisors, $L_{13}$ be the line on $S$.  
\end{itemize}
such that the dual graph of the $(-1)$-curves  and $(-2)$-curves on $S$ is given the picture below. 
\begin{figure}[h!]
    \centering
   \includegraphics[width=15cm]{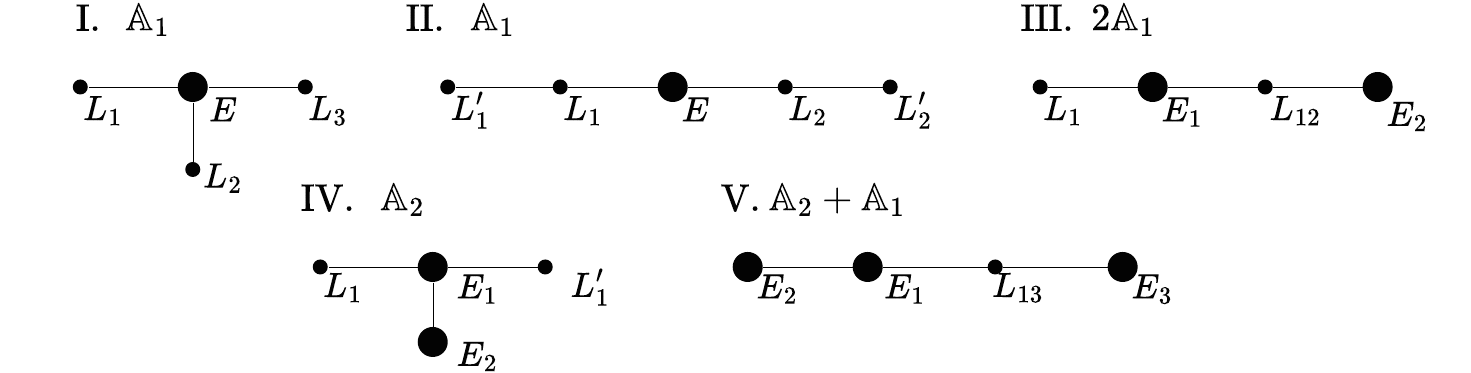}
    \caption{Du Val del Pezzo surfaces with $(-K_S)^2=6$}
\end{figure}
 
Then
\begin{itemize}
    \item[I.] $\delta(X)=\frac{3}{4}$ since depending on the position of point $P\in S$ we have 
\begin{table}[h!]
\hspace*{1cm}\begin{tabular}{ | c || c | c | c |}
 \hline 
 $P$ & $E$ & $(L_1\cup L_2\cup L_3)\backslash E$ & o/w \\\hline
$\delta_P(S)$ & $\frac{3}{4}$ & $\frac{9}{10}$ & $\ge 1$\\
\hline
    \end{tabular}
    \caption{Local $\delta$-invariants: $(-K_S)^2=6$ and  $\DA_1$ singularity, 3 lines}
\end{table}
    \item[II.] $\delta(X)=\frac{9}{11}$ since since depending on the position of point $P\in S$ we have 
\begin{table}[h!]
\hspace*{1cm}\begin{tabular}{ | c || c | c | c |}
 \hline 
 $P$ & $ E\cup L_1\cup L_2$ & $(L_1'\cup L_2')\backslash (L_1\cup L_2)$ & o/w \\\hline
$\delta_P(S)$ & $\frac{9}{10}$ & $1$ & $\ge 1$\\\hline
    \end{tabular}
    \caption{Local $\delta$-invariants: $(-K_S)^2=6$ and $\DA_1$ singularity, 4 lines}
\end{table}
\item[III.]  $\delta(X)=\frac{9}{14}$ since since depending on the position of point $P\in S$ we have 
\begin{table}[h!]
\hspace*{1cm}
\begin{tabular}{ | c || c | c | c | c |  c |}
 \hline 
 $P$ & $L_{12}$ & $E_1\backslash L_{12}$ & $E_2\backslash L_{12} $ & $L_1\backslash  E_1$ & o/w \\\hline
$\delta_P(S)$ & $\frac{9}{14}$ & $\frac{3}{4}$ & $\frac{9}{11}$ & $\frac{9}{10}$ &  $\ge 1$\\\hline
    \end{tabular}
   \caption{Local $\delta$-invariants: $(-K_S)^2=6$ and  $2\DA_1$ singularities}
\end{table}
    \item[IV.] $\delta(X)=\frac{3}{5}$ since since depending on the position of point $P\in S$ we have 
\begin{table}[h!]
\hspace*{1cm} \begin{tabular}{ | c || c | c | c |c |}
 \hline 
 $P$ & $ E_1$ & $ E_2\backslash E_1$ & $(L_{1}\cup L_1')\backslash E_1$ & o/w \\\hline
$\delta_P(S)$ & $\frac{3}{5}$ & $\frac{3}{4}$ & $\frac{4}{5}$  & $\ge 1$\\\hline
    \end{tabular}
 \caption{Local $\delta$-invariants: $(-K_S)^2=6$ and  $\DA_2$ singularity}
\end{table}
 \item[V.] $\delta(X)=\frac{1}{2}$ since since depending on the position of point $P\in S$ we have 
\begin{table}[h!]
\hspace*{1cm}
\begin{tabular}{ | c || c | c | c |c |}
 \hline 
 $P$ & $L_{13}$ & $E_1\backslash L_{13}$ & $(E_2\cup E_3)\backslash (E_1\backslash L_{13})$ & o/w \\\hline
$\delta_P(S)$ & $\frac{1}{2}$ & $\frac{3}{5}$ & $\frac{3}{4}$  & $\ge 1$\\\hline
    \end{tabular}
\caption{Local $\delta$-invariants: $(-K_S)^2=6$ and  $\DA_2\DA_1$ singularities}
\end{table}
\end{itemize}

\newpage
\begin{proof}
We prove each case separately using lemmas from the previous section.
    \begin{itemize}
        \item[I.]  If $P\in E$ the assertion follows from Lemma \ref{deg6-34_1_3_points}.
 If $P\in  (L_1\cup L_2\cup L_3)\backslash E$ the assertion follows from Lemma \ref{deg6-910_2_points}.
 If $P$ is a general point the assertion follows from Lemma \ref{deg6-generalpoint}.  
    \item[II.] If $P\in E$ the assertion follows from Lemma \ref{deg6-911_1_2_points} [a).].
 If $P\in  (L_1\cup L_2)\backslash E$ the assertion follows from Lemma \ref{deg6-911_1_2_3_points}.
 If $P\in  (L_1'\cup L_2')\backslash (L_1\cup L_2)$ the assertion follows from Lemma \ref{deg6-1_1_2_point} [b).].
 If $P$ is a general point the assertion follows from Lemma \ref{deg6-generalpoint}.
  \item[III.]  If $P\in L_{12}$ the assertion follows from Lemma \ref{deg6-914_2_4_points}.
 If $P\in  E_1\backslash L_{12}$ the assertion follows from Lemma \ref{deg6-34_1_3_points} [b).].
 If $P\in  E_2\backslash L_{12}$ the assertion follows from Lemma \ref{deg6-911_1_2_points} [b).].
 If $P\in  L_1\backslash  E_1$ the assertion follows from Lemma \ref{deg6-910_2_points}.
 If $P$ is a general point the assertion follows from Lemma \ref{deg6-generalpoint}. 
 \item[IV.]  If $P\in E_1$ the assertion follows from Lemma \ref{deg6-35_1_4_points} [a).]. 
 If $P\in  E_2\backslash E_1$ the assertion follows from Lemma \ref{deg6-34_2_points}.
 If $P\in  (L_{1}\cup L_1')\backslash E_1$ the assertion follows from Lemma \ref{deg6-45_32_3_points}.
 If $P$ is a general point the assertion follows from Lemma \ref{deg6-generalpoint}. 
  \item[V.] If $P\in L_{12}$ the assertion follows from Lemma \ref{deg6-12_6_points}.
 If $P\in E_1\backslash L_{13}$ the assertion follows from Lemma \ref{deg6-35_1_4_points} [b).].
 If $P\in  E_2\backslash E_1$ the assertion follows from Lemma \ref{deg6-34_2_points}.
 If $P\in E_2\backslash L_{13}$ the assertion follows from Lemma \ref{deg6-34_1_3_points}.
 If $P$ is a general point the assertion follows from Lemma \ref{deg6-generalpoint}. 
    \end{itemize}
\end{proof}

\section{Du Val Del Pezzo surfaces of degree $5$}
\noindent In \cite[Lemma 2.11]{Fano21} it was proven that $\delta(X)=\frac{15}{13}$ when $X$ is a smooth del Pezzo surface of degree $5$. In this chapter, we compute  $\delta$-invariants  of singular  Du Val del Pezzo surfaces of degree $5$. 
\begin{maintheorem*}
Let $X$ be a singular Du Val  del Pezzo surface of degree $5$. Then the  $\delta$-invariant of $X$ is uniquely determined by the degree of $X$, the number of lines on $X$, and the type of singularities on $X$ which is given in the following table:
\begin{table}[h]
    \centering\renewcommand{\arraystretch}{1.1}
    \begin{tabular}{ | c | c | c | c | }
   \hline
   $K_X^2$ & $\#$ lines & $\mathrm{Sing}(X)$ & $\delta(X)$\\
  \hline\hline
 $5$ & $7$ & $\DA_1$ & $\frac{15}{17}$\\
\hline
  $5$ & $5$ & $2\DA_1$ & $\frac{15}{19}$\\
\hline
$5$ & $4$ & $\DA_2$ & $\frac{5}{7}$\\
\hline
  \end{tabular}\renewcommand{\arraystretch}{1.1}
    \begin{tabular}{ | c | c | c | c | }
   \hline
   $K_X^2$ & $\#$ lines & $\mathrm{Sing}(X)$ & $\delta(X)$\\
  \hline\hline
  $5$ & $3$ & $\DA_2+\DA_1$ & $\frac{15}{23}$\\
\hline
$5$ & $2$ & $\DA_3$ & $\frac{5}{9}$\\
\hline 
 $5$ & $1$ & $\DA_4$ & $\frac{3}{7}$\\
\hline
  \end{tabular}
\end{table}
 \end{maintheorem*}
\subsection{General results for degree $5$}
Let $X$ be a del Pezzo surface of degree $5$ with at most Du Val singularities, $S$ be a minimal resolution of $X$ and $P$ is a point on $S$.  Then: 
\begin{lemma}\label{deg5-generalpoint}
Suppose is a general point on $S$. Consider the blowup $\sigma:\widetilde{S}\to S$ of $S$ at $P$ with the exceptional divisor $A$. There exist $(-1)$-curves and $(-2)$-curves   which form  one of the following dual graphs:
\begin{figure}[h!]
    \centering
   \includegraphics[width=12cm]{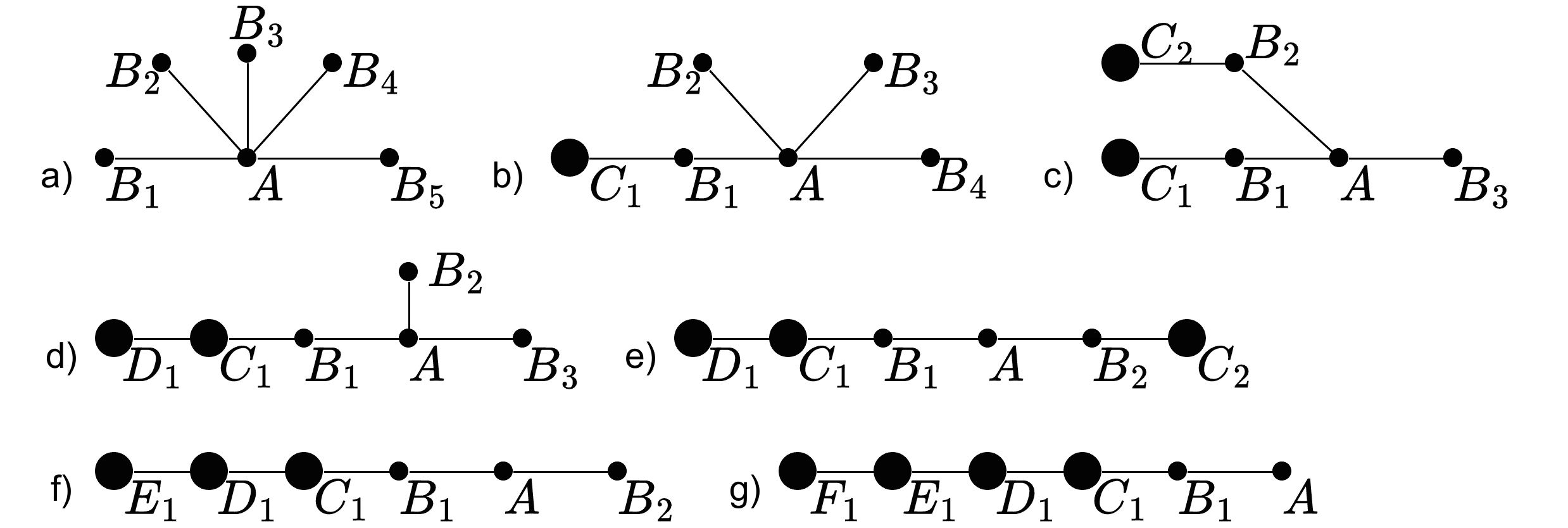}
    \caption{Dual graph: $(-K_S)^2=5$ for a general point}
\end{figure}
Then $\tau(A)=\frac{5}{2}$ and the Zariski decomposition of the divisor $\sigma^*(-K_S)-vE_P$ is given by:
{\allowdisplaybreaks \begin{align*}
&{\text{\bf a). }}&P(v)=\begin{cases}
\sigma^*(-K_{S})-vA\text{ if }v\in[0,2],\\
\sigma^*(-K_{S})-vA-(v-2)(B_1+B_2+B_3+B_4+B_5)\text{ if }v\in\big[2,\frac{5}{2}\big].
\end{cases}\\&&
N(v)=\begin{cases}
0\text{ if }v\in[0,2],\\
(v-2)(B_1+B_2+B_3+B_4+B_5)\text{ if }v\in\big[2,\frac{5}{2}\big].
\end{cases}\\
&{\text{\bf b). }}&P(v)=\begin{cases}
\sigma^*(-K_{S})-vA\text{ if }v\in[0,2],\\
\sigma^*(-K_{S})-vA-(v-2)(2B_1+C_1+B_2+B_3+B_4)\text{ if }v\in\big[2,\frac{5}{2}\big].
\end{cases}\\&&
N(v)=\begin{cases}
0\text{ if }v\in[0,2],\\
(v-2)(2B_1+C_1+B_2+B_3+B_4)\text{ if }v\in\big[2,\frac{5}{2}\big].
\end{cases}\\
&{\text{\bf c). }}&P(v)=\begin{cases}
\sigma^*(-K_{S})-vA\text{ if }v\in[0,2],\\
\sigma^*(-K_{S})-vA-(v-2)(2B_1+C_1+2B_2+C_2+B_3)\text{ if }v\in\big[2,\frac{5}{2}\big].
\end{cases}\\&&
N(v)=\begin{cases}
0\text{ if }v\in[0,2],\\
(v-2)(2B_1+C_1+2B_2+C_2+B_3)\text{ if }v\in\big[2,\frac{5}{2}\big].
\end{cases}\\
&{\text{\bf d). }}&P(v)=\begin{cases}
\sigma^*(-K_{S})-vA\text{ if }v\in[0,2],\\
\sigma^*(-K_{S})-vA-(v-2)(3B_1+2C_1+D_1+B_2+B_3)\text{ if }v\in\big[2,\frac{5}{2}\big].
\end{cases}\\&&
N(v)=\begin{cases}
0\text{ if }v\in[0,2],\\
(v-2)(3B_1+2C_1+D_1+B_2+B_3)\text{ if }v\in\big[2,\frac{5}{2}\big].
\end{cases}\\
&{\text{\bf e). }}&P(v)=\begin{cases}
\sigma^*(-K_{S})-vA\text{ if }v\in[0,2],\\
\sigma^*(-K_{S})-vA-(v-2)(3B_1+2C_1+D_1+2B_2+C_2)\text{ if }v\in\big[2,\frac{5}{2}\big].
\end{cases}\\&&
N(v)=\begin{cases}
0\text{ if }v\in[0,2],\\
(v-2)(3B_1+2C_1+D_1+2B_2+C_2)\text{ if }v\in\big[2,\frac{5}{2}\big].
\end{cases}\\
&{\text{\bf f). }}&P(v)=\begin{cases}
\sigma^*(-K_{S})-vA\text{ if }v\in[0,2],\\
\sigma^*(-K_{S})-vA-(v-2)(4B_1+3C_1+2D_1+E_1+B_2)\text{ if }v\in\big[2,\frac{5}{2}\big].
\end{cases}\\&&
N(v)=\begin{cases}
0\text{ if }v\in[0,2],\\
(v-2)(4B_1+3C_1+2D_1+E_1+B_2)\text{ if }v\in\big[2,\frac{5}{2}\big].
\end{cases}\\
&{\text{\bf g). }}&P(v)=\begin{cases}
\sigma^*(-K_{S})-vA\text{ if }v\in[0,2],\\
\sigma^*(-K_{S})-vA-(v-2)(5B_1+4C_1+3D_1+2E_1+F_1)\text{ if }v\in\big[2,\frac{5}{2}\big].
\end{cases}\\&&
N(v)=\begin{cases}
0\text{ if }v\in[0,2],\\
(v-2)(5B_1+4C_1+3D_1+2E_1+F_1)\text{ if }v\in\big[2,\frac{5}{2}\big].
\end{cases}
\end{align*}
}
Moreover, 
$$(P(v))^2=
\begin{cases}
5-v^2\text{ if }v\in[0,2],\\
(2v-5)^2\text{ if }v\in\big[2,\frac{5}{2}\big].
\end{cases}
P(v)\cdot A=
\begin{cases}
v\text{ if }v\in[0,2],\\
2(5-2v)\text{ if }v\in\big[2,\frac{5}{2}\big].
\end{cases}$$
In this case  $\delta_P(S)\ge\frac{6}{5}$.
\end{lemma}
\begin{proof}
The Zariski Decomposition in part a). follows from 
$$\sigma^*(-K_S)-vA\sim_{\DR} \Big(\frac{5}{2}-v\Big)A+\frac{1}{2}\Big(B_1+B_2+B_3+B_4+B_5\Big).$$ 
A similar statement holds in other parts. 
We have $S_{S}(A)=\frac{3}{2}$. Thus, $\delta_P(S)\le \frac{2}{3/2}=\frac{4}{3}$. Moreover,
$$h(v)\le\begin{cases}
    \frac{v^2}{2} \text{ if }v\in[0,2],\\
    2 (5 - 2 v) (3 v - 5)\text{ if }v\in\big[2,\frac{5}{2}\big].
    \end{cases}
    $$
So 
$S(W_{\bullet,\bullet}^{A};O)\le\frac{5}{6}$.
We get that $\delta_P(S)\ge\frac{6}{5}$.
\end{proof}
\begin{lemma}\label{deg5-1513_1_2_points}
Suppose $P$ belongs to a $(-1)$-curve $A$ and there exist $(-1)$-curves and $(-2)$-curves   which form the following dual graph:
\begin{figure}[h!]
    \centering
   \includegraphics[width=12cm]{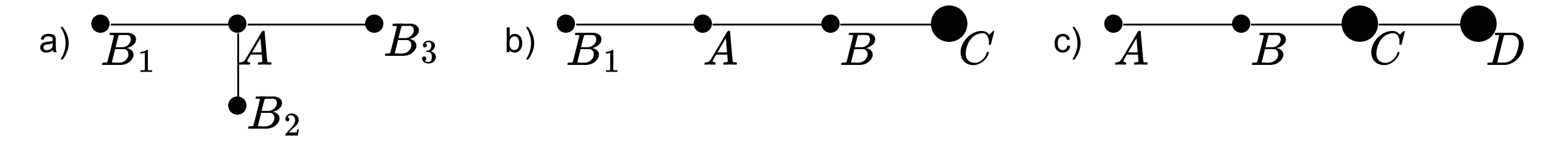}
    \caption{Dual graph: $(-K_S)^2=5$ and $\delta_P(S)=\frac{15}{13}$}
\end{figure}
 \par Then $\tau(A)=2$ and the Zariski Decomposition of the divisor $-K_S-vA$ is given by:
{\allowdisplaybreaks \begin{align*}
&{\text{\bf a). }}&P(v)=\begin{cases}
-K_{S}-vA\text{ if }v\in[0,1],\\
-K_{S}-vA-(v-1)(B_1+B_2+B_3)\text{ if }v\in[1,2].
\end{cases}
\\&& N(v)=\begin{cases}
0\text{ if }v\in[0,1],\\
(v-1)(B_1+B_2+B_3)\text{ if }v\in[1,2].
\end{cases}\\
&{\text{\bf b). }}&P(v)=\begin{cases}
-K_{S}-vA\text{ if }v\in[0,1],\\
-K_{S}-vA-(v-1)(B_1+2B+C)\text{ if }v\in[1,2].
\end{cases}
\\&& N(v)=\begin{cases}
0\text{ if }v\in[0,1],\\
(v-1)(B_1+2B+C)\text{ if }v\in[1,2].
\end{cases}\\
&{\text{\bf c). }}&P(v)=\begin{cases}
-K_{S}-vA\text{ if }v\in[0,1],\\
-K_{S}-vA-(v-1)(3B+2C+D)\text{ if }v\in[1,2].
\end{cases}
\\&& N(v)=\begin{cases}
0\text{ if }v\in[0,1],\\
(v-1)(3B+2C+D)\text{ if }v\in[1,2].
\end{cases}
\end{align*}}
Moreover, 
$$(P(v))^2=\begin{cases}
5-2v-v^2\text{ if }v\in[0,1],\\
2(2-v)^2\text{ if }v\in[1,2].
\end{cases}
P(v)\cdot A=\begin{cases}
v+1\text{ if }v\in[0,1],\\
2(2-v)\text{ if }v\in[1,2].
\end{cases}$$
In this case  $\delta_P(S)=\frac{15}{13}\text{ if }P\in A\backslash B$.
\end{lemma}

\begin{proof}
The Zariski Decomposition in part a). follows from
$-K_S-vA\sim_{\DR} (2-v)A+B_1+B_2+B_3$.
A similar statement holds in other parts. We have
$S_{S}(A)=\frac{15}{13}$.
Thus, $\delta_P(S)\le \frac{15}{13}$ for $P\in A$. Moreover for $P\in A\backslash B$,
$$h(v)\le\begin{cases}
\frac{(v+1)^2}{2}\text{ if }v\in[0,1],\\
2(2-v)\text{ if }v\in[1,2].
\end{cases}$$
So
$S(W_{\bullet,\bullet}^{A};P)\le\frac{15}{13}$.
We get that $\delta_P(S)=\frac{15}{13}$ for $P\in A\backslash B$.
\end{proof}
\begin{lemma}\label{deg5-1_1_2_points}
Suppose $P$ belongs to a $(-1)$-curve $A$ and there exist $(-1)$-curves and $(-2)$-curves   which form the following dual graph:
\begin{figure}[h!]
    \centering
   \includegraphics[width=3.7cm]{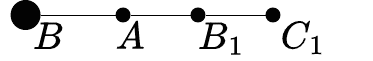}
    \caption{Dual graph: $(-K_S)^2=5$ and $\delta_P(S)=1$}
\end{figure}
\par Then $\tau(A)=2$ and the Zariski Decomposition of the divisor $-K_S-vA$ is given by:
{\allowdisplaybreaks \begin{align*}
P(v)=\begin{cases}
-K_S-vA-\frac{v}{2}B\text{ if }v\in[0,1],\\
-K_S-vA-\frac{v}{2}B-(v-1)B_1\text{ if }v\in[1,2].
\end{cases} N(v)=\begin{cases}
\frac{v}{2}B\text{ if }v\in[0,1],\\
\frac{v}{2}B+(v-1)B_1\text{ if }v\in[1,2].
\end{cases}
\end{align*}}
Moreover, 
$$(P(v))^2=\begin{cases}
5-2v-\frac{v^2}{2}\text{ if }v\in[0,1],\\
\frac{(2-v)(6-v)}{2}\text{ if }v\in[1,2].
\end{cases}P(v)\cdot A=
\begin{cases}
1+\frac{v}{2}\text{ if }v\in[0,1],\\
2-\frac{v}{2}\text{ if }v\in[1,2].
\end{cases}$$
In this case  $\delta_P(S)=1\text{ if }P\in A\backslash B$.
\end{lemma}

\begin{proof}
The Zariski Decomposition follows from
$-K_S-vA\sim_{\DR} (2-v)A+2B_1+C_1+B$. We have
$S_S(A)=1$.
Thus, $\delta_P(S)\le 1$ for $P\in A$.  Note that for $P\in A\backslash B$:
$$h(v)\le \begin{cases}
     \frac{(v + 2)^2}{8} \text{ if }v\in[0,1],\\
     \frac{3 (4 - v) v}{8} \text{ if }v\in[1,2].
    \end{cases}$$
 So $S(W_{\bullet,\bullet}^{A};P)\le \frac{13}{15}<1$.
Thus, $\delta_P(S)=1$ if $P\in A\backslash B$.
\end{proof}
\begin{lemma}\label{deg5-3031_32_2_points}
Suppose $P$ belongs to a $(-1)$-curve $A$ and there exist $(-1)$-curves and $(-2)$-curves   which form the following dual graph:
\begin{figure}[h!]
    \centering
   \includegraphics[width=3.5cm]{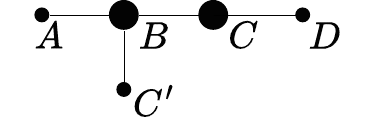}
    \caption{Dual graph: $(-K_S)^2=5$ and $\delta_P(S)=\frac{30}{31}$}
\end{figure}
\par Then $\tau(A)=2$ and the Zariski Decomposition of the divisor $-K_S-vA$ is given by:
{\allowdisplaybreaks \begin{align*}
&&P(v)=\begin{cases}
-K_S-vA-\frac{v}{3}(2B+C)\text{ if }v\in\big[0,\frac{3}{2}\big],\\
-K_S-vA-(v-1)(2B+C)-(2v-3)C'\text{ if }v\in\big[\frac{3}{2},2\big].
\end{cases}\\&&N(v)=\begin{cases}
\frac{v}{3}(2B+C)\text{ if }v\in\big[0,\frac{3}{2}\big],\\
(v-1)(2B+C)+(2v-3)C'\text{ if }v\in\big[\frac{3}{2},2\big].
\end{cases}
\end{align*}}
Moreover, 
$$(P(v))^2=\begin{cases}
5-2v-\frac{v^2}{3} \text{ if }v\in\big[0,\frac{3}{2}\big],\\
(2-v)(4-v)\text{ if }v\in\big[\frac{3}{2},2\big].
\end{cases}
P(v)\cdot A=\begin{cases}
1+\frac{v}{3} \text{ if }v\in\big[0,\frac{3}{2}\big],\\
3-v\text{ if }v\in\big[\frac{3}{2},2\big].
\end{cases}$$
In this case  $\delta_P(S)=\frac{30}{31}\text{ if }P\in A\backslash B$.
\end{lemma}
\begin{proof}
The Zariski Decomposition follows from
$-K_S-vA\sim_{\DR} (2-v)A+3B+2C+D+2C'$. We have
$S_{S}(A)=\frac{31}{30}$.
Thus, $\delta_P(S)\le \frac{30}{31}$ for $P\in A$. Moreover if $P\in A\backslash B$:
$$h(v)=\begin{cases}
\frac{ (v + 3)^2}{18}\text{ if }v\in\big[0,\frac{3}{2}\big],\\
\frac{(3-v)^2}{2}\text{ if }v\in\big[\frac{3}{2},2\big].
\end{cases}$$
So
$S(W_{\bullet,\bullet}^{A};P)\le\frac{19}{30}\le\frac{31}{30}$.
We get that $\delta_P(S)=\frac{30}{31}$ for $P\in A\backslash B$.
\end{proof}
\begin{lemma}\label{deg5-1516_2_points}
Suppose $P$ belongs to a $(-1)$-curve $A$ and there exist $(-1)$-curves and $(-2)$-curves   which form the following dual graph:
\begin{figure}[h!]
    \centering
   \includegraphics[width=4.3cm]{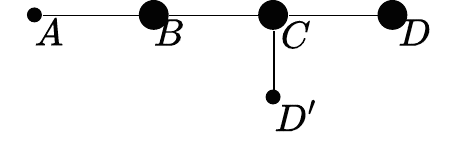}
    \caption{Dual graph: $(-K_S)^2=5$ and $\delta_P(S)=\frac{15}{16}$}
\end{figure}
\par Then $\tau(A)=2$ and the Zariski Decomposition of the divisor $-K_S-vA$ is given by:
$$P(v)=-K_S-vA-\frac{v}{4}(3B+2C+D)\text{ and }N(v)=\frac{v}{4}(3B+2C+D)\text{ if }v\in[0,2].$$
Moreover, 
$$(P(v))^2=\frac{(2-v)(10+v)}{2 }\text{ and }
P(v)\cdot A=1+\frac{v}{4}\text{ if }v\in[0,2].$$
In this case  $\delta_P(S)=\frac{15}{16}\text{ if }P\in A\backslash B$.
\end{lemma}
\begin{proof}
The Zariski Decomposition follows from
$-K_S-vA\sim_{\DR} (2-v)A+3B+4C+2D+3D'$. We have $S_{S}(A)=\frac{16}{15}$.
Thus, $\delta_P(S)\le \frac{15}{16}$ for $P\in A$. Moreover if $P\in A\backslash B$:
$$h(v)=\frac{(v+4)^2}{32}\text{ if }v\in[0,2].$$
So
$S(W_{\bullet,\bullet}^{A};P)\le\frac{19}{30}\le\frac{16}{15}$.
We get that $\delta_P(S)=\frac{15}{16}$ for $P\in A\backslash B$.
\end{proof}
\begin{lemma}\label{deg5-1517_1_2_points}
Suppose $P$ belongs to a $(-2)$-curve $A$ and there exist $(-1)$-curves and $(-2)$-curves   which form the following dual graph:
\begin{figure}[h!]
    \centering
   \includegraphics[width=13.3cm]{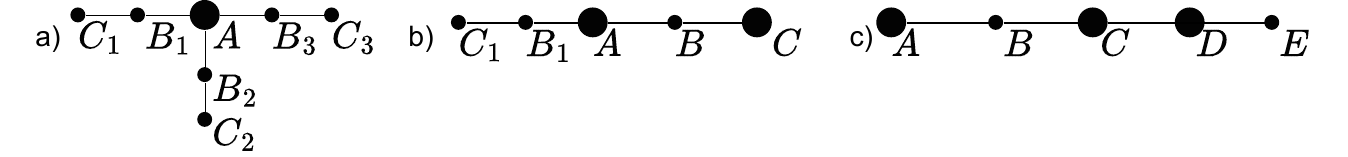}
    \caption{Dual graph: $(-K_S)^2=5$ and $\delta_P(S)=\frac{15}{17}$ with $\tau(A)=2$}
\end{figure}
 \par Then $\tau(A)=2$ and the Zariski Decomposition of the divisor $-K_S-vA$ is given by:{
{\allowdisplaybreaks \begin{align*}
&{\text{\bf a). }}&P(v)=\begin{cases}
-K_{S}-vA\text{ if }v\in[0,1],\\
-K_{S}-vA-(v-1)(B_1+B_2+B_3)\text{ if }v\in[1,2].
\end{cases} 
\\&&N(v)=\begin{cases}
0\text{ if }v\in[0,1],\\
(v-1)(B_1+B_2+B_3)\text{ if }v\in[1,2].
\end{cases}\\
&{\text{\bf b). }}&P(v)=\begin{cases}
-K_{S}-vA\text{ if }v\in[0,1],\\
-K_{S}-vA-(v-1)(B_1+2B+C)\text{ if }v\in[1,2].
\end{cases} 
\\&&N(v)=\begin{cases}
0\text{ if }v\in[0,1],\\
(v-1)(B_1+2B+C)\text{ if }v\in[1,2].
\end{cases}\\
&{\text{\bf c). }}&P(v)=\begin{cases}
-K_{S}-vA\text{ if }v\in[0,1],\\
-K_{S}-vA-(v-1)(3B+2C+D)\text{ if }v\in[1,2].
\end{cases} 
\\&&N(v)=\begin{cases}
0\text{ if }v\in[0,1],\\
(v-1)(3B+2C+D)\text{ if }v\in[1,2].
\end{cases}
\end{align*}}}
Moreover, 
$$(P(v))^2=\begin{cases}5-2v^2\text{ if }v\in[0,1],\\
(2-v)(4-v)\text{ if }v\in[1,2].
\end{cases}P(v)\cdot A=\begin{cases}2v\text{ if }v\in[0,1],\\
3-v\text{ if }v\in[1,2].
\end{cases}$$
In this case $\delta_P(S)=\frac{15}{17}\text{ if }P\in A\backslash B$.
\end{lemma}
\begin{proof}
The Zariski Decomposition in part a). follows from
$$-K_S-vA\sim_{\DR} (2-v)A+\frac{1}{3}\Big(4B_1+4B_2+4B_3+C_1+C_2+C_3\Big).$$
A similar statement holds in other parts. We have
$S_S(A)=\frac{17}{15}$.
Thus, $\delta_P(S)\le \frac{15}{17}$ for $P\in A$. Note that we have:
$$h(v)\le \begin{cases} 2v^2\text{ if }v\in[0,1],\\
 \frac{(3 - v) (v + 1)}{2} \frac{}{}\text{ if }v\in[1,2].
\end{cases}$$
So 
$S(W_{\bullet,\bullet}^{A};P)\le1\le \frac{17}{15}$.
Thus, $\delta_P(S)=\frac{15}{17}$ if $P\in A\backslash B$.
\end{proof}
\begin{lemma}\label{deg5-1517_1_3_points}
Suppose $P$ belongs to a $(-1)$-curve $A$ and there exist $(-1)$-curves and $(-2)$-curves   which form the following dual graph:
\begin{figure}[h!]
    \centering
   \includegraphics[width=3.5cm]{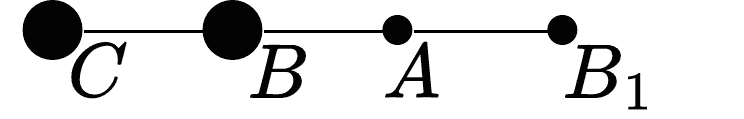}
    \caption{Dual graph: $(-K_S)^2=5$ and $\delta_P(S)=\frac{15}{17}$ with $\tau(A)=3$ }
\end{figure}
\par Then $\tau(A)=3$ and the Zariski Decomposition of the divisor $-K_S-vA$ is given by:
\begin{align*}
&&P(v)=\begin{cases}
-K_S-vA-\frac{v}{3}(2B+C)\text{ if }v\in[0,1],\\
-K_S-vA-\frac{v}{3}(2B+C)-(v-1)B_1\text{ if }v\in[1,3].
\end{cases}
\\&& N(v)=\begin{cases}
\frac{v}{3}(2B+C)\text{ if }v\in[0,1],\\
\frac{v}{3}(2B+C)+(v-1)B_1\text{ if }v\in[1,3].
\end{cases}
\end{align*}
Moreover, 
$$(P(v))^2=\begin{cases}
5-2v-\frac{v^2}{3} \text{ if }v\in[0,1],\\
\frac{2(3-v)^2}{3} \text{ if }v\in[1,3].
\end{cases}
P(v)\cdot A=\begin{cases}
1+\frac{v}{3} \text{ if }v\in[0,1],\\
2(1-\frac{v}{3}) \text{ if }v\in[1,3].
\end{cases}$$
In this case  $\delta_P(S)=\frac{15}{17}\text{ if }P\in A\backslash B$.
\end{lemma}
\begin{proof}
The Zariski Decomposition follows from
$-K_S-vA\sim_{\DR} (3-v)A+2B+C+2B_1$. We have
$S_{S}(A)=\frac{17}{15}$.
Thus, $\delta_P(S)\le \frac{15}{17}$ for $P\in A$. Moreover, for $P\in A\backslash B$:
$$h(v)\le \begin{cases}
\frac{ (v + 3)^2}{18}\text{ if }v\in[0,1],\\
\frac{4(3 -v )}{9} \text{ if }v\in[1,3].
\end{cases}$$
So
$S(W_{\bullet,\bullet}^{A};P)\le \frac{13}{15}\le  \frac{17}{15}$.
We get that $\delta_P(S)=\frac{15}{17}$ for $P\in A\backslash B$.
\end{proof}
\begin{lemma}\label{deg5-1519_2_3_points}
Suppose $P$ belongs to a $(-1)$-curve $A$ and there exist $(-1)$-curves and $(-2)$-curves   which form the following dual graph:
\begin{figure}[h!]
    \centering
   \includegraphics[width=10cm]{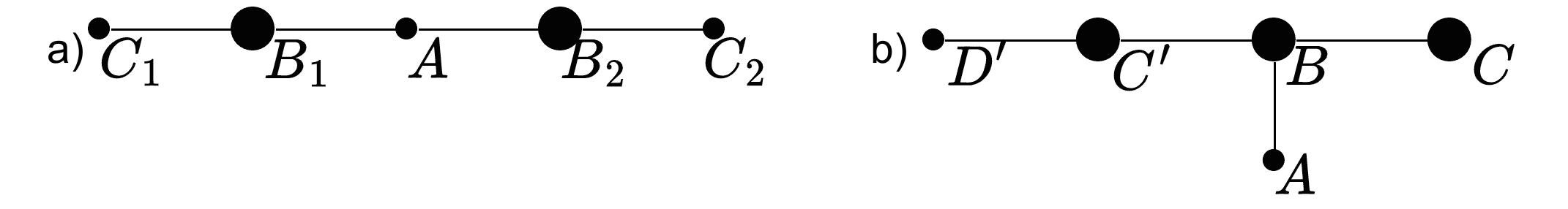}
    \caption{Dual graph: $(-K_S)^2=5$ and $\delta_P(S)=\frac{15}{19}$ with $\tau(A)=3$}
\end{figure}
 \par Then $\tau(A)=3$ and the Zariski Decomposition of the divisor $-K_S-vA$ is given by:
{\allowdisplaybreaks \begin{align*}
&{\text{\bf a). }}&P(v)=\begin{cases}
-K_S-vA-\frac{v}{2}(B_1+B_2)\text{ if }v\in[0,2],\\
-K_S-vA-(v-1)(B_1+B_2)-(v-2)(C_1+C_2)\text{ if }v\in[2,3].
\end{cases}\\&&N(v)=\begin{cases}
\frac{v}{2}(B_1+B_2)\text{ if }v\in[0,2],\\
(v-1)(B_1+B_2)+(v-2)(C_1+C_2)\text{ if }v\in[2,3].
\end{cases}\\
&{\text{\bf b). }}&P(v)=\begin{cases}
-K_S-vA-\frac{v}{2}(2B+C+C')\text{ if }v\in[0,2],\\
-K_S-vA-(v-1)(2B+C)-(2v-3)C'-(2v-4)D'\text{ if }v\in[2,3].
\end{cases}\\&&N(v)=\begin{cases}
\frac{v}{2}(2B+C+C')\text{ if }v\in[0,2],\\
(v-1)(2B+C)+(2v-3)C'+(2v-4)D'\text{ if }v\in[2,3].
\end{cases}
\end{align*}}
Moreover, 
$$(P(v))^2=\begin{cases}
5-2v\text{ if }v\in[0,2],\\
(3-v)^2\text{ if }v\in[2,3].
\end{cases}
P(v)\cdot A=\begin{cases}
1\text{ if }v\in[0,2],\\
3-v\text{ if }v\in[2,3].
\end{cases}$$
In this case  $\delta_P(S)=\frac{15}{19}\text{ if }P\in A\backslash B$.
\end{lemma}
\begin{proof}
The Zariski Decomposition in part a). follows from
$-K_S-vA\sim_{\DR} (3-v)A+2B_1+C_1+2B_2+C_2$. A similar statement holds in other parts.
We have $S_{S}(A)=\frac{19}{15}$.
Thus, $\delta_P(S)\le \frac{15}{19}$ for $P\in A$. Moreover,
$$h(v)\le \begin{cases}
\frac{1+v}{2}\text{ if }v\in[0,2],\\
 \frac{ (3 - v) (v + 1)}{2} \text{ if }v\in[2,3].
\end{cases}$$
So 
$S(W_{\bullet,\bullet}^{A};P)\le \frac{17}{15}\le\frac{19}{15} $.
We get that $\delta_P(S)=\frac{15}{19}$ for $P\in A$.
\end{proof}
\begin{lemma}\label{deg5-1519_1_2_points}
Suppose $P$ belongs to a $(-2)$-curve $A$ and there exist $(-1)$-curves and $(-2)$-curves   which form the following dual graph:
\begin{figure}[h!]
    \centering
   \includegraphics[width=3.5cm]{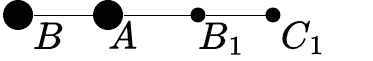}
    \caption{Dual graph: $(-K_S)^2=5$ and $\delta_P(S)=\frac{15}{19}$ with $\tau(A)=2$ }
\end{figure}
 \par Then $\tau(A)=2$ and the Zariski Decomposition of the divisor $-K_S-vA$ is given by:
$$P(v)=\begin{cases}
-K_S-vA-\frac{v}{2}B\text{ if }v\in[0,1],\\
-K_S-vA-\frac{v}{2}B-(v-1)B_1\text{ if }v\in[1,2].
\end{cases}N(v)=\begin{cases}
\frac{v}{2}B\text{ if }v\in[0,1],\\
\frac{v}{2}B+(v-1)B_1\text{ if }v\in[1,2].
\end{cases}$$
Moreover, 
$$(P(v))^2=\begin{cases}
5-\frac{3v^2}{2} \text{ if }v\in[0,1],\\
\frac{(2-v)(6+v)}{2}\text{ if }v\in[1,2].
\end{cases}
P(v)\cdot A=\begin{cases}
\frac{3v}{2} \text{ if }v\in[0,1],\\
1+\frac{v}{2}\text{ if }v\in[1,2].
\end{cases}$$
In this case  $\delta_P(S)=\frac{15}{19}\text{ if }P\in A\backslash B$.
\end{lemma}
\begin{proof}
The Zariski Decomposition follows from
$-K_S-vA\sim_{\DR} (2-v)A+3B_1+2C_1+B$. We have
$S_{S}(A)=\frac{19}{15}$.
Thus, $\delta_P(S)\le \frac{15}{19}$ for $P\in A$. Moreover, if $P\in A\backslash B$ then:
$$h(v)\le \begin{cases}
\frac{9v^2}{8}\text{ if }v\in[0,1],\\
\frac{(v + 2) (5 v - 2)}{8} \text{ if }v\in[1,2].
\end{cases}$$
So
$S(W_{\bullet,\bullet}^{A};P)\le\frac{17}{15}\le\frac{19}{15}$.
We get that $\delta_P(S)=\frac{15}{19}$ for $P\in A\backslash B$.
\end{proof}
\begin{lemma}\label{deg5-1013_32_2_points}
Suppose $P$ belongs to a $(-2)$-curve $A$ and there exist $(-1)$-curves and $(-2)$-curves   which form the following dual graph:
\begin{figure}[h!]
    \centering
   \includegraphics[width=4cm]{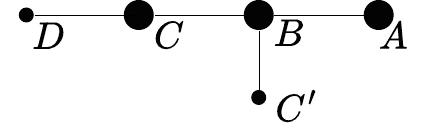}
    \caption{Dual graph: $(-K_S)^2=5$ and $\delta_P(S)=\frac{10}{13}$}
\end{figure}
\par Then $\tau(A)=2$ and the Zariski Decomposition of the divisor $-K_S-vA$ is given by:
{\allowdisplaybreaks \begin{align*}
&&P(v)=\begin{cases}
-K_S-vA-\frac{v}{3}(2B+C)\text{ if }v\in\big[0,\frac{3}{2}\big],\\
-K_S-vA-(v-1)(2B+C)-(2v-3)C'\text{ if }v\in\big[\frac{3}{2},2\big].
\end{cases}\\&&N(v)=\begin{cases}
\frac{v}{3}(2B+C)\text{ if }v\in\big[0,\frac{3}{2}\big],\\
(v-1)(2B+C)+(2v-3)C'\text{ if }v\in\big[\frac{3}{2},2\big].
\end{cases}
\end{align*}}
Moreover, 
$$(P(v))^2=\begin{cases}
5-\frac{4v^2}{3}\text{ if }v\in\big[0,\frac{3}{2}\big],\\
4(2-v)\text{ if }v\in\big[\frac{3}{2},2\big].
\end{cases}
P(v)\cdot A=\begin{cases}
\frac{4v}{3}\text{ if }v\in\big[0,\frac{3}{2}\big],\\
2\text{ if }v\in\big[\frac{3}{2},2\big].
\end{cases}$$
In this case  $\delta_P(S)=\frac{10}{13}\text{ if }P\in A\backslash B$.
\end{lemma}
\begin{proof}
The Zariski Decomposition follows from
$-K_S-vA\sim_{\DR} (2-v)A+4B+3C+2D+3C'$.
We have
$S_{S}(A)=\frac{13}{10}$.
Thus, $\delta_P(S)\le \frac{10}{13}$ for $P\in A$. Moreover if $P\in A\backslash B$:
$$h(v)=\begin{cases}
\frac{8v^2}{9}\text{ if }v\in\big[0,\frac{3}{2}\big],\\
2\text{ if }v\in\big[\frac{3}{2},2\big].
\end{cases}$$
So
$S(W_{\bullet,\bullet}^{A};P)\le\frac{4}{5}\le\frac{13}{10} $.
We get that $\delta_P(S)=\frac{10}{13}$ for $P\in A\backslash B$.
\end{proof}
\begin{lemma}\label{deg5-34_2_points}
Suppose $P$ belongs to a $(-2)$-curve $A$ and there exist $(-1)$-curves and $(-2)$-curves   which form the following dual graph:
\begin{figure}[h!]
    \centering
   \includegraphics[width=4cm]{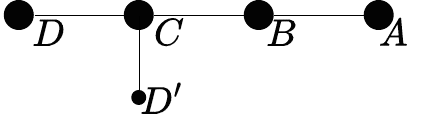}
    \caption{Dual graph: $(-K_S)^2=5$ and $\delta_P(S)=\frac{3}{4}$}
\end{figure}
\par Then $\tau(A)=2$ and the Zariski Decomposition of the divisor $-K_S-vA$ is given by:
$$P(v)=-K_S-vA-\frac{v}{4}(3B+2C+D)\text{ and }N(v)=\frac{v}{4}(3B+2C+D)\text{ if }v\in[0,2].$$
Moreover, 
$$(P(v))^2=\frac{5(2-v)(2+v)}{2}\text{ and }
P(v)\cdot A=\frac{5v}{4}\text{ if }v\in[0,2].$$
In this case  $\delta_P(S)=\frac{3}{4}\text{ if }P\in A\backslash B$.
\end{lemma}

\begin{proof}
The Zariski Decomposition follows from
$-K_S-vA\sim_{\DR} (2-v)A+4B+6C+3D+5D'$. 
We have $S_{S}(A)=\frac{4}{3}$.
Thus, $\delta_P(S)\le \frac{3}{4}$ for $P\in A$. Moreover if $P\in A\backslash B$ then 
$h(v)=\frac{25v^2}{32}\text{ if }v\in[0,2]$.
So
$S(W_{\bullet,\bullet}^{A};P)\le \frac{5}{6}\le\frac{4}{3}$.
We get that $\delta_P(S)=\frac{3}{4}$ for $P\in A\backslash B$.
\end{proof}
\begin{lemma}\label{deg5-57_1_2_3_points}
Suppose $P$ belongs to a $(-2)$-curve $A$ and there exist $(-1)$-curves and $(-2)$-curves   which form the following dual graph:
\begin{figure}[h!]
    \centering
   \includegraphics[width=9.3cm]{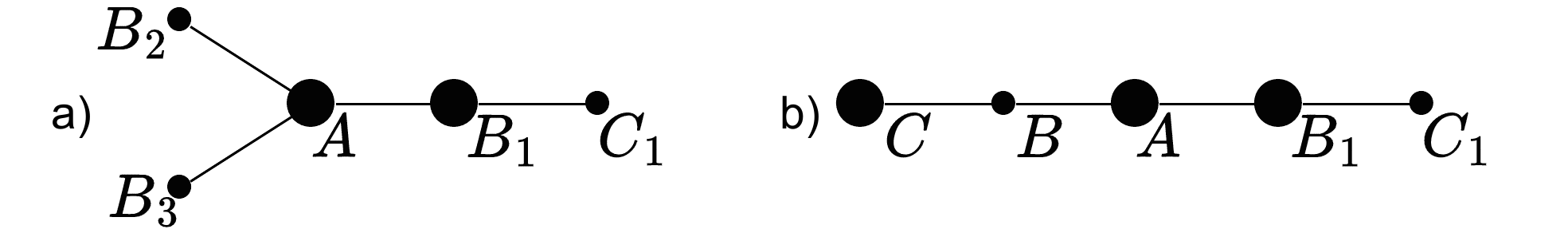}
    \caption{Dual graph: $(-K_S)^2=5$ and $\delta_P(S)=\frac{5}{7}$}
\end{figure}
 \par Then $\tau(A)=3$ and the Zariski Decomposition of the divisor $-K_S-vA$ is given by:
{\allowdisplaybreaks \begin{align*}
&{\text{\bf a). }}&P(v)=\begin{cases}
-K_S-vA-\frac{v}{2}B_1\text{ if }v\in[0,1],\\
-K_S-vA-\frac{v}{2}B_1-(v-1)(B_2+B_3)\text{ if }v\in[1,2],\\
-K_S-vA-(v-1)(B_1+B_2+B_3)-(v-2)C_1\text{ if }v\in[2,3].
\end{cases}\\&&N(v)=\begin{cases}
\frac{v}{2}B_1\text{ if }v\in[0,1],\\
\frac{v}{2}B_1+(v-1)(B_2+B_3)\text{ if }v\in[1,2],\\
(v-1)(B_1+B_2+B_3)+(v-2)C_1\text{ if }v\in[2,3].
\end{cases}
\\&{\text{\bf b). }}&P(v)=\begin{cases}
-K_S-vA-\frac{v}{2}B_1\text{ if }v\in[0,1],\\
-K_S-vA-\frac{v}{2}B_1-(v-1)(2B_2+C_2)\text{ if }v\in[1,2],\\
-K_S-vA-(v-1)(B_1+2B_2+C_2)-(v-2)C_1\text{ if }v\in[2,3].
\end{cases}\\&&N(v)=\begin{cases}
\frac{v}{2}B_1\text{ if }v\in[0,1],\\
\frac{v}{2}B_1+(v-1)(2B_2+C_2)\text{ if }v\in[1,2],\\
(v-1)(B_1+2B_2+C_2)+(v-2)C_1\text{ if }v\in[2,3].
\end{cases}
\end{align*}}
Moreover, 
$$(P(v))^2=\begin{cases}
5-\frac{3v^2}{2} \text{ if }v\in[0,1],\\
7-4v+\frac{v^2}{2} \text{ if }v\in[1,2],\\
(3-v)^2\text{ if }v\in[2,3].
\end{cases}
P(v)\cdot A=\begin{cases}
\frac{3v}{2} \text{ if }v\in[0,1],\\
2-\frac{v}{2}\text{ if }v\in[1,2],\\
3-v\text{ if }v\in[2,3].
\end{cases}$$
In this case  $\delta_P(S)=\frac{5}{7}\text{ if }P\in A\backslash B$.
\end{lemma}
\begin{proof}
The Zariski Decomposition in part a). follows from
$-K_S-vA\sim_{\DR} (3-v)A+2B_1+C_1+2B_2+2B_3$. A similar statement holds in other parts.
We have $S_{S}(A)=\frac{7}{5}$.
Thus, $\delta_P(S)\le \frac{5}{7}$ for $P\in A$. Moreover,  if $P\in A\cap B_1$ or if $P\in A\backslash B_1$ then:
$$h(v)=\begin{cases}
\frac{15v^2}{8}\text{ if }v\in[0,1],\\
\frac{  (4-v) (v + 4)}{8}\text{ if }v\in[1,2],\\
 \frac{ (3 - v) (v + 1)}{2} \text{ if }v\in[2,3].
\end{cases}
\text{ or }
h(v)\le\begin{cases}
\frac{9v^2}{8}\text{ if }v\in[0,1],\\
 \frac{3 (4 - v) v}{8} \text{ if }v\in[1,2],\\
 \frac{ (3 - v) (v + 1)}{2} \text{ if }v\in[2,3].
\end{cases}
$$
So  $S(W_{\bullet,\bullet}^{A};P)\le \frac{19}{15}\le  \frac{7}{5}$ or $S(W_{\bullet,\bullet}^{A};P)\le \frac{31}{30}\le  \frac{7}{5}$.
We get that $\delta_P(S)=\frac{5}{7}$ for $P\in A$.
\end{proof}
\begin{lemma}\label{deg5-3043_1_32_3_points}
Suppose $P$ belongs to a $(-2)$-curve $A$ and there exist $(-1)$-curves and $(-2)$-curves   which form the following dual graph:
\begin{figure}[h!]
    \centering
   \includegraphics[width=4cm]{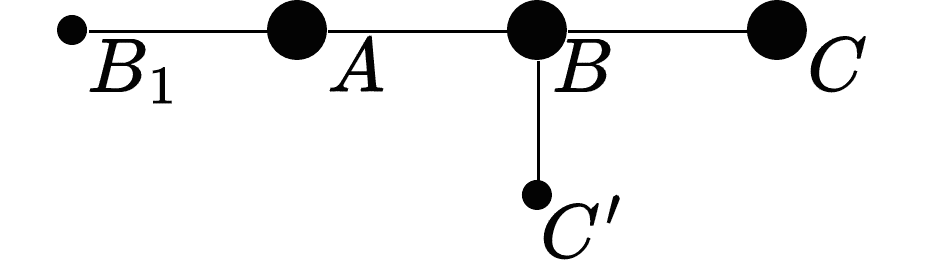}
    \caption{Dual graph: $(-K_S)^2=5$ and $\delta_P(S)=\frac{30}{43}$}
\end{figure}
\par Then $\tau(A)=3$ and the Zariski Decomposition of the divisor $-K_S-vA$ is given by:
{\allowdisplaybreaks \begin{align*}
&&P(v)=\begin{cases}
-K_S-vA-\frac{v}{3}(2B+C)\text{ if }v\in[0,1],\\
-K_S-vA-\frac{v}{3}(2B+C)-(v-1)B_1\text{ if }v\in\big[1,\frac{3}{2}\big],\\
-K_S-vA-(v-1)(2B+C+B_1)-(2v-3)C'\text{ if }v\in\big[\frac{3}{2},3\big].
\end{cases}\\&&N(v)=\begin{cases}
\frac{v}{3}(2B+C)\text{ if }v\in[0,1],\\
\frac{v}{3}(2B+C)+(v-1)B_1\text{ if }v\in\big[1,\frac{3}{2}\big],\\
(v-1)(2B+C+B_1)+(2v-3)C'\text{ if }v\in\big[\frac{3}{2},3\big].
\end{cases}
\end{align*}}
Moreover, 
$$(P(v))^2=\begin{cases}5-\frac{4v^2}{3}\text{ if }v\in[0,1],\\
6-2v-\frac{v^2}{3}\text{ if }v\in\big[1,\frac{3}{2}\big],\\
(3-v)^2\text{ if }v\in\big[\frac{3}{2},3\big].
\end{cases}P(v)\cdot A=
\begin{cases}
\frac{4v}{3}\text{ if }v\in[0,1],\\
1+\frac{v}{3} \text{ if }v\in\big[1,\frac{3}{2}\big],\\
3-v\text{ if }v\in\big[\frac{3}{2},3\big].
\end{cases}$$
In this case  $\delta_P(S)=\frac{30}{43}\text{ if }P\in A\backslash B$.
\end{lemma}
\begin{proof}
The Zariski Decomposition follows from
$-K_S-vA\sim_{\DR} (3-v)A+4B+2C+3C'+2B_1$.
We have $S_S(A)=\frac{43}{30}$.
Thus, $\delta_P(S)\le \frac{30}{43}$ for $P\in A$.  Note that if $P\in A\backslash B$
$$h(v)=\begin{cases}
\frac{8v^2}{9}\text{ if }v\in[0,1],\\
\frac{(v + 3) (7 v - 3)}{18}\text{ if }v\in\big[1,\frac{3}{2}\big],\\
 \frac{ (3 - v) (v + 1)}{2} \text{ if }v\in\big[\frac{3}{2},3\big].
    \end{cases}$$
 So $S(W_{\bullet,\bullet}^{A};P)\le\frac{16}{15}\le  \frac{43}{30}$.
Thus, $\delta_P(S)=\frac{30}{43}$ if $P\in A\backslash B$.
\end{proof}
\begin{lemma}\label{deg5-913_43_3_points}
Suppose $P$ belongs to a $(-2)$-curve $A$ and there exist $(-1)$-curves and $(-2)$-curves   which form the following dual graph:
\begin{figure}[h!]
    \centering
   \includegraphics[width=4cm]{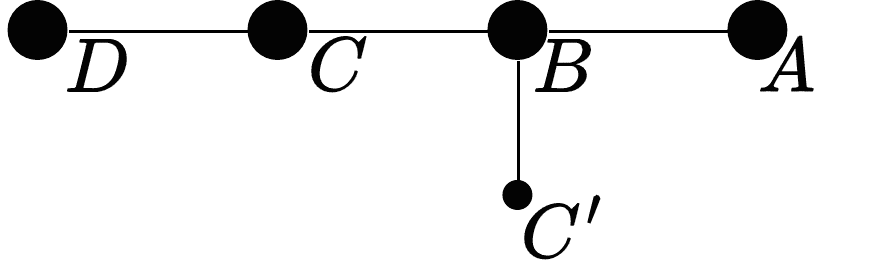}
    \caption{Dual graph: $(-K_S)^2=5$ and $\delta_P(S)=\frac{9}{13}$}
\end{figure}
\par  Then $\tau(A)=3$ and the Zariski Decomposition of the divisor $-K_S-vA$ is given by:
{\allowdisplaybreaks \begin{align*}
&&P(v)=\begin{cases}
-K_S-vA-\frac{v}{4}(3B+2C+D)\text{ if }v\in\big[0,\frac{4}{3}\big],\\
-K_S-vA-(v-1)(3B+2C+D)-(3v-4)C'\text{ if }v\in\big[\frac{4}{3},3\big].
\end{cases}\\&&N(v)=\begin{cases}
\frac{v}{4}(3B+2C+D)\text{ if }v\in\big[0,\frac{4}{3}\big],\\
(v-1)(3B+2C+D)+(3v-4)C'\text{ if }v\in\big[\frac{4}{3},3\big].
\end{cases}
\end{align*}}
Moreover, 
$$(P(v))^2=\begin{cases}
\frac{5(2-v)(2+v)}{2}\text{ if }v\in\big[0,\frac{4}{3}\big],\\
(3-v)^2\text{ if }v\in\big[\frac{4}{3},3\big].
\end{cases}
P(v)\cdot A=\begin{cases}
\frac{5v}{4}\text{ if }v\in\big[0,\frac{4}{3}\big],\\
3-v\text{ if }v\in\big[\frac{4}{3},3\big].
\end{cases}$$
In this case  $\delta_P(S)=\frac{9}{13}\text{ if }P\in A\backslash B$.
\end{lemma}

\begin{proof}
The Zariski Decomposition follows from
$-K_S-vA\sim_{\DR} (3-v)A+6B+4C+5C'+2D$.
We have
$S_{S}(A)=\frac{13}{9}$.
Thus, $\delta_P(S)\le \frac{9}{13}$ for $P\in A$. Moreover if $P\in A\backslash B$:
$$h(v)=\begin{cases}
\frac{25}{32}v^2\text{ if }v\in\big[0,\frac{4}{3}\big],\\
\frac{(3-v)^2}{2}\text{ if }v\in\big[\frac{4}{3},3\big].
\end{cases}$$
So $S(W_{\bullet,\bullet}^{A};P)\le\frac{5}{9}\le\frac{13}{9}$. We get that $\delta_P(S)=\frac{9}{13}$ for $P\in A\backslash B$.
\end{proof}
\begin{lemma}\label{deg5-1523_3_4_points}
Suppose $P$ belongs to a $(-1)$-curve $A$ and there exist $(-1)$-curves and $(-2)$-curves   which form the following dual graph:
\begin{figure}[h!]
    \centering
   \includegraphics[width=5.5cm]{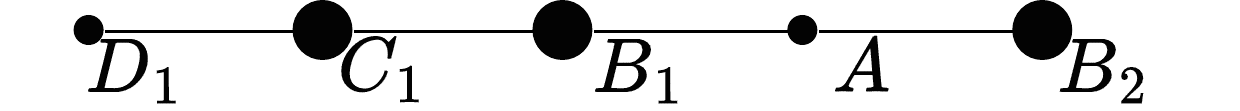}
    \caption{Dual graph: $(-K_S)^2=5$ and $\delta_P(S)=\frac{15}{23}$}
\end{figure}
\par Then $\tau(A)=4$ and the Zariski Decomposition of the divisor $-K_S-vA$ is given by:
{\allowdisplaybreaks \begin{align*}
&&P(v)=\begin{cases}
-K_S-vA-\frac{v}{3}(2B_1+C_1)-\frac{v}{2}B_2\text{ if }v\in[0,3],\\
-K_S-vA-(v-1)B_1-(v-2)C_1-(v-3)D_1-\frac{v}{2}B_2\text{ if }v\in[3,4].
\end{cases}\\&&N(v)=\begin{cases}
\frac{v}{3}(2B_1+C_1)+\frac{v}{2}B_2\text{ if }v\in[0,3],\\
(v-1)B_1+(v-2)C_1+(v-3)D_1+\frac{v}{2}B_2\text{ if }v\in[3,4].
\end{cases}
\end{align*}}
Moreover, 
$$(P(v))^2=\begin{cases}
5-2v+\frac{v^2}{6}\text{ if }v\in[0,3],\\
\frac{(4-v)^2}{2}\text{ if }v\in[3,4].
\end{cases}
P(v)\cdot A=\begin{cases}
1-\frac{v}{6}\text{ if }v\in[0,3],\\
2-\frac{v}{2}\text{ if }v\in[3,4].
\end{cases}$$
In this case  $\delta_P(S)=\frac{15}{23}\text{ if }P\in A$.
\end{lemma}
\begin{proof}
The Zariski Decomposition follows from
$-K_S-vA\sim_{\DR} (4-v)A+3B_1+2C_1+D_1+2B_2$. We have
$S_{S}(A)=\frac{23}{15}$.
Thus, $\delta_P(S)\le \frac{15}{23}$ for $P\in A$. Moreover,
if $P\in A\backslash B_1$ or if $P\in A\cap B_1$:
$$h(v)\le \begin{cases}
\frac{ (6 - v) (5 v + 6)}{72}\text{ if }v\in[0,3],\\
\frac{  (4-v) (v + 4)}{8}\text{ if }v\in[3,4].
\end{cases}\text{ or }h(v)\le \begin{cases}
\frac{(6 - v) (7 v + 6)}{72}\text{ if }v\in[0,3],\\
 \frac{3 (4 - v) v}{8} \text{ if }v\in[3,4].
\end{cases}$$
So $S(W_{\bullet,\bullet}^{A};P)\le\frac{17}{15}\le\frac{23}{15} $ or $S(W_{\bullet,\bullet}^{A};P)\le\frac{7}{5}\le\frac{23}{15}$.
We get that $\delta_P(S)=\frac{15}{23}$ for $P\in A$.
\end{proof}
\begin{lemma}\label{deg5-35_4_points}
Suppose $P$ belongs to a $(-1)$-curve $A$ and there exist $(-1)$-curves and $(-2)$-curves   which form the following dual graph:
\begin{figure}[h!]
    \centering
   \includegraphics[width=4cm]{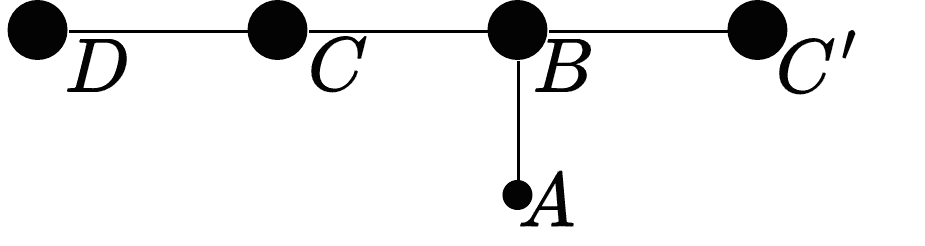}
    \caption{Dual graph: $(-K_S)^2=5$ and $\delta_P(S)=\frac{3}{5}$}
\end{figure}
\par Then $\tau(A)=5$ and the Zariski Decomposition of the divisor $-K_S-vA$ is given by:
$$P(v)=-K_S-vA-\frac{v}{5}(2D+4C+6B+3C')\text{ and }N(v)=\frac{v}{5}(2D+4C+6B+3C')\text{ if }v\in[0,5].$$
Moreover, 
$$(P(v))^2=\frac{(5-v)^2}{2}\text{ and }
P(v)\cdot A=1-\frac{v}{5}\text{ if }v\in[0,5].$$
In this case  $\delta_P(S)=\frac{3}{5}\text{ if }P\in A\backslash B$.
\end{lemma}

\begin{proof}
The Zariski Decomposition follows from
$-K_S-vA\sim_{\DR} (5-v)A+3C'+6B+4C+2D$. We have
$S_{S}(A)=\frac{5}{3}$.
Thus, $\delta_P(S)\le \frac{3}{5}$ for $P\in A$. Moreover if $P\in A\backslash B$ then
$h(v)=\frac{ (5 - v)^2}{50}\text{ if }v\in[0,5]$. So 
$S(W_{\bullet,\bullet}^{A};P)\le \frac{1}{3}\le\frac{5}{3}$.
We get that $\delta_P(S)=\frac{3}{5}$ for $P\in A\backslash B$.
\end{proof}
\begin{lemma}\label{deg5-59_1_2_4_points}
Suppose $P$ belongs to a $(-2)$-curve $A$ and there exist $(-1)$-curves and $(-2)$-curves   which form the following dual graph:
\begin{figure}[h!]
    \centering
   \includegraphics[width=4.5cm]{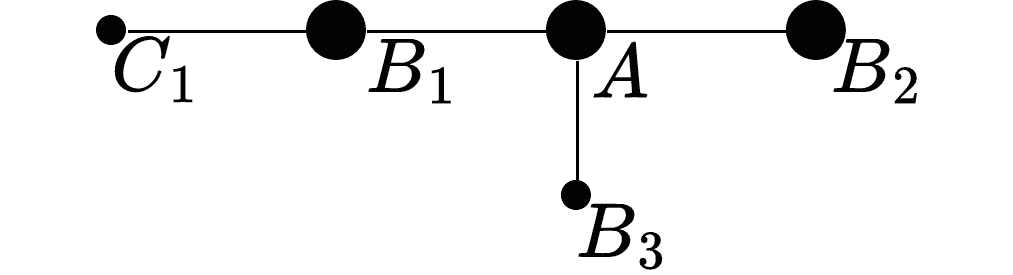}
    \caption{Dual graph: $(-K_S)^2=5$ and $\delta_P(S)=\frac{5}{9}$}
\end{figure}
\par  Then $\tau(A)=4$ and the Zariski Decomposition of the divisor $-K_S-vA$ is given by:
{\allowdisplaybreaks \begin{align*}
&&P(v)=\begin{cases}
-K_S-vA-\frac{v}{2}(B_1+B_2)\text{ if }v\in[0,1],\\
-K_S-vA-\frac{v}{2}(B_1+B_2)-(v-1)B_3\text{ if }v\in[1,2],\\
-K_S-vA-(v-1)(B_1+B_3)-\frac{v}{2}B_2-(v-2)C_1\text{ if }v\in[2,4].
\end{cases}\\&&N(v)=\begin{cases}
\frac{v}{2}(B_1+B_2)\text{ if }v\in[0,1],\\
\frac{v}{2}(B_1+B_2)+(v-1)B_3\text{ if }v\in[1,2],\\
(v-1)(B_1+B_3)+\frac{v}{2}B_2+(v-2)C_1\text{ if }v\in[2,4].
\end{cases}
\end{align*}}
Moreover, 
$$(P(v))^2=\begin{cases}
5-v^2\text{ if }v\in[0,1],\\
2(3-v)\text{ if }v\in[1,2],\\
\frac{(4-v)^2}{2}\text{ if }v\in[2,4].
\end{cases}
P(v)\cdot A=\begin{cases}
v\text{ if }v\in[0,1],\\
1\text{ if }v\in[1,2],\\
2-\frac{v}{2}\text{ if }v\in[2,4].
\end{cases}$$
In this case  $\delta_P(S)=\frac{5}{9}\text{ if }P\in A$.
\end{lemma}
\begin{proof}
The Zariski Decomposition follows from
$-K_S-vA\sim_{\DR} (4-v)A+3B_1+2C_1+2B_2+3B_3$.  We have
$S_{S}(A)=\frac{9}{5}$.
Thus, $\delta_P(S)\le \frac{5}{9}$ for $P\in A$. Moreover,
if $P\in A\backslash (B_2\cup B_3)$ or if  $P\in A\cap B_2$ or if $P\in A\cap B_3$:
$$\hspace*{-0.7cm}h(v)\le\begin{cases}
v^2\text{ if }v\in[0,1],\\
\frac{(v+1)}{2} \text{ if }v\in[1,2],\\
 \frac{3 (4 - v) v}{8} \text{ if }v\in[2,4].
\end{cases}
\text{or }
h(v)=\begin{cases}
v^2\text{ if }v\in[0,1],\\
\frac{(v+1)}{2} \text{ if }v\in[1,2],\\
\frac{  (4-v) (v + 4)}{8}\text{ if }v\in[2,4].
\end{cases}
\text{or }
h(v)=\begin{cases}
\frac{v^2}{2} \text{ if }v\in[0,1],\\
v-\frac{1}{2} \text{ if }v\in[1,2],\\
 \frac{3 (4 - v) v}{8} \text{ if }v\in[2,4].
\end{cases}$$
So $S(W_{\bullet,\bullet}^{A};P)\le \frac{43}{30}\le\frac{9}{5} $ or $S(W_{\bullet,\bullet}^{A};P)\le \frac{13}{10}\le\frac{9}{5} $ or $S(W_{\bullet,\bullet}^{A};P)\le \frac{19}{15}\le\frac{9}{5}$.
We get that $\delta_P(S)=\frac{5}{9}$ for $P\in A$.
\end{proof}
\begin{lemma}\label{deg5-611_32_4_points}
Suppose $P$ belongs to a $(-2)$-curve $A$ and there exist $(-1)$-curves and $(-2)$-curves   which form the following dual graph:
\begin{figure}[h!]
    \centering
   \includegraphics[width=4.5cm]{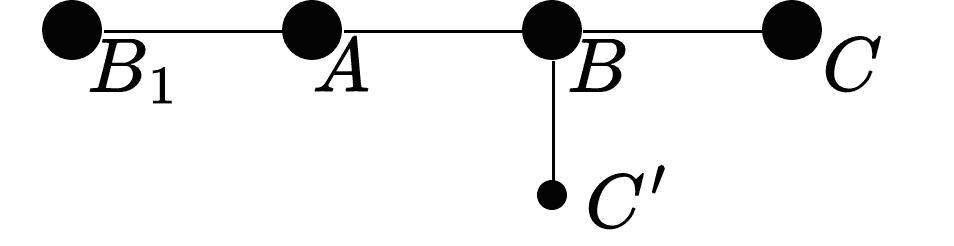}
    \caption{Dual graph: $(-K_S)^2=5$ and $\delta_P(S)=\frac{6}{11}$}
\end{figure}
\par  Then $\tau(A)=4$ and the Zariski Decomposition of the divisor $-K_S-vA$ is given by:
{\allowdisplaybreaks \begin{align*}
&&P(v)=\begin{cases}
-K_S-vA-\frac{v}{2}B_1-\frac{v}{3}(2B+C)\text{ if }v\in\big[0,\frac{3}{2}\big],\\
-K_S-vA-\frac{v}{2}B_1-(v-1)(2B+C)-(2v-3)C'\text{ if }v\in\big[\frac{3}{2},4\big].
\end{cases}\\&&N(v)=\begin{cases}
\frac{v}{2}B_1+\frac{v}{3}(2B+C)\text{ if }v\in\big[0,\frac{3}{2}\big],\\
\frac{v}{2}B_1+(v-1)(2B+C)+(2v-3)C'\text{ if }v\in\big[\frac{3}{2},4\big].
\end{cases}
\end{align*}}
Moreover, 
$$(P(v))^2=\begin{cases}
5-\frac{5v^2}{2} \text{ if }v\in\big[0,\frac{3}{2}\big],\\
\frac{(4-v)^2}{2}\text{ if }v\in\big[\frac{3}{2},4\big].
\end{cases}
P(v)\cdot A=\begin{cases}
\frac{5v}{6}\text{ if }v\in\big[0,\frac{3}{2}\big],\\
2-\frac{v}{2}\text{ if }v\in \big[\frac{3}{2},4\big].
\end{cases}$$
In this case  $\delta_P(S)=\frac{6}{11}\text{ if }P\in A\backslash B$.
\end{lemma}

\begin{proof}
The Zariski Decomposition follows from
$-K_S-vA\sim_{\DR} (4-v)A+2B_1+6B+3C+5C'$. We have
$S_{S}(A)=\frac{11}{6}$.
Thus, $\delta_P(S)\le \frac{6}{11}$ for $P\in A$. Moreover,
if $P\in A\backslash B$ then:
$$h(v)\le \begin{cases}
\frac{55v^2}{72}\text{ if }v\in\big[0,\frac{3}{2}\big],\\
\frac{  (4-v) (v + 4)}{8}\text{ if }v\in\big[\frac{3}{2},4\big].
\end{cases}$$
So
$S(W_{\bullet,\bullet}^{A};P)\le \frac{4}{3}\le\frac{11}{6}$.
We get that $\delta_P(S)=\frac{6}{11}$ for $P\in  A\backslash B$.
\end{proof}
\begin{lemma}\label{deg5-37_1_6_points}
Suppose $P$ belongs to a $(-2)$-curve $A$ and there exist $(-1)$-curves and $(-2)$-curves   which form the following dual graph:
\begin{figure}[h!]
    \centering
   \includegraphics[width=4.3cm]{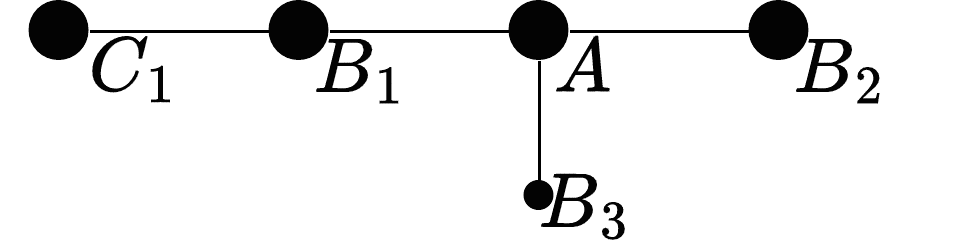}
    \caption{Dual graph: $(-K_S)^2=5$ and $\delta_P(S)=\frac{3}{7}$}
\end{figure}
\par Then $\tau(A)=6$ and the Zariski Decomposition of the divisor $-K_S-vA$ is given by:
{\allowdisplaybreaks \begin{align*}
&&P(v)=\begin{cases}
-K_S-vA-\frac{v}{3}(2B_1+C_1)-\frac{v}{2}B_2\text{ if }v\in[0,1],\\
-K_S-vA-\frac{v}{3}(2B_1+C_1)-\frac{v}{2}B_2-(v-1)B_3\text{ if }v\in[1,6].
\end{cases}\\&&N(v)=\begin{cases}
\frac{v}{3}(2B_1+C_1)+\frac{v}{2}B_2\text{ if }v\in[0,1],\\
\frac{v}{3}(2B_1+C_1)+\frac{v}{2}B_2+(v-1)B_3\text{ if }v\in[1,6].
\end{cases}
\end{align*}}
Moreover, 
$$(P(v))^2=\begin{cases}
5-5\frac{v^2}{2} \text{ if }v\in[0,1],\\
\frac{(6-v)^2}{2}\text{ if }v\in[1,6].
\end{cases}
P(v)\cdot A=\begin{cases}
\frac{5v}{6}\text{ if }v\in[0,1],\\
1-\frac{v}{6}\text{ if }v\in[1,6].
\end{cases}$$
In this case  $\delta_P(S)=\frac{3}{7}\text{ if }P\in A$.
\end{lemma}
\begin{proof}
The Zariski Decomposition follows from
$-K_S-vA\sim_{\DR} (6-v)A+4B_1+2C_1+3B_2+5B_3$.
We have $S_{S}(A)=\frac{7}{3}$.
Thus, $\delta_P(S)\le \frac{3}{7}$ for $P\in A$. Moreover, if $P\in A\cap B_3$ or if $P\in A\backslash B_3$ then:
$$h(v)=\begin{cases}
\frac{25v^2}{72}\text{ if }v\in[0,1],\\
\frac{ (6 - v) (11 v - 6)}{72}\text{ if }v\in[1,6].
\end{cases}\text{ or }
h(v)=\begin{cases}
\frac{65v^2}{72}\text{ if }v\in[0,1],\\
\frac{(6 - v) (7 v + 6)}{72}\text{ if }v\in[1,6].
\end{cases}$$
So
$S(W_{\bullet,\bullet}^{A};P)\le\frac{5}{3}\le\frac{7}{3}$
or
$S(W_{\bullet,\bullet}^{A};P)\le\frac{11}{6}\le\frac{7}{3}$.
We get that $\delta_P(S)=\frac{3}{7}$ for $P\in A$.
\end{proof}


\subsection{Finding $\delta$-invariants for degree $5$}

Let $X$ be a singular del Pezzo surface of degree $5$ with and $S$ be a minimal resolution of $X$. Then there are several possible cases:

\begin{itemize}
    \item[I.] $X$ has an $\DA_1$ singularity and contains $7$ lines. In this case, we let $E$ be the exceptional divisor,  $L_{123}$, $L_{i}$, $L_{i}'$ for $i\in\{1,2,3\}$ be the lines on $S$,  
    \item[II.] $X$ has two $\DA_1$ singularities and contains $5$ lines. In this case, we let  $E_i$ for $i\in\{1,2\}$ be the exceptional divisors,  $L_{12}$, $L_{i}$, $L_{i}'$ for $i\in\{1,2\}$ be the lines on $S$,  
    \item[III.] $X$ has an $\DA_2$ singularity and contains $4$ lines. In this case, we let  $E_i$ for $i\in\{1,2\}$ be the exceptional divisors, $L_{i}$, $L_{i}'$ for $i\in\{1,2\}$ be the lines on $S$,  
     \item[IV.] $X$ has  $\DA_2$ and $\DA_1$ singularities and contains $3$ lines. In this case, we let $E_i$ for $i\in\{1,2,3\}$ be the exceptional divisors, $L_{13}$
    $L_{2}$, $L_{2}'$  be the lines on $S$,  
    \item[V.] $X$ has an $\DA_3$ singularity and contains $2$ lines. In this case, we let $E_i$ for $i\in\{1,2,3\}$ be the exceptional divisors,  $L_{1}$ and $L_2$ be the lines on $S$,  
     \item[VI.] $X$ has an $\DA_4$ singularity and contains $1$ line. In this case, we let  $E_i$ for $i\in\{1,2,3,4\}$ be the exceptional divisors,  $L_{3}$ be the line on $S$.  
\end{itemize}
such that the dual graph of the $(-1)$-curves  and $(-2)$-curves on $S$ is given on the picture below:
\begin{figure}[h!]
    \centering
   \includegraphics[width=17cm]{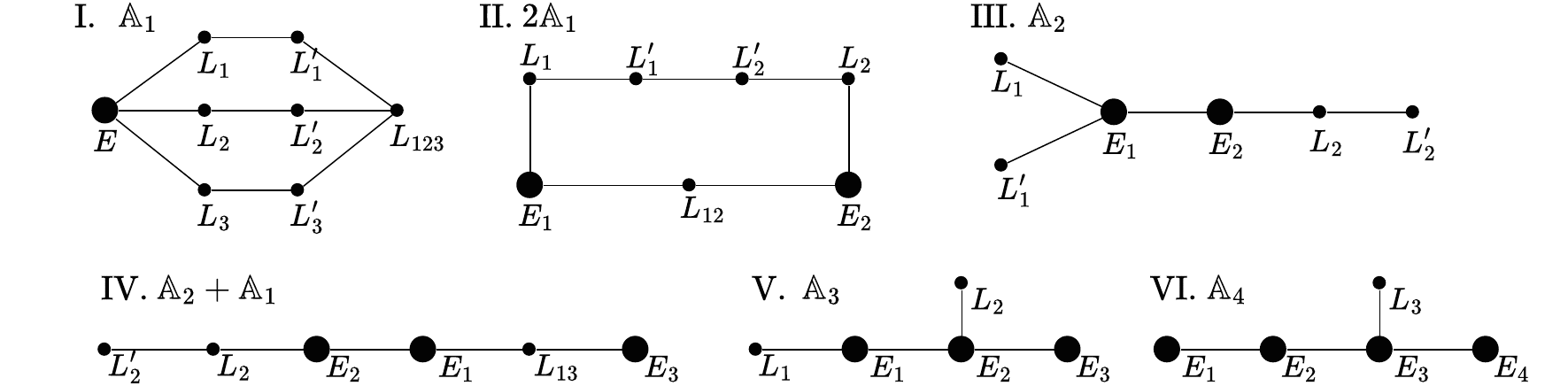}
    \caption{Du Val del Pezzo surfaces with $(-K_S)^2=5$}
\end{figure}
\newpage One has:
\begin{itemize}
    \item[I.] $\delta(X)=\frac{15}{17}$ since depending on the position of point $P\in S$ we have 
\begin{table}[h!]
\hspace*{1cm}
\begin{tabular}{ | c || c | c | c | c |}
 \hline 
 $P$ & $E$ & $(L_1\cup L_2\cup L_3)\backslash E$ & $(L_1'\cup L_2'\cup L_3'\cup L_{123})\backslash (L_1\cup L_2\cup L_3)$ & o/w \\\hline
$\delta_P(S)$ & $\frac{15}{17}$ & $1$ & $\frac{15}{13}$ & $\ge \frac{6}{5}$\\\hline
    \end{tabular}
\caption{Local $\delta$-invariants: $(-K_S)^2=5$ and  $\DA_1$ singularity}
\end{table}
\item[II.] $\delta(X)=\frac{15}{19}$ since  depending on the position of point $P\in S$ we have 
\begin{table}[h!]
\hspace*{1cm}
\begin{tabular}{ | c || c | c | c | c | c |}
 \hline 
 $P$ & $L_{12}$ & $(E_1\cup E_2)\backslash L_{12}$ & $(L_1\cup L_2)\backslash (E_1\cup E_2)$ & $(L_1'\cup L_2')\backslash (L_1\cup L_2)$ & o/w \\\hline
$\delta_P(S)$ & $\frac{15}{19}$ & $\frac{15}{17}$ & $1$ & $\frac{15}{13}$ & $\ge \frac{6}{5}$\\\hline
    \end{tabular}
\caption{Local $\delta$-invariants: $(-K_S)^2=5$ and  $2\DA_1$ singularities}
\end{table}
\item[III.] $\delta(X)=\frac{5}{7}$ since depending on the position of point $P\in S$ we have
\begin{table}[h!]
\hspace*{1cm}
\begin{tabular}{ | c || c | c | c | c | c | c |}
 \hline 
 $P$ & $E_1$ & $E_2\backslash E_1$ & $L_2\backslash E_2$ & $(L_{1}\cup L_1')\backslash E_1$ & $L_2'\backslash L_{2} $ & o/w \\\hline
$\delta_P(S)$ & $\frac{5}{7}$ & $\frac{15}{19}$ & $\frac{15}{17}$ & $\frac{30}{31}$ & $\frac{15}{13}$ & $\ge \frac{6}{5}$\\\hline
    \end{tabular}
\caption{Local $\delta$-invariants: $(-K_S)^2=5$ and  $\DA_2$ singularity}
\end{table}
\item[IV.] $\delta(X)=\frac{15}{23}$ since depending on the position of point $P\in S$ we have
\begin{table}[h!]
\hspace*{1cm}
\begin{tabular}{ | c || c | c | c | c | c | c |}
 \hline 
 $P$ & $L_{13}$ & $E_1\backslash L_{13}$ & $E_2\backslash E_1$ & $(E_3\cup L_2)\backslash (L_{13}\cup E_2)$ & $L_2' \backslash L_2$ & o/w \\\hline
$\delta_P(S)$ & $\frac{15}{23}$ & $\frac{5}{7}$ & $\frac{15}{19}$ & $\frac{15}{17}$  & $\frac{15}{13}$ & $\ge \frac{6}{5}$\\\hline
    \end{tabular}
\caption{Local $\delta$-invariants: $(-K_S)^2=5$ and  $\DA_2\DA_1$ singularities}
\end{table}
\item[V.] $\delta(X)=\frac{5}{9}$ since depending on the position of point $P\in S$ we have
\begin{table}[h!]
\hspace*{1cm}
\begin{tabular}{ | c || c | c | c | c | c | c |}
 \hline 
 $P$ & $E_2$ & $E_1\backslash E_2$ & $ E_3\backslash E_2$ & $ L_2\backslash E_2$ & $L_1\backslash E_1$ & o/w \\\hline
$\delta_P(S)$ & $\frac{5}{9}$ & $\frac{30}{43}$ & $\frac{10}{13}$ & $\frac{15}{19}$  & $\frac{15}{16}$ & $\ge \frac{6}{5}$\\\hline
    \end{tabular}
\caption{Local $\delta$-invariants: $(-K_S)^2=5$ and  $\DA_3$ singularity}
\end{table}
\item[VI.] $\delta(X)=\frac{3}{7}$ since depending on the position of point $P\in S$ we have
\begin{table}[h!]
\hspace*{1cm}
\begin{tabular}{ | c || c | c | c | c | c | c |}
 \hline 
 $P$ & $E_3$ & $E_2\backslash E_3$ & $ L_3\backslash E_3$ & $E_4\backslash E_3$ & $E_1 \backslash E_2$ & o/w \\\hline
$\delta_P(S)$ & $\frac{3}{7}$ & $\frac{6}{11}$ & $\frac{3}{5}$ & $ \frac{9}{13}$  & $\frac{3}{4}$ & $\ge \frac{6}{5}$\\\hline
    \end{tabular}
\caption{Local $\delta$-invariants: $(-K_S)^2=5$ and  $\DA_4$ singularity}
\end{table}
\end{itemize}
\begin{proof}
We prove each case separately using lemmas from the previous section.
    \begin{itemize}
        \item[I.] If $P\in E$ the assertion follows from Lemma \ref{deg5-1517_1_2_points}.
 If $P\in  (L_1\cup L_2\cup L_3)\backslash E$ the assertion follows from Lemma \ref{deg5-1_1_2_points}.
 If $P\in L_{123}\backslash (L_1'\cup L_2'\cup L_3')$ the assertion follows from Lemma \ref{deg5-1513_1_2_points} [a).].
 If $P\in (L_1'\cup L_2'\cup L_3')\backslash (L_1\cup L_2\cup L_3)$ the assertion follows from Lemma \ref{deg5-1513_1_2_points} [b).]. 
 If $P$ is a general point the assertion follows from Lemma \ref{deg5-generalpoint}. 
 \item[II.] If $P\in L_{12}$ the assertion follows from Lemma \ref{deg5-1519_2_3_points}.
 If $P\in (E_1\cup E_2)\backslash L_{12}$ the assertion follows from Lemma \ref{deg5-1517_1_2_points} [b).].
 If $P\in  (L_1\cup L_2)\backslash (E_1\cup E_2)$ the assertion follows from Lemma \ref{deg5-1_1_2_points}.
 If $P\in (L_1'\cup L_2')\backslash (L_1\cup L_2)$ the assertion follows from Lemma \ref{deg5-1513_1_2_points} [b).].
 If $P$ is a general point the assertion follows from Lemma \ref{deg5-generalpoint}. 
\item[III.] If $P\in E_1$ the assertion follows from Lemma \ref{deg5-57_1_2_3_points} [a).].
 If $P\in E_2\backslash E_1$ the assertion follows from Lemma \ref{deg5-1519_1_2_points}.
 If $P\in  L_{2}\backslash E_2$ the assertion follows from Lemma \ref{deg5-1517_1_3_points}.
 If $P\in  (L_{1}\cup L_1')\backslash E_1$ the assertion follows from Lemma \ref{deg5-3031_32_2_points}.
 If $P\in   L_2'\backslash L_{2}$ the assertion follows from Lemma \ref{deg5-1513_1_2_points} [c).].
 If $P$ is a general point the assertion follows from Lemma \ref{deg5-generalpoint}. 
 \item[IV.] If $P\in L_{13}$ the assertion follows from Lemma \ref{deg5-1523_3_4_points}.
 If $P\in E_1\backslash L_{13}$ the assertion follows from Lemma \ref{deg5-57_1_2_3_points} [b).].
 If $P\in  E_2\backslash E_1$ the assertion follows from Lemma \ref{deg5-1519_1_2_points}.
 If $P\in  L_{2}\backslash E_2$ the assertion follows from Lemma \ref{deg5-1517_1_3_points}.
 If $P\in  E_3\backslash L_{13}$ the assertion follows from Lemma \ref{deg5-1517_1_2_points}.
 If $P\in   L_2'\backslash L_{2}$ the assertion follows from Lemma \ref{deg5-1513_1_2_points} [c).].
 If $P$ is a general point the assertion follows from Lemma \ref{deg5-generalpoint}. 
 \item[V.]  If $P\in E_2$ the assertion follows from Lemma \ref{deg5-59_1_2_4_points}.
 If $P\in E_1\backslash E_2$ the assertion follows from Lemma \ref{deg5-3043_1_32_3_points}.
 If $P\in E_3\backslash E_2$ the assertion follows from Lemma \ref{deg5-1013_32_2_points}.
 If $P\in L_2\backslash E_2$ the assertion follows from Lemma \ref{deg5-1519_2_3_points} [b)].
 If $P\in L_1\backslash E_1$ the assertion follows from Lemma \ref{deg5-1516_2_points}.
 If $P$ is a general point the assertion follows from Lemma \ref{deg5-generalpoint}. 
 \item[VI.] If $P\in E_3$ the assertion follows from Lemma \ref{deg5-37_1_6_points}.
 If $P\in E_2\backslash E_3$ the assertion follows from Lemma \ref{deg5-611_32_4_points}.
 If $P\in L_3\backslash E_3$ the assertion follows from Lemma \ref{deg5-35_4_points}.
 If $P\in  E_4\backslash E_3$ the assertion follows from Lemma \ref{deg5-913_43_3_points}.
 If $P\in E_1 \backslash E_2$ the assertion follows from Lemma \ref{deg5-34_2_points}.
 If $P$ is a general point the assertion follows from Lemma \ref{deg5-generalpoint}. 
 \end{itemize}
 
\end{proof}

\section{Du Val del Pezzo Surfaces of Degree $4$}
\noindent In \cite[Lemma 2.12]{Fano21} it was proven that $\delta(X)=\frac{4}{3}$ when $X$ is a smooth del Pezzo surface of degree $4$. In this chapter, we compute  $\delta$-invariants  of singular Du Val del Pezzo surfaces of degree $4$.  
\begin{maintheorem*}
Let $X$ be a singular Du Val  del Pezzo surface of degree $4$. Then the  $\delta$-invariant of $X$ is uniquely determined by the degree of $X$, the number of lines on $X$, and the type of singularities on $X$ which is given in the following table:
\\{
 \begin{minipage}{5.5cm}
 \renewcommand{\arraystretch}{1.1}
  \begin{longtable}{ | c | c | c | c | }
   \hline
   $K_X^2$ & $\#$ lines & $\mathrm{Sing}(X)$ & $\delta(X)$\\
  \hline\hline
\endhead 
 $4$ & $12$ & $\DA_1$ & $1$\\
\hline
  $4$ & $9$ & $2\DA_1$ & $1$\\
\hline
  $4$ & $8$ & $2\DA_1$ & $1$\\
\hline
 $4$ & $6$ & $3\DA_1$ & $1$\\
\hline
 $4$ & $4$ & $4\DA_1$ & $1$\\
 \hline
  \end{longtable}
  \end{minipage}
  \begin{minipage}{6.5cm}
 \renewcommand{\arraystretch}{1.1}
    \begin{longtable}{ | c | c | c | c | }
   \hline
   $K_X^2$ & $\#$ lines & $\mathrm{Sing}(X)$ & $\delta(X)$\\
  \hline\hline
\endhead 
 $4$ & $8$ & $\DA_2$ & $\frac{6}{7}$\\
\hline
 $4$ & $6$ & $\DA_2+\DA_1$ & $\frac{6}{7}$\\
\hline
$4$ & $4$ & $\DA_2+2\DA_1$ & $\frac{6}{7}$\\
\hline 
$4$ & $5$ & $\DA_3$ & $\frac{2}{3}$\\
\hline 
$4$ & $4$ & $\DA_3$ & $\frac{3}{4}$\\
\hline
  \end{longtable}
  \end{minipage}
    \begin{minipage}{5.5cm}
 \renewcommand{\arraystretch}{1.1}
    \begin{longtable}{ | c | c | c | c | }
   \hline
   $K_X^2$ & $\#$ lines & $\mathrm{Sing}(X)$ & $\delta(X)$\\
  \hline\hline
\endhead 
  $4$ & $3$ & $\DA_3+\DA_1$ & $\frac{3}{4}$\\
\hline
 $4$ & $2$ & $\DA_3+2\DA_1$ & $\frac{3}{4}$\\
\hline
 $4$ & $3$ & $\DA_4$ & $\frac{6}{11}$\\
\hline
$4$ & $2$ & $\mathbb{D}_4$ & $\frac{1}{2}$\\
\hline
$4$ & $1$ & $\mathbb{D}_5$ & $ \frac{3}{8}$\\
\hline
  \end{longtable}
  \end{minipage}}
   \end{maintheorem*}
  
\subsection{General results for degree $4$}
Let $X$ be a del Pezzo surface of degree $4$ with at most Du Val singularities, $S$ be a minimal resolution of $X$ and $P$ is a point on $S$.  Then: 
\begin{lemma}\label{deg4-generalpoint}
 If $P$ is a general point on $S$. Consider the blowup $\sigma:\widetilde{S}\to S$ of $S$ at $P$ with the exceptional divisor $E_P$.  Such that $\tau(E_P)=2$ and the Zariski decomposition of the divisor $\sigma^*(-K_{S})-vE_P$ is given by:
$P(v)=\sigma^*(-K_{S})-vE_P$ and $N(v)=0\text{ if }v\in[0,2]$ Moreover,
$P(v)^2=(2-v)(2+v)$ and $P(v)\cdot E_P=v \text{ if }v\in[0,2]$.
In this case $\delta_P(S)=\frac{3}{2}.$
\end{lemma}
\begin{proof}
The Zariski Decomposition follows from $\sigma^*(-K_S)-vE_P\sim_{\DR}  (2 - v)E_P+Q+L$ where $Q$ and $L$ are $(-1)$-curves which are strict transforms of a conic and a line on $\DP^2$ respectively. We have
$S_S(E_P)=\frac{4}{3}$.
Thus, $\delta_P(S)\le \frac{2}{4/3}=\frac{3}{2}$ Moreover if $O\in E_P$ we have
    $h(v)=\frac{v^2}{2}\text{ if }v\in[0,2].$
So 
$S(W_{\bullet,\bullet}^{E_P};O) \le \frac{2}{3}$
Thus, $\delta_P(S)=\frac{3}{2}$.
\end{proof}
\begin{lemma}\label{deg4-conic}
Suppose  $P\in Q$ where $Q$ is a $(-1)$-curve on $S$ such that it does not intersect $(-2)$-curves. Let $L_P$ be a $(1)$-curve which intersects with $Q$ at this point with multiplicity $2$ (one can see that such curve exists by considering explicit models). Consider the blowup of $S$ at $P$ with the exceptional divisor $E_P$. After that we blow up  the intersection of strict transforms $Q$, $L_P$, and $E_P$ with the exceptional divisor $E$, and contract a $(-2)$-curve. We call the resulting surface $\overline{S}$. Let $\overline{L}_P$, $\overline{Q}$, $\overline{E}$, and $\overline{P}$ be the images of strict transforms of $L_P$, $Q$, $E$ and $E_P$ respectively. Note that $\overline{P}$ is a point. We denote this weighted blowup by $\sigma:\overline{S}\to S$.
\begin{figure}[h!]
    \centering
  \includegraphics[width=15cm]{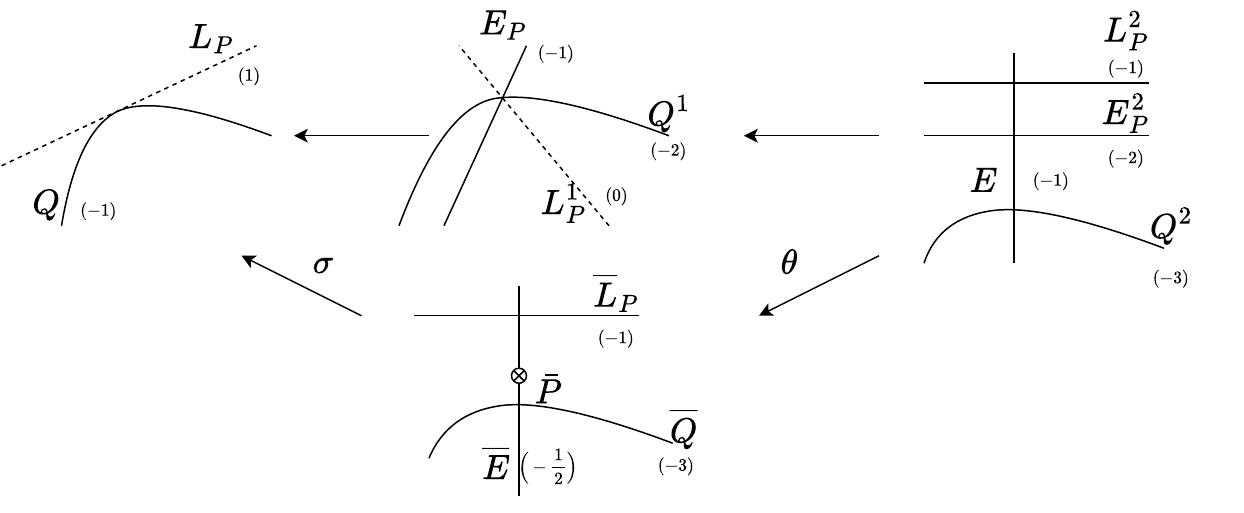}
    \caption{$(-K_S)^2=4$,  a general point on a $(-1)$-curve}
\end{figure}
\par The intersections are given by
\begin{table}[h!]
    \centering
    \begin{tabular}{ |c||c|c|c|  }
 \hline
& $\overline{Q}$ & $\overline{E}$ & $\overline{L}_P$  \\ 
 \hline\hline
 $\overline{Q}$ & \cellcolor[gray]{0.7} $-3$ & \cellcolor[gray]{0.9}$1$ & $0$  \\ 
 \hline
$\overline{E}$  & \cellcolor[gray]{0.9} $1$ & \cellcolor[gray]{0.7} $-\frac{1}{2}$ & \cellcolor[gray]{0.9} $1$  \\ 
 \hline
  $\overline{L}_P$ & $0$  & \cellcolor[gray]{0.9} $1$ & \cellcolor[gray]{0.7} $-1$ \\ 
 \hline
\end{tabular}
    \caption{$(-K_S)^2=4$, a general point on a $(-1)$-curve}
\end{table}

Then $\tau(\overline{E})=4$ and the Zariski decomposition of the divisor $\sigma^*(-K_{S})-v\overline{E}$  is given by:
{\allowdisplaybreaks\begin{align*}
&&P(v)=\begin{cases}
\sigma^*(-K_{S})-v\overline{E}\text{ if }v\in[0,1],\\
\sigma^*(-K_{S})-v\overline{E}-\frac{(v-1)}{3}\overline{Q}\text{ if }v\in[1,3],\\
\sigma^*(-K_{S})-v\overline{E}-\frac{(v-1)}{3}\overline{Q}-(v-3)\overline{L}_P\text{ if }v\in[3,4].
\end{cases}
\\&&N(v)=\begin{cases}
0\text{ if }v\in[0,1],\\
\frac{(v-1)}{3}\overline{Q}\text{ if }v\in[1,3],\\
\frac{(v-1)}{3}\overline{Q}+(v-3)\overline{L}_P\text{ if }v\in[3,4].
\end{cases}
\end{align*}}
Moreover
$$P(v)^2=\begin{cases}
4-\frac{v^2}{2}\text{ if }v\in[0,1],\\
\frac{26-4v-v^2}{6}\text{ if }v\in[1,3],\\
\frac{5(4-v)^2}{6}\text{ if }v\in[3,4].
\end{cases}
P(v)\cdot \overline{E}=\begin{cases}
\frac{v}{2}\text{ if }v\in[0,1],\\
\frac{2+v}{6}\text{ if }v\in[1,3],\\
\frac{5(4-v)}{6}\text{ if }v\in[3,4].
\end{cases}$$
In this case $\delta_P(S)=\frac{18}{13}$.
\end{lemma}
\begin{proof}
The Zariski Decomposition follows from $\sigma^*(-K_S)-v\overline{E}\sim \overline{L}_P+\overline{Q}+(4-v)\overline{E}$. We have
$S_S(\overline{E})=\frac{13}{6}$.
Thus, since  $A_S(\overline{E})=3$ then $\delta_P(S)\le \frac{3}{13/6}=\frac{18}{13}$.
Moreover, if $O\in\overline{E}\backslash(\overline{Q}\cup \overline{L}_P)$ or  $O \in \overline{E}\cap \overline{Q}$ or if  $O\in \overline{E}\cap \overline{L}_P$:
$$
h(v)\le \begin{cases}
\frac{v^2}{8}\text{ if }v\in[0,1],\\
 \frac{(v + 2)^2}{72}\text{ if }v\in[1,3],\\
\frac{25 (v - 4)^2}{72}\text{ if }v\in[3,4].
\end{cases}\text{ or }h(v)\le \begin{cases}
\frac{v^2}{8}\text{ if }v\in[0,1],\\
\frac{(v + 2) (5 v - 2)}{72}\text{ if }v\in[1,3],\\
\frac{5 (4 - v) (16 - v)}{72}\text{ if }v\in[3,4].
\end{cases}$$$$\text{ or }h(v)\le \begin{cases}
\frac{v^2}{8}\text{ if }v\in[0,1],\\
\frac{ (v + 2)^2}{72}\text{ if }v\in[1,3],\\
\frac{5 (4 - v) (7 v - 16)}{72}\text{ if }v\in[3,4].
\end{cases}$$
Thus if $O\in\overline{E}\backslash(\overline{Q}\cup \overline{L}_P)$ then $S(W_{\bullet,\bullet}^{\overline{E}};O)=\frac{11}{36}$, if $O \in \overline{E}\cap \overline{Q}$ then $S(W_{\bullet,\bullet}^{\overline{E}};O)=\frac{17}{24}$, if $O\in E_P\cap L_P$ then $S(W_{\bullet,\bullet}^{\overline{E}};O)=\frac{3}{8}$.
Now, to get a~lower bound for $\delta_P(S)$, we use Corollary~\ref{estimation2} that gives
$$
\delta_P(S)\geqslant\mathrm{min}\Bigg\{\frac{18}{13},\inf_{O\in\overline{E}}\frac{A_{\overline{E},\Delta_{\overline{E}}}(O)}{S\big(W^{\overline{E}}_{\bullet,\bullet};O\big)}\Bigg\},
$$
where $\Delta_{\overline{E}}=\frac{1}{2}\theta(\overline{E})$.
So $\delta_P(S)\ge\min\Big\{\frac{18}{13},\frac{36}{11},\frac{36}{22},\frac{24}{17},\frac{8}{3}\Big\}=\frac{18}{13}$ thus $\delta_P(S)=\frac{18}{13}$.
\end{proof}
\begin{lemma}\label{deg4-twolines}
Suppose $A_1$ and $A_2$ are $(-1)$-curves on $S$ which intersect at point $P$.  Consider the blowup $\sigma:\widetilde{S}\to S$ of $S$ at $P$ with the exceptional divisor $E_P$. Suppose $\MA_1$ and $\MA_2$ are strict transforms of  $A_1$ and $A_2$ on $\widetilde{S}$ and $\MA_3$ is another $(-1)$-curve on $S$ such that the intersections are given by
\begin{table}[]
    \centering
    \begin{tabular}{ |c||c|c|c| c | } 
 \hline
& $E_P$ & $\MA_1$ & $\MA_2$ &$\MA_3$ \\ 
 \hline\hline
 $E_P$&\cellcolor[gray]{0.7}  $-1$ & \cellcolor[gray]{0.9} $1$ & \cellcolor[gray]{0.9} $1$ & \cellcolor[gray]{0.9} $1$ \\ 
 \hline
 $\MA_1$  & \cellcolor[gray]{0.9} $1$ &\cellcolor[gray]{0.7} $-2$ & $0$ & $0$ \\ 
 \hline
  $\MA_2$ & \cellcolor[gray]{0.9} $1$  & $0$ &\cellcolor[gray]{0.7} $-2$ & $0$ \\ 
 \hline
  $\MA_3$ &\cellcolor[gray]{0.9} $1$ & $0$ & $0$ & \cellcolor[gray]{0.7}$-1$ \\ 
 \hline
\end{tabular}
\caption{$(-K_S)^2=4$, the intersection of two $(-1)$-curves}
\end{table}
Then $\tau(E_P)=3$ and the Zariski decomposition of the divisor $\sigma^*(-K_{S})-vE_P$  is given by:
\begin{align*}
 &&P(v)=\begin{cases}
\sigma^*(-K_S)-vE_P\text{ if }v\in[0,1],\\
\sigma^*(-K_S)-vE_P-\frac{(v-1)}{2}(\MA_1+\MA_2)\text{ if }v\in[1,2],\\
\sigma^*(-K_S)-vE_P-\frac{(v-1)}{2}(\MA_1+\MA_2)-(v-2)\MA_3\text{ if }v\in[2,3].
\end{cases}\\&& N(v)=\begin{cases}
0\text{ if }v\in[0,1],\\
\frac{(v-1)}{2}(\MA_1+\MA_2)\text{ if }v\in[1,2],\\
\frac{(v-1)}{2}(\MA_1+\MA_2)+(v-2)\MA_3\text{ if }v\in[2,3].
\end{cases}
\end{align*}
with
$$P(v)^2=\begin{cases}
(2-v)(2+v)\text{ if }v\in[0,1],\\
5-2v\text{ if }v\in[1,2],\\
(3-v)^2\text{ if }v\in[2,3].
\end{cases}
P(v)\cdot E_P=\begin{cases}
v\text{ if }v\in[0,1],\\
1\text{ if }v\in[1,2],\\
3-v\text{ if }v\in[2,3].
\end{cases}$$
In this case $\delta_P(S)=\frac{4}{3}$.
\end{lemma}
\begin{proof}
The Zariski Decomposition follows from $\sigma^*(-K_S)-vE_P\sim_{\DR}  (3 - v)E_P+\MA_1+\MA_2+\MA_3$.
We have
$S_S(E_P)=\frac{3}{2}$.
Thus $\delta_P(S)\le \frac{2}{3/2}=\frac{4}{3}$.
Moreover, if $O\in E_P\backslash \MA_3$ or if $O\in E_P\cap \MA_3$ then:
$$h(v)\le \begin{cases}
\frac{v^2}{2}\text{ if }v\in[0,1],\\
\frac{v}{2}\text{ if }v\in[1,2],\\
3-v\text{ if }v\in[2,3].
\end{cases}
\text{ or }
h(v)\le\begin{cases}
\frac{v^2}{2}\text{ if }v\in[0,1],\\
\frac{1}{2}\text{ if }v\in[1,2],\\
\frac{(3 - v) (v - 1)}{2}\text{ if }v\in[2,3].
\end{cases}
$$
So 
$S(W_{\bullet,\bullet}^{E_P};O)\le \frac{17}{24}\le\frac{3}{4}$
or
$S(W_{\bullet,\bullet}^{E_P};O)\le\frac{1}{2}\le\frac{3}{4}$
Thus, $\delta_P(S)=\frac{4}{3}$ for $A_1\cap A_2$.
\end{proof}
\begin{lemma}\label{deg4-nearA1points}
 If $P$ belongs to a $(-1)$-curve $A$ and there exist $(-1)$-curves and $(-2)$-curves   which form the following dual graph:
\begin{figure}[h!]
    \centering
 \includegraphics[width=10cm]{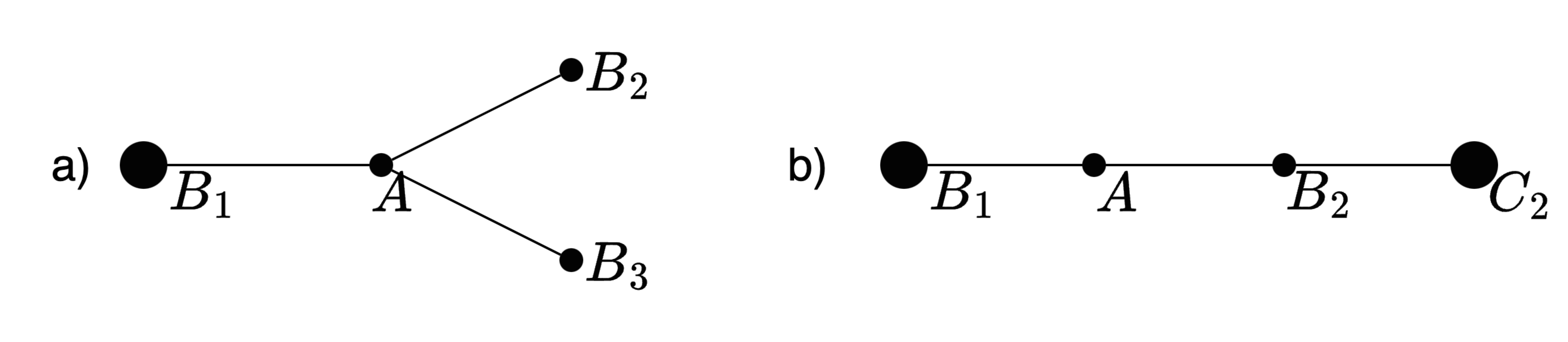}
    \caption{Dual graph: $(-K_S)^2=4$ and $\delta_P(S)=\frac{6}{5}$}
\end{figure}
\par   Then $\tau(A)=2$ and the Zariski Decomposition of the divisor $-K_S-vA$ is given by:
   {\allowdisplaybreaks \begin{align*}
  &  {\text{\bf a). }}&P(v)=\begin{cases}
-K_S-vA-\frac{v}{2}B_1\text{ if }v\in[0,1],\\
-K_S-vA-\frac{v}{2}B_1-(v-1)(B_2+B_3)\text{ if }v\in[1,2].
\end{cases}\\&&N(v)=\begin{cases}
\frac{v}{2}B_1\text{ if }v\in[0,1],\\
\frac{v}{2}B_1+(v-1)(B_2+B_3)\text{ if }v\in[1,2].
\end{cases}\\
&{\text{\bf b). }}&P(v)=\begin{cases}
-K_S-vA-\frac{v}{2}B_1\text{ if }v\in[0,1],\\
-K_S-vA-\frac{v}{2}B_1-(v-1)(2B_2+C_2)\text{ if }v\in[1,2].
\end{cases}\\&&N(v)=\begin{cases}
\frac{v}{2}B_1\text{ if }v\in[0,1],\\
\frac{v}{2}B_1+(v-1)(2B_2+C_2)\text{ if }v\in[1,2].
\end{cases}
\end{align*}}
Moreover:
$$(P(v))^2=\begin{cases}
4-2v-\frac{v^2}{2}\text{ if }v\in[0,1],\\
\frac{3(2-v)^2}{2}\text{ if }v\in[1,2].
\end{cases}
P(v)\cdot A=\begin{cases}
 1+\frac{v}{2}\text{ if }v\in[0,1],\\
3-\frac{3v}{2}\text{ if }v\in[1,2].
\end{cases}$$
In this case $\delta_P(S)=\frac{6}{5}\text{ if }P\in A\backslash B_1$.
\end{lemma}
\begin{proof}
The Zariski Decomposition in part a). follows from $-K_S-vA\sim_{\DR} (2-v)A+B_1+B_2+B_3$. A similar statement holds in other parts.
We have
$S_{S}(A)=\frac{5}{6}$. Thus, $\delta_P(S)\le \frac{6}{5}$ for $P\in A$. Moreover,
$$h(v)=\begin{cases}
\frac{(v + 2)^2}{8}\text{ if }v\in[0,1],\\
\frac{3 (2 - v) (5 v - 2)}{8}\text{ if }v\in[1,2].
\end{cases}$$
So
$S(W_{\bullet,\bullet}^{A};P)\le \frac{5}{6}$.
We get that $\delta_P(S)=\frac{6}{5}$ for $P\in A\backslash B_1$.
\end{proof}
\begin{lemma}\label{deg4-nearA2points}
 If $P$ belongs to a $(-1)$-curve $A$ and there exist $(-1)$-curves and $(-2)$-curves   which form the following dual graph:
\begin{figure}[h!]
    \centering
  \includegraphics[width=4.2cm]{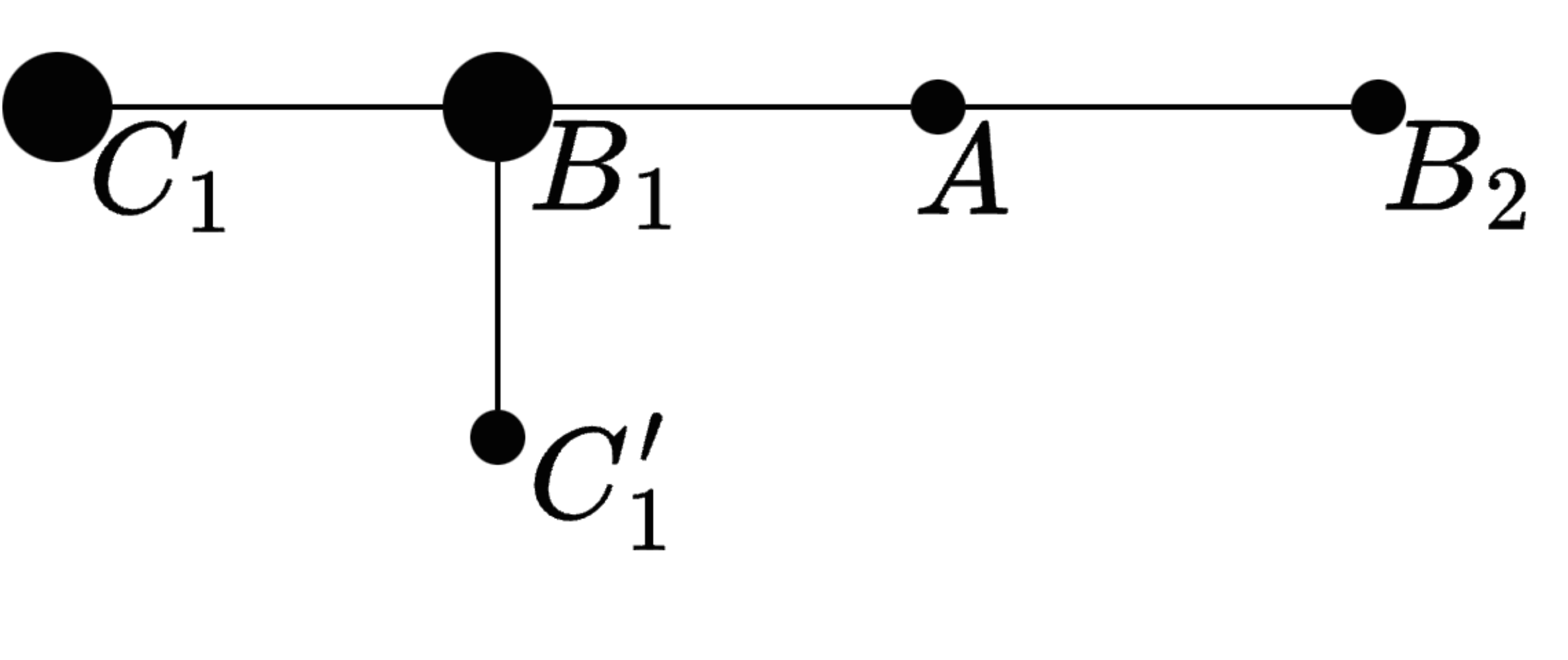}
    \caption{Dual graph: $(-K_S)^2=4$ and $\delta_P(S)=\frac{8}{7}$}
\end{figure}
\par  Then $\tau(A)=2$ and the Zariski Decomposition of the divisor $-K_S-vA$ is given by:
{\allowdisplaybreaks \begin{align*}
&&P(v)=\begin{cases}-K_S-v A-\frac{v}{3}(2B_1+C_1)\text{ if }v\in[0,1],\\
-K_S-v A-\frac{v}{3}(2B_1+C_1)-(v-1)B_2\text{ if }v\in\big[1,\frac{3}{2}\big],\\
-K_S-v A-(v-1)(2B_1+C_1+B_2)-(2v-3)C_1'\text{ if }v\in\big[\frac{3}{2},2\big].
\end{cases}\\
&&N(v)=\begin{cases}\frac{v}{3}(2B_1+C_1)\text{ if }v\in[0,1],\\
\frac{v}{3}(2B_1+C_1)+(v-1)B_2\text{ if }v\in\big[1,\frac{3}{2}\big],\\
(v-1)(2B_1+C_1+B_2)+(2v -3 )C_1'\text{ if }v\in\big[\frac{3}{2},2\big].
\end{cases}
\end{align*}}
Moreover, 
$$(P(v))^2=\begin{cases} 
4- 2v-\frac{v^2}{3} \text{ if }v\in[0,1],\\
5 -4v + \frac{2v^2}{3} \text{ if }v\in\big[1,\frac{3}{2}\big],\\
2(v - 2)^2 \text{ if }v\in\big[\frac{3}{2},2\big].
\end{cases}P(v)\cdot  A=\begin{cases}
1+\frac{v}{3}\text{ if }v\in[0,1],\\
2(1-\frac{v}{3})\text{ if }v\in\big[1,\frac{3}{2}\big],\\
4-2v\text{ if }v\in\big[\frac{3}{2},2\big].
\end{cases}$$
In this case $\delta_P(S)=\frac{8}{7}\text{ if }P\in A\backslash B_1$.
\end{lemma}
\begin{proof}
The Zariski Decomposition follows from $-K_S-vA\sim_{\DR} (2-v)A+2B_1+C_1+C_1'+B_2$.
We have 
$S_S(A)=\frac{7}{8}$. Thus, $\delta_P(S)\le \frac{8}{7}$ for $P\in  A$. Note that for $P\in A\backslash B_1$ we have:
$$h(v)=\begin{cases}
\frac{ (v + 3)^2}{18}\text{ if }v\in[0,1],\\
\frac{4 v (3-v)}{9}\text{ if }v\in\big[1,\frac{3}{2}\big],\\
 2 (2 - v)\text{ if }v\in\big[\frac{3}{2},2\big].
\end{cases}$$
So $S(W_{\bullet,\bullet}^{A};P)\le\frac{17}{24}\le\frac{7}{8}$. Thus, $\delta_P(S)=\frac{8}{7}$ if $P\in A\backslash B_1$. 
\end{proof}
\begin{lemma}\label{deg4-nearA3points}
 If $P$ belongs to a $(-1)$-curve $A$ and there exist $(-1)$-curves and $(-2)$-curves   which form the following dual graph:
\begin{figure}[h!]
    \centering
 \includegraphics[width=4cm]{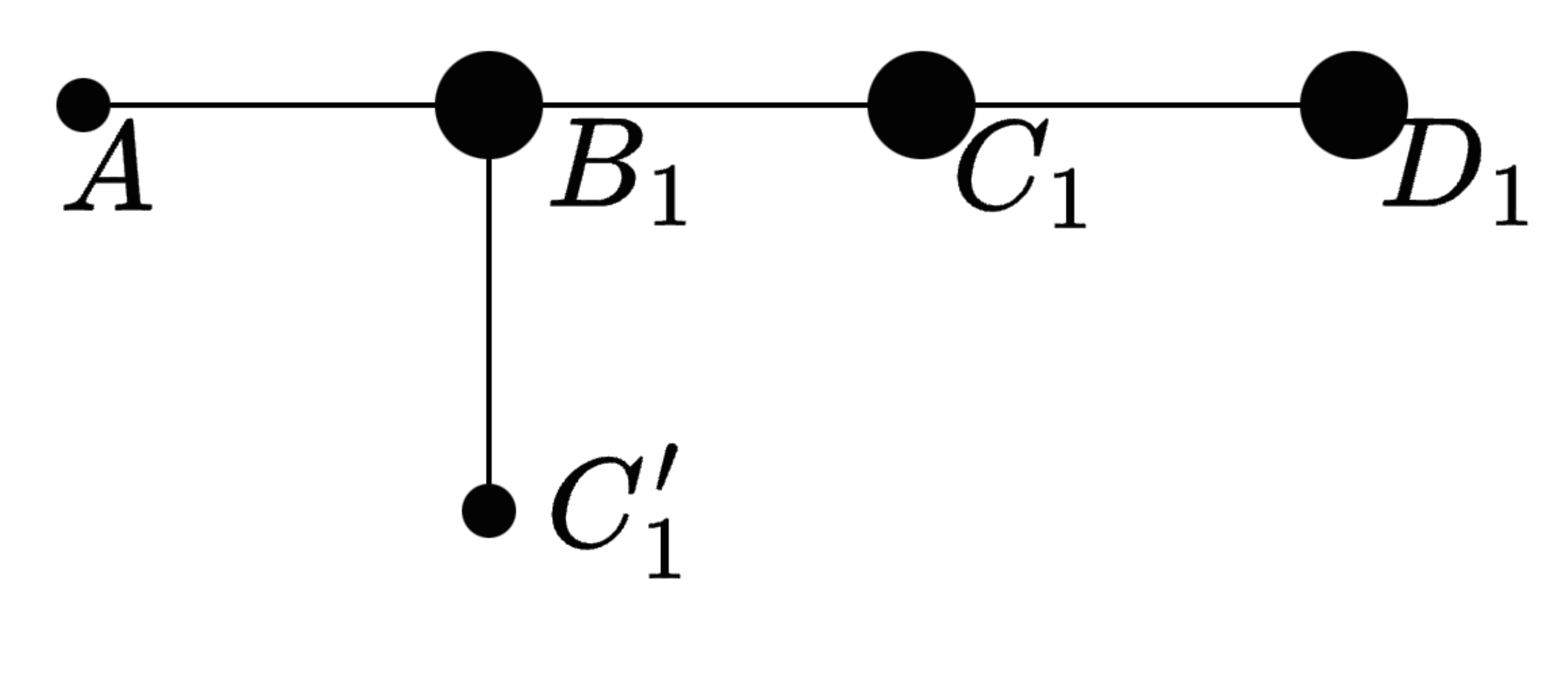}
    \caption{Dual graph: $(-K_S)^2=4$ and $\delta_P(S)=\frac{9}{8}$}
\end{figure}
 \par  Then $\tau(A)=2$ and the Zariski Decomposition of the divisor $-K_S-vA$ is given by:
{{\allowdisplaybreaks \begin{align*}
&&P(v)=\begin{cases}
-K_S-vA-\frac{v}{4}( 3B_1+2C_1+D_1)\text{ if }v\in\big[0,\frac{4}{3}\big],\\
-K_S-vA-(v-1)(3B_1+2C_1+D_1)-(3v-4)C_1'\text{ if }v\in\big[\frac{4}{3},2\big].
\end{cases}\\&&
N(v)=\begin{cases}\frac{v}{4}( 3B_1+2C_1+D_1)\text{ if }v\in\big[0,\frac{4}{3}\big],\\
(v-1)(3B_1+2C_1+D_1)+(3v-4)C_1'\text{ if }v\in\big[\frac{4}{3},2\big].
\end{cases}\end{align*}}}
Moreover, 
$$(P(v))^2=\begin{cases} 4-2v-\frac{v^2}{4} \text{ if }v\in\big[0,\frac{4}{3}\big],\\
2(v - 2)^2\text{ if }v\in\big[\frac{4}{3},2\big].
\end{cases}P(v)\cdot  A=\begin{cases}1+\frac{v}{4} \text{ if } v  \in\big[0,\frac{4}{3}\big],\\
2(2-v) \text{ if }v\in\big[\frac{4}{3},2\big].
\end{cases}$$
In this case 
$\delta_P(S)=\frac{9}{8}\text{ if }P\in A\backslash B_1$.
\end{lemma}
\begin{proof}
The Zariski Decomposition follows from $-K_S-vA\sim_{\DR} (2-v)A+3B_1+2C_1+D_1+2C_1'$. Thus,
$S_S(A)=\frac{8}{9}$.
Thus, $\delta_P(S)\le \frac{9}{8}$ for $P\in  A$. Note that we have that if $P\in  A\backslash B_1$ then
$$h(v)=\begin{cases}
\frac{ (v + 4)^2}{32}\text{ if }v\in\big[0,\frac{4}{3}\big],\\
2(2-v)^2\text{ if }v\in\big[\frac{4}{3},2\big].
\end{cases}$$
So for $S(W_{\bullet,\bullet}^{A};P)\le \frac{5}{9}\le\frac{8}{9}$.
Thus, $\delta_P(S)=\frac{9}{8}$ if $P\in A\backslash B_1$. 
\end{proof}
\begin{lemma}\label{deg4-1211_1_2_points}
 If $P$ belongs to a $(-1)$-curve $A$ and there exist $(-1)$-curves and $(-2)$-curves   which form the following dual graph:
\begin{figure}[h!]
    \centering
\includegraphics[width=5cm]{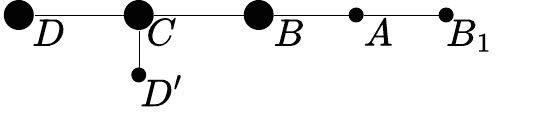}
    \caption{Dual graph: $(-K_S)^2=4$ and $\delta_P(S)=\frac{12}{11}$}
\end{figure}
 \par    Then $\tau(A)=2$ and the Zariski Decomposition of the divisor $-K_S-vA$ is given by:
{\allowdisplaybreaks \begin{align*}
&&P(v)=\begin{cases}
-K_S-vA-\frac{v}{4}( 3B+2C+D)\text{ if }v\in[0,1],\\
-K_S-vA-\frac{v}{4}( 3B+2C+D)-(v-1)B_1\text{ if }v\in[1,2].
\end{cases}\\&&N(v)=\begin{cases}\frac{v}{4}( 3B+2C+D)\text{ if }v\in[0,1],\\
\frac{v}{4}( 3B+2C+D)+(v-1)B_1\text{ if }v\in[1,2].
\end{cases}
\end{align*}}
Moreover, 
$$(P(v))^2=\begin{cases}4-2v-\frac{v^2}{4} \text{ if }v\in[0,1],\\
\frac{(2-v)(10-3v)}{4}\text{ if }v\in[1,2].
\end{cases}P(v)\cdot  A=\begin{cases}2v\text{ if }    \in[0,1]\\
2-\frac{3v}{4} \text{ if }v\in[1,2].
\end{cases}$$
In this case $\delta_P(S)=\frac{12}{11}\text{ if }P\in A\backslash B$.
\end{lemma}
\begin{proof}
The Zariski Decomposition follows from 
$-K_S-vA\sim_{\DR} (2-v)A+2B+2C+D+D'+B_1$.
We have
$S_S(A)=\frac{11}{12}$.
Thus, $\delta_P(S)\le \frac{12}{11}$ for $P\in  A$. Note that for $P\in A\backslash B$ we have:
$$h(v)\le
\begin{cases}
\frac{(v + 4)^2}{32}\text{ if }v\in[0,1],\\
\frac{5 (8 - 3 v) v}{32}\text{ if }v\in[1,2].
\end{cases}$$
So $S(W_{\bullet,\bullet}^{A};P)\le \frac{17}{24}\le\frac{11}{12}$.
Thus, $\delta_P(S)=\frac{12}{11}$ if $P\in  A\backslash B$.
\end{proof}
\begin{lemma}\label{deg4-2423_1_52_points}
 If $P$ belongs to a $(-1)$-curve $A$ and there exist $(-1)$-curves and $(-2)$-curves   which form the following dual graph:
\begin{figure}[h!]
    \centering
\includegraphics[width=7cm]{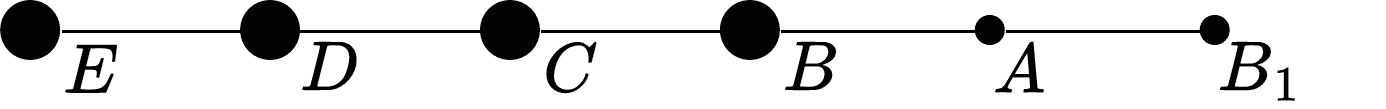}
    \caption{Dual graph: $(-K_S)^2=4$ and $\delta_P(S)=\frac{24}{23}$}
\end{figure}
 \par   Then $\tau(A)=\frac{5}{2}$ and the Zariski Decomposition of the divisor $-K_S-vA$ is given by:
 {\allowdisplaybreaks \begin{align*}
&&P(v)=\begin{cases}
-K_S-vA-\frac{v}{5}(4B+3C+2D+E)\text{ if }v\in[0,1],\\
-K_S-vA-\frac{v}{5}(4B+3C+2D+E)-(v-1)B_1\text{ if }v\in\big[1,\frac{5}{2}\big].
\end{cases}\\&&N(v)=\begin{cases}\frac{v}{5}(4B+3C+2D+E)\text{ if }v\in[0,1],\\
\frac{v}{5}(4B+3C+2D+E)+(v-1)B_1\text{ if }v\in\big[1,\frac{5}{2}\big].
\end{cases}
\end{align*}}
Moreover, 
$$(P(v))^2=\begin{cases}
4 -2v -\frac{v^2}{5}  \text{ if }v\in[0,1],\\
\frac{(2v-5)^2}{5}\text{ if }v\in\big[1,\frac{5}{2}\big].
\end{cases}P(v)\cdot  A=
\begin{cases}
1+\frac{v}{5}\text{ if } v  \in[0,1],\\
2-\frac{4v}{5}\text{ if }v\in\big[1,\frac{5}{2}\big].
\end{cases}$$
In this case $\delta_P(S)=\frac{24}{23}\text{ if }P\in A\backslash B$.
\end{lemma}

\begin{proof}
The Zariski Decomposition  follows from 
 $$-K_S-vA\sim_{\DR} \Big(\frac{5}{2}-v\Big)A+\frac{1}{2}\Big(4B+3C+2D+E+3B_1\Big).$$
We have
$S_S(A)=\frac{23}{24}$.
Thus, $\delta_P(S)\le \frac{24}{23}$ for $P\in  A\backslash B$. Note that we have:
$$h(v)\le \begin{cases}\frac{ (v + 5)^2}{50}\text{ if }v\in[0,1],\\
\frac{ 6 v (5 - 2 v)}{25}\text{ if }v\in\big[1,\frac{5}{2}\big].
\end{cases}$$
So $S(W_{\bullet,\bullet}^{A};P)\le \frac{73}{120}\le\frac{23}{24}$.
Thus, $\delta_P(S)=\frac{24}{23}$ if $P\in A\backslash B$.
\end{proof}
\begin{lemma}\label{deg4-A1points}
 If $P$ belongs to a $(-2)$-curve $A$ and there exist $(-1)$-curves and $(-2)$-curves   which form the following dual graph:
\begin{figure}[h!]
    \centering
  \includegraphics[width=9cm]{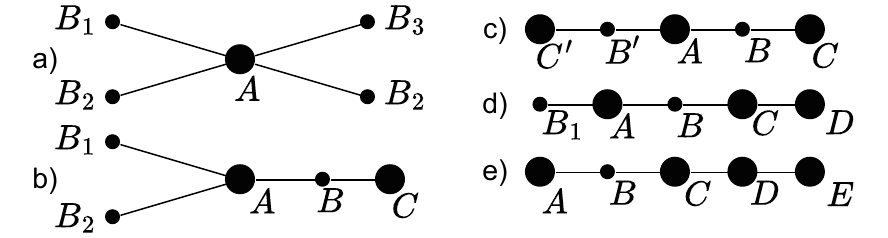}
    \caption{Dual graph: $(-K_S)^2=4$ and $\delta_P(S)=1$ with $-K_S-vA$ nef on $[0,1]$}
\end{figure}
\par Then $\tau(A)=2$ and the Zariski Decomposition of the divisor $-K_S-vA$ is given by:
   {\allowdisplaybreaks \begin{align*}
&{\text{\bf a). }} & P(v)=\begin{cases}-K_S-vA\text{ if }v\in[0,1],\\
-K_S-vA-(v-1)(B_1+B_2+B_3+B_4)\text{ if }v\in[1,2].
\end{cases}\\&&N(v)=\begin{cases}0\text{ if }v\in[0,1],\\
(v-1)(B_1+B_2+B_3+B_4)\text{ if }v\in[1,2].
\end{cases}\\
&{\text{\bf b). }} &P(v)=\begin{cases}-K_S-vA\text{ if }v\in[0,1],\\
-K_S-vA-(v-1)(2B+C+B_1+B_2)\text{ if }v\in[1,2].
\end{cases}\\&&N(v)=\begin{cases}0\text{ if }v\in[0,1],\\
(v-1)(2B+C+B_1+B_2)\text{ if }v\in[1,2].
\end{cases}\\
&{\text{\bf c). }} &P(v)=\begin{cases}-K_S-vA\text{ if }v\in[0,1],\\
-K_S-vA-(v-1)(2B+C+2B'+C')\text{ if }v\in[1,2].
\end{cases}\\&&N(v)=\begin{cases}0\text{ if }v\in[0,1],\\
(v-1)(2B+C+2B'+C')\text{ if }v\in[1,2].
\end{cases}\\
&{\text{\bf d). }}& P(v)=\begin{cases}-K_S-vA\text{ if }v\in[0,1],\\
-K_S-vA-(v-1)(3B+2C+D+B_1)\text{ if  }v\in[1,2]
\end{cases}\\&&N(v)=\begin{cases}0\text{ if  }v\in[0,1]\\
(v-1)(3B+2C+D+B_1)\text{ if }v\in[1,2].
\end{cases}\\
&{\text{\bf e). }} & P(v)=\begin{cases}-K_S-vA\text{ if  }v\in[0,1],\\
-K_S-vA-(v-1)(4B+3C+2D+E)\text{ if  }v\in[1,2].
\end{cases}\\&&N(v)=\begin{cases}0\text{ if  }v\in[0,1],\\
(v-1)(4B+3C+2D+E)\text{ if  }v\in[1,2].
\end{cases}
\end{align*}}
Moreover:
$$(P(v))^2=\begin{cases}4-2v^2\text{ if }v\in[0,1],\\
2(2-v)^2\text{ if }v\in[1,2].
\end{cases}P(v)\cdot A=\begin{cases}2v\text{ if }v\in[0,1],\\
2(2-v)\text{ if }v\in[1,2].
\end{cases}$$
In this case $\delta_P(S)=1$ if $P\in A\backslash B$.
\end{lemma}
\begin{proof}
The Zariski Decomposition in part a). follows from $-K_S-vA\sim_{\DR} (2-v)A+B_1+B_2+B_3+B_4$. A similar statement holds in other parts.  Consider a point $P$ described above. We have
$S_S(A)=1$.
Thus, $\delta_P(S)\le 1$ for $P\in A$. Note that we have:
$$h(v)\le\begin{cases} 2v^2\text{ if }v\in[0,1],\\
2v(2-v)\text{ if }v\in[1,2].
\end{cases}$$
So 
$S(W_{\bullet,\bullet}^{A};P)\le1$.
Thus, $\delta_P(S)=1$.
\end{proof}
\begin{lemma}\label{deg4-near2A1points}
 If $P$ belongs to a $(-1)$-curve $A$ and there exist $(-1)$-curves and $(-2)$-curves   which form the following dual graph:
\begin{figure}[h!]
    \centering
  \includegraphics[width=11cm]{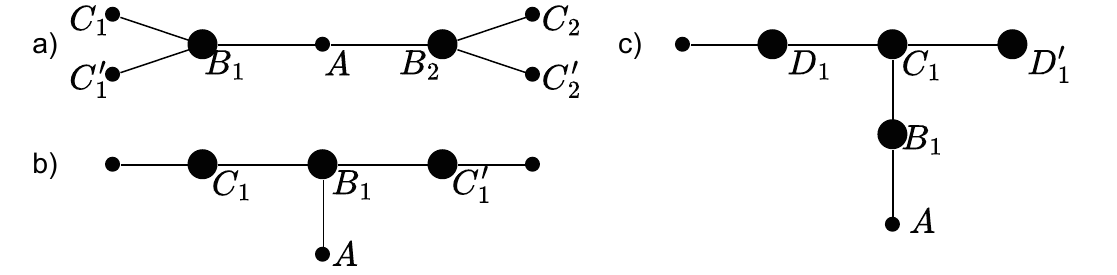}
    \caption{Dual graph: $(-K_S)^2=4$ and $\delta_P(S)=1$ with $-K_S-vA$ nef on $[0,2]$}
\end{figure}
  \par Then $\tau(A)=2$ and the Zariski Decomposition of the divisor $-K_S-vA$ is given by:
{\allowdisplaybreaks \begin{align*}
&{\text{\bf a). }}&P(v)=-K_S-vA-\frac{v}{2}(B_1+B_2) \text{ if }v\in[0,2].\\&&N(v)=\frac{v}{2}(B_1+B_2)\text{ if }v\in[0,2].\\
&{\text{\bf b). }}&P(v)=-K_S-v A-\frac{v}{2}( 2B_1+C_1+C_1') \text{ if }v\in[0,2].\\&& N(v)=\frac{v}{2}( 2B_1+C_1+C_1')\text{ if }v\in[0,2].\\
&{\text{\bf c). }}&P(v)=-K_S-v A-\frac{v}{2}( 2B_1+2C_1+D_1+D_1')\text{ if }v\in[0,2].\\&& N(v)=\frac{v}{2}( 2B_1+2C_1+D_1+D_1')\text{ if }v\in[0,2].
\end{align*}}
Moreover, 
$$(P(v))^2=4 - 2v \text{ and }
P(v)\cdot A=1\text{ if }v\in[0,2].$$
In this case $\delta_P(S)=1$ if $P\in A\backslash (B_1\cup B_2)$.
\end{lemma}
\begin{proof}
The Zariski Decomposition in part a). follows from $-K_S-vA\sim_{\DR} (2-v)A+\frac{1}{2}\Big(3B_1+C_1+C_1'+3B_2+C_2+C_2'\Big)$. A similar statement holds in other parts. We have
$S_{S}(A)= 1$.
Thus, $\delta_P(S)\le 1$ for $P\in A$. Moreover, if $P\in A\backslash (B_1\cup B_2)$ we have
$h(v)=1/2\text{ if }v\in[0,2].$
So
$S(W_{\bullet,\bullet}^{A};P)\le \frac{1}{2}\le 1 $.
We get that $\delta_P(S)= 1$ for $P\in A\backslash (B_1\cup B_2)$.
\end{proof}
\begin{lemma}\label{deg4-A2points}
 If $P$ belongs to a $(-2)$-curve $A$ and there exist $(-1)$-curves and $(-2)$-curves   which form the following dual graph:
\begin{figure}[h!]
    \centering
  \includegraphics[width=10.4cm]{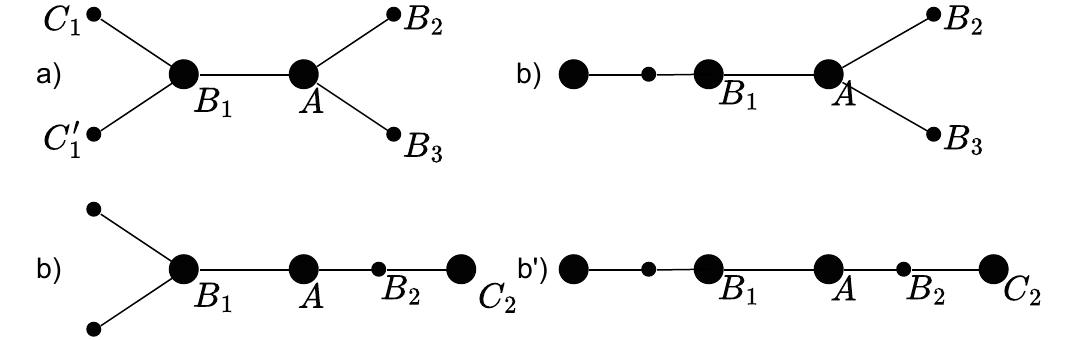}
    \caption{Dual graph: $(-K_S)^2=4$ and $\delta_P(S)=\frac{6}{7}$ with $\tau(A)=2$}
\end{figure}
\par   Then $\tau(A)=2$ and the Zariski Decomposition of the divisor $-K_S-vA$ is given by:
{\allowdisplaybreaks \begin{align*}
&{\text{\bf a). }}&P(v)=\begin{cases}-K_S-v A-\frac{v}{2}B_1\text{ if }v\in[0,1],\\
-K_S-v A-\frac{v}{2}B_1-(v-1)(B_2+B_3)\text{ if }v\in[1,2].
\end{cases}\\&&N(v)=\begin{cases}\frac{v}{2}B_1\text{ if }v\in[0,1],\\
\frac{v}{2}B_1+(v-1)(B_2+B_3)\text{ if }v\in[1,2].
\end{cases}\\
&{\text{\bf b). }}&P(v)=\begin{cases}-K_S-v A-\frac{v}{2}B_1\text{ if }v\in[0,1],\\
-K_S-v A-\frac{v}{2}B_1-(v-1)(2B_2+C_2)\text{ if }v\in[1,2].
\end{cases}\\&&N(v)=\begin{cases}\frac{v}{2}B_1\text{ if }v\in[0,1],\\
\frac{v}{2}B_1+(v-1)(2B_2+C_2)\text{ if }v\in[1,2].
\end{cases}
\end{align*}}
Moreover, 
$$(P(v))^2=\begin{cases}4-\frac{3v^2}{2} \text{ if }v\in[0,1],\\
\frac{(v - 2)(v - 6)}{2}\text{ if }v\in[1,2].
\end{cases}P(v)\cdot  A=\begin{cases}2v\text{ if }\frac{3v}{2}\in[0,1]\\
2-\frac{v}{2} \text{ if }v\in[1,2].
\end{cases}$$
In this case $\delta_P(S)=\frac{6}{7}\text{ if }P\in A$.
\end{lemma}
\begin{proof}
The Zariski Decomposition in part a). follows from $-K_S-vA\sim_{\DR} (2-v)A+2B_1+C_1+C_1'+B_2+B_3$. A similar statement holds in other parts.
We have
$S_S(A)=\frac{7}{6}$.
Thus, $\delta_P(S)\le \frac{6}{7}$ for $P\in  A$. Note that we have if $P\in A\backslash B_1$ or if $P\in A\cap B_1$:
$$h(v)\le \begin{cases}
\frac{9v^2}{8}\text{ if }v\in[0,1],\\
\frac{(4 - v) (7 v - 4)}{8}\text{ if }v\in[1,2].
\end{cases}
\text{ or }
\begin{cases}
\frac{15v^2}{8}\text{ if }v\in[0,1],\\
\frac{(4 - v) (4+v)}{8}\text{ if }v\in[1,2].
\end{cases}$$
So 
$S(W_{\bullet,\bullet}^{A};P)\le\frac{7}{6}$.
Thus, $\delta_P(S)=\frac{6}{7}$ if $P\in A$.
\end{proof}
\begin{lemma}\label{deg4-67-2_3_points}
 If $P$ belongs to a $(-1)$-curve $A$ and there exist $(-1)$-curves and $(-2)$-curves   which form the following dual graph:
\begin{figure}[h!]
    \centering
 \includegraphics[width=5.5cm]{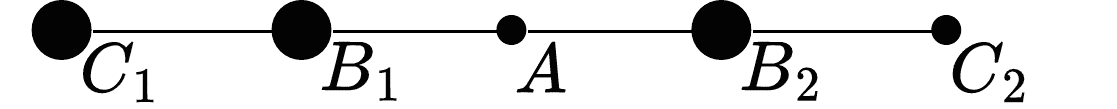}
    \caption{Dual graph: $(-K_S)^2=4$ and $\delta_P(S)=\frac{6}{7}$ with $\tau(A)=3$}
\end{figure}
 \par  
   Then $\tau(A)=3$ and the Zariski Decomposition of the divisor $-K_S-vA$ is given by:
{\allowdisplaybreaks \begin{align*}
&&P(v)=\begin{cases}
-K_S-v A - \frac{v}{3}(2B_1 + C_1)-\frac{v}{2}B_2\text{ if }v\in[0,2],\\
-K_S-v A-\frac{v}{3}(2B_1 + C_1)-(v-1)B_2-(v-2)C_2\text{ if }v\in[2,3].
\end{cases}\\&&N(v)=\begin{cases}
\frac{v}{3}(2B_1 + C_1)+\frac{v}{2}B_2\text{ if }v\in[0,2],\\
\frac{v}{3}(2B_1 + C_1)+(v-1)B_2+(v-2)C_2\text{ if }v\in[2,3].
\end{cases}\end{align*}}
Moreover, 
$$(P(v))^2=\begin{cases} 4-2v+\frac{v^2}{6} \text{ if }v\in[0,2],\\
\frac{2(v - 3)^2}{3}\text{ if }v\in[2,3].
\end{cases}P(v)\cdot  A=\begin{cases}2v\text{ if }1-\frac{v}{6}  \in[0,2],\\
2(1-\frac{v}{3})\text{ if }v\in[2,3].
\end{cases}$$
In this case $\delta_P(S)=\frac{6}{7}\text{ if }P\in A\backslash B_1$.
\end{lemma}
\begin{proof}
The Zariski Decomposition follows from $-K_S-vA\sim_{\DR} (3-v)A+2B_1+C_1+2B_2+C_2$.
We have
$S_S(A)=\frac{7}{6}$. Thus, $\delta_P(S)\le \frac{6}{7}$ for $P\in  A$. Note that for $P\in A\backslash B_1$ we have:
$$h(v)\le
\begin{cases}
\frac{(6 - v) (5 v + 6)}{72}\text{ if }v\in[0,2],\\
\frac{4v (3-v)}{9}\text{ if }v\in[2,3].
\end{cases}$$
So $ S(W_{\bullet,\bullet}^{A};P)\le 1\le\frac{7}{6}$. Thus, $\delta_P(S)=\frac{6}{7}$ if $P\in A\backslash B_1$.
\end{proof}
\begin{lemma}\label{deg4-2429_1_32_2_points}
 If $P$ belongs to a $(-2)$-curve $A$ and there exist $(-1)$-curves and $(-2)$-curves   which form the following dual graph:
\begin{center}
\includegraphics[width=5.5cm]{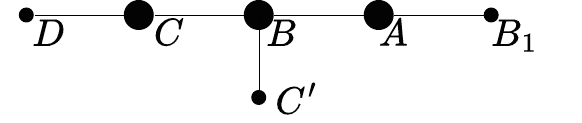}
   \end{center}
   Then $\tau(A)=2$ and the Zariski Decomposition of the divisor $-K_S-vA$ is given by:
{\allowdisplaybreaks \begin{align*}
&&P(v)=\begin{cases}
-K_S-vA-\frac{v}{3}(2B+C)\text{ if }v\in[0,1],\\
-K_S-vA-\frac{v}{3}(2B+C)-(v-1)B_1\text{ if }v\in\big[1,\frac{3}{2}\big],\\
-K_S-vA-(v-1)(2B+C+B_1)-(2v-3)C'\text{ if }v\in\big[\frac{3}{2},2\big].
\end{cases}\\&&N(v)=\begin{cases}\frac{v}{3}(2B+C)\text{ if }v\in[0,1],\\
\frac{v}{3}(2B+C)+(v-1)B_1\text{ if }v\in\big[1,\frac{3}{2}\big],\\
(v-1)(2B+C+B_1)+(2v-3)C'\text{ if }v\in\big[\frac{3}{2},2\big].
\end{cases}
\end{align*}}
Moreover, 
$$(P(v))^2=\begin{cases} 
4-\frac{4v^2}{3} \text{ if }v\in[0,1],\\
5-2v-\frac{v^2}{3} \text{ if }v\in\big[1,\frac{3}{2}\big],\\
(2-v)(4-v) \text{ if }v\in\big[\frac{3}{2},2\big].
\end{cases}P(v)\cdot  A=\begin{cases}
\frac{4v}{3}\text{ if }v\in[0,1],\\
1+\frac{v}{3} \text{ if }v\in\big[1,\frac{3}{2}\big],\\
3-v\text{ if }v\in\big[\frac{3}{2},2\big].
\end{cases}$$
In this case $\delta_P(S)=\frac{24}{29}\text{ if }P\in A\backslash B$.
\end{lemma}
\begin{proof}
The Zariski Decomposition follows from $-K_S-vA\sim_{\DR} (2-v)A+3B+2C+D+2C'+B_1$.
We have 
$S_S(A)=\frac{29}{24}$.
Thus, $\delta_P(S)\le \frac{24}{29}$ for $P\in  A$. Note that for $P\in A\backslash B$ we have:
$$h(v)=\begin{cases}
\frac{8v^2}{9}\text{ if }v\in[0,1],\\
 \frac{(v + 3) (7 v - 3)}{18}\text{ if }v\in\big[1,\frac{3}{2}\big],\\
\frac{(3 - v) (v + 1)}{2}\text{ if }v\in\big[\frac{3}{2},2\big].
\end{cases}$$
So $S(W_{\bullet,\bullet}^{A};P)\le \frac{11}{12}\le\frac{29}{24}$. Thus, $\delta_P(S)=\frac{24}{29}$ if $P\in A\backslash B$.
\end{proof}
\begin{lemma}\label{deg4-2429_52_3_points}
 If $P$ belongs to a $(-1)$-curve $A$ and there exist $(-1)$-curves and $(-2)$-curves   which form the following dual graph:
\begin{figure}[h!]
    \centering
\includegraphics[width=6cm]{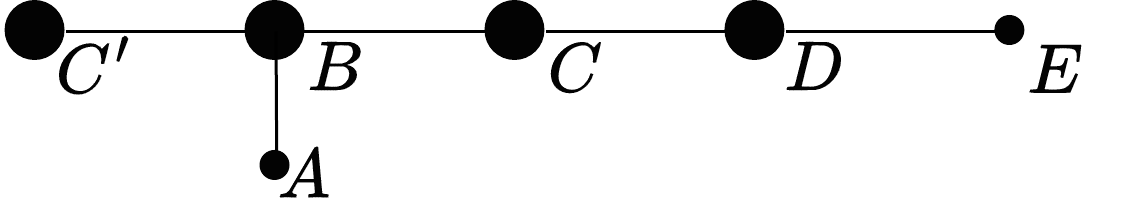}
    \caption{Dual graph: $(-K_S)^2=4$ and $\delta_P(S)=\frac{24}{29}$}
\end{figure}
 \par
   Then $\tau(A)=3$ and the Zariski Decomposition of the divisor $-K_S-vA$ is given by:
{\allowdisplaybreaks \begin{align*}
&&P(v)=\begin{cases}
-K_S-vA- \frac{v}{5}(3C' + 6B + 4C + 2D)\text{ if }v\in\big[0,\frac{5}{2}\big],\\
-K_S-vA-(v-1)(2B+C')-(2v-3)C-(2v-4)D-(2v - 5)E\text{ if }v\in\big[\frac{5}{2},3\big].
\end{cases}\\&&N(v)=
\begin{cases} \frac{v}{5}(3C' + 6B + 4C + 2D)\text{ if }v\in\big[0,\frac{5}{2}\big],\\
(v-1)(2B+C')+(2v-3)C+(2v-4)D+(2v - 5)E\text{ if }v\in\big[\frac{5}{2},3\big].
\end{cases}
\end{align*}}
Moreover, 
$$(P(v))^2=\begin{cases}
4 -2v +\frac{v^2}{5}  \text{ if }v\in\big[0,\frac{5}{2}\big],\\
(3 - v)^2\text{ if }v\in\big[\frac{5}{2},3\big].
\end{cases}P(v)\cdot A=\begin{cases}1-\frac{v}{5}\text{ if } v  \in\big[0,\frac{5}{2}\big],\\
3-v\text{ if }v\in\big[\frac{5}{2},3\big].
\end{cases}$$
In this case $\delta_P(S)=\frac{24}{29}\text{ if }P\in A\backslash B$.
\end{lemma}

\begin{proof}
The Zariski Decomposition follows from $-K_S-vA\sim_{\DR} (3-v)A+2C'+4B+3C+2D+E$. We have
$S_S(A)=\frac{29}{24}$.
Thus, $\delta_P(S)\le \frac{24}{29}$ for $P\in  A$. Note that if $P\in  A\backslash B$ then
$$h(v)=\begin{cases}
\frac{(5 - v)^2}{50}\text{ if }v\in\big[0,\frac{5}{2}\big],\\
\frac{(3-v)^2}{2}\text{ if }v\in\big[\frac{5}{2},3\big].
\end{cases}$$
So $S(W_{\bullet,\bullet}^{A};P)\le \frac{3}{8}\le\frac{29}{24}$.
Thus, $\delta_P(S)=\frac{24}{29}$ if $P\in A\backslash B$.
\end{proof}
\begin{lemma}\label{deg4-911_43_2_points}
 If $P$ belongs to a $(-2)$-curve $A$ and there exist $(-1)$-curves and $(-2)$-curves   which form the following dual graph:
\begin{figure}[h!]
    \centering
\includegraphics[width=5cm]{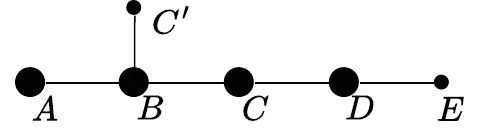}
    \caption{Dual graph: $(-K_S)^2=4$ and $\delta_P(S)=\frac{9}{11}$}
\end{figure}
 \par    Then $\tau(A)=2$ and the Zariski Decomposition of the divisor $-K_S-vA$ is given by:
{\allowdisplaybreaks \begin{align*}
&&P(v)=\begin{cases}
-K_S-vA-\frac{v}{4}(3B+2C+D)\text{ if }v\in\big[0,\frac{4}{3}\big],\\
-K_S-vA-(v-1)(3B+2C+D)-(3v-4)C'\text{ if }v\in\big[\frac{4}{3},2\big].
\end{cases}\\&&N(v)=\begin{cases}\frac{v}{4}( 3B+2C+D)\text{ if }v\in\big[0,\frac{4}{3}\big],\\
(v-1)(3B+2C+D)+(3v-4)C'\text{ if }v\in\big[\frac{4}{3},2\big].
\end{cases}
\end{align*}}
Moreover, 
$$(P(v))^2=\begin{cases} 4-\frac{5v^2}{4} \text{ if }v\in\big[0,\frac{4}{3}\big],\\
(2 - v)(4-v)\text{ if }v\in\big[\frac{4}{3},2\big].
\end{cases}P(v)\cdot  A=\begin{cases}\frac{5v}{4}\text{ if } v  \in\big[0,\frac{4}{3}\big],\\
3-v\text{ if }v\in\big[\frac{4}{3},2\big].
\end{cases}$$
In this case $\delta_P(S)=\frac{9}{11}\text{ if }P\in A\backslash B$.
\end{lemma}
\begin{proof}
The Zariski Decomposition follows from $-K_S-vA\sim_{\DR} (2-v)A+2C'+3B+3C+2D+E$. We have
$S_S(A)=\frac{11}{9}$. Thus, $\delta_P(S)\le \frac{9}{11}$ for $P\in  A$. Note that we have that if $P\in  A\backslash B$ then
$$h(v)\le\begin{cases}
\frac{25v^2}{32}\text{ if }v\in\big[0,\frac{4}{3}\big],\\
\frac{(3-v)^2}{2}\text{ if }v\in\big[\frac{4}{3},2\big].
\end{cases}$$
So  $S(W_{\bullet,\bullet}^{A};P)\le \frac{11}{18}\le\frac{11}{9}$. Thus, $\delta_P(S)=\frac{9}{11}$ if $P\in A\backslash B$.
\end{proof}
\begin{lemma}\label{deg4-45_1_2_points}
 If $P$ belongs to a $(-2)$-curve $A$ and there exist $(-1)$-curves and $(-2)$-curves   which form the following dual graph:
\begin{figure}[h!]
    \centering
\includegraphics[width=5.3cm]{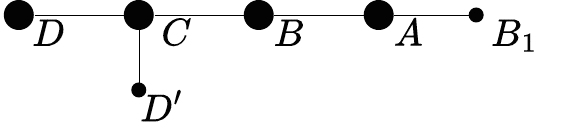}
    \caption{Dual graph: $(-K_S)^2=4$ and $\delta_P(S)=\frac{4}{5}$}
\end{figure}
 \par  Then $\tau(A)=2$ and the Zariski Decomposition of the divisor $-K_S-vA$ is given by:
{\allowdisplaybreaks \begin{align*}
&&P(v)=\begin{cases}
-K_S-vA-\frac{v}{4}( 3B+2C+D)\text{ if }v\in[0,1],\\
-K_S-vA-\frac{v}{4}(3B+2C+D)-(v-1)B_1\text{ if }v\in[1,2].
\end{cases}\\&&N(v)=\begin{cases}\frac{v}{4}( 3B+2C+D)\text{ if }v\in[0,1],\\
\frac{v}{4}(3B+2C+D)+(v-1)B_1\text{ if }v\in[1,2].
\end{cases}
\end{align*}}
Moreover, 
$$(P(v))^2=\begin{cases} 
4-\frac{5v^2}{4} \text{ if }v\in[0,1],\\
\frac{(v + 10)(2 - v)}{4}\text{ if }v\in[1,2].
\end{cases}P(v)\cdot  A=\begin{cases}\frac{5v}{4}\text{ if } v  \in[0,1],\\
1+\frac{v}{4} \text{ if }v\in[1,2].
\end{cases}$$
In this case $\delta_P(S)=\frac{4}{5}\text{ if }P\in A\backslash B$.
\end{lemma}

\begin{proof}
The Zariski Decomposition follows from $-K_S-vA\sim_{\DR} (2-v)A+3B+4C+2D+3D'+B_1$. We have
$S_S(A)=\frac{5}{4}$.
Thus, $\delta_P(S)\le \frac{4}{5}$ for $P\in  A$. Note that for $P\in A\backslash B$ we have:
$$h(v)\le \begin{cases}\frac{25v^2}{32}\text{ if }v\in[0,1],\\
\frac{(v + 4) (9 v - 4)}{32}\text{ if }v\in[1,2].
\end{cases}$$
So $ S(W_{\bullet,\bullet}^{A};P)\le \frac{23}{24}\le\frac{5}{4}$.
Thus, $\delta_P(S)=\frac{4}{5}$ if $P\in A\backslash B$.
\end{proof}
\begin{lemma}\label{deg4-centerA3points}
 If $P$ belongs to a $(-2)$-curve $A$ and there exist $(-1)$-curves and $(-2)$-curves   which form the following dual graph:
\begin{figure}[h!]
    \centering
 \includegraphics[width=10cm]{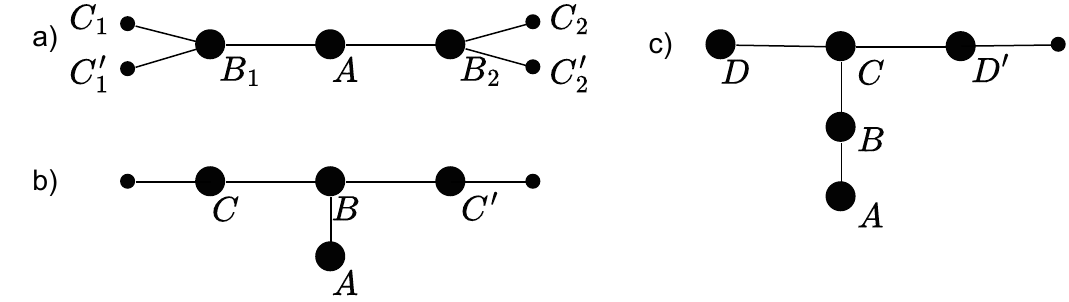}
    \caption{Dual graph: $(-K_S)^2=4$ and $\delta_P(S)=\frac{3}{4}$ with $\tau(A)=2$ }
\end{figure}
\par    Then $\tau(A)=2$ and the Zariski Decomposition of the divisor $-K_S-vA$ is given by:
{\allowdisplaybreaks \begin{align*}
&{\text{\bf a). }}&P(v)=-K_S-vA-\frac{v}{2}(B_1+B_2) \text{ and }N(v)=\frac{v}{2}(B_1+B_2)\text{ if }v\in[0,2].\\
&{\text{\bf b). }}&P(v)=-K_S-vA-\frac{v}{2}( 2B+C+C') \text{ and }N(v)=\frac{v}{2}(2B+C+C')\text{ if }v\in[0,2].\\
&{\text{\bf c). }}&P(v)=-K_S-vA-\frac{v}{2}(2B+2C+D+D') \text{ and }N(v)=\frac{v}{2}( 2B+2C+D+D')\text{ if }v\in[0,2].
\end{align*}}
Moreover, 
$$(P(v))^2=(2 - v)(v + 2) \text{ and }
P(v)\cdot A=v\text{ if }v\in[0,2].$$
In this case $\delta_P(S)=\frac{3}{4}$ if $P\in A\backslash B$.
\end{lemma}
\begin{proof}
The Zariski Decomposition in part a). follows from $-K_S-vA\sim_{\DR} (2-v)A+2B_1+C_1+C_1'+2B_2+C_2+C_2'$. A similar statement holds in other parts.
We have
$S_{S}(A)=\frac{4}{3}$.
Thus, $\delta_P(S)\le \frac{3}{4}$ for $P\in A$. Moreover $h(v)\le v^2\text{ if }v\in[0,2].$
So
$S(W_{\bullet,\bullet}^{A};P)\le \frac{4}{3}$.
We get that $\delta_P(S)= \frac{3}{4}$  $P\in A\backslash B$.
\end{proof}
\begin{lemma}\label{deg4-pointA3bound}
 If $P$ belongs to a $(-2)$-curve $A$ and there exist $(-1)$-curves and $(-2)$-curves   which form the following dual graph:
\begin{figure}[h!]
    \centering
 \includegraphics[width=12cm]{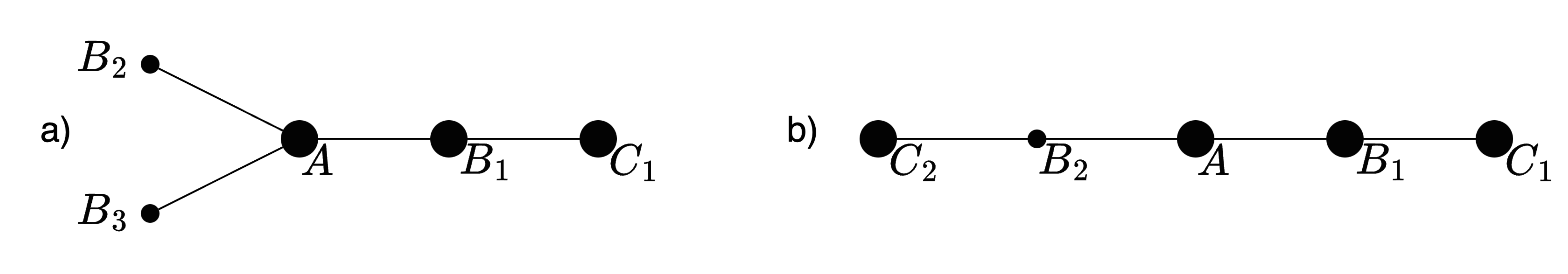}
    \caption{Dual graph: $(-K_S)^2=4$ and $\delta_P(S)=\frac{3}{4}$ with $\tau(A)=3$ }
\end{figure}
\par    Then $\tau(A)=3$ and the Zariski Decomposition of the divisor $-K_S-vA$ is given by:
{{\allowdisplaybreaks \begin{align*}
&{\text{\bf a). }}&P(v)=\begin{cases}
-K_S-vA-\frac{v}{3}(2B_1+C_1)\text{ if }v\in[0,1],\\
-K_S-vA-\frac{v}{3}(2B_1+C_1)-(v-1)(B_2+B_3)\text{ if }v\in[1,3].
\end{cases}\\&&N(v)=\begin{cases}\frac{v}{3}(2B_1+C_1)\text{ if }v\in[0,1],\\
\frac{v}{3}(2B_1+C_1)+(v-1)(B_2+B_3)\text{ if }v\in[1,3].
\end{cases}\\
&{\text{\bf b). }}&P(v)=\begin{cases}
-K_S-vA-\frac{v}{3}(2B_1+C_1)\text{ if }v\in[0,1],\\
-K_S-vA-\frac{v}{3}(2B_1+C_1)-(v-1)(2B_2+C_2)\text{ if }v\in[1,3].
\end{cases}\\&&N(v)=\begin{cases}\frac{v}{3}(2B_1+C_1)\text{ if }v\in[0,1],\\
\frac{v}{3}(2B_1+C_1)+(v-1)(2B_2+C_2)\text{ if }v\in[1,3].
\end{cases}
\end{align*}}}
Moreover, 
$$(P(v))^2=\begin{cases} 4-\frac{4v^2}{3}\text{ if }v\in[0,1],\\
\frac{2(3-v)^2}{3}\text{ if }v\in[1,3].
\end{cases}P(v)\cdot  A=\begin{cases}\frac{4v}{3} \text{ if } v \in[0,1],\\
2(1-\frac{v}{3})   \text{ if }v\in[1,3].
\end{cases}$$
In this case 
$\delta_P(S)=\frac{3}{4}$ if $P\in A\backslash B_1$.
\end{lemma}
\begin{proof}
The Zariski Decomposition in part a). follows from $-K_S-vA\sim_{\DR} (3-v)A+2B_1+C_1+2B_2+2B_3$. A similar statement holds in other parts.
$S_S(A)=\frac{4}{3}$.
Thus, $\delta_P(S)\le \frac{3}{4}$ for $P\in  A$. Note that for $P\in A\backslash B_1 $ we have:
$$h(v)\le\begin{cases}
 \frac{8v^2}{9}\text{ if }v\in[0,1],\\
\frac{ 2 (3 - v) (5 v - 3)}{9}\text{ if }v\in[1,3].
\end{cases}$$
So $S(W_{\bullet,\bullet}^{A};P)\le \frac{4}{3}$.
Thus, $\delta_P(S)=\frac{3}{4}$ if $P\in A \backslash B_1$.
\end{proof}
\begin{lemma}\label{deg4-nearA1A3points}\label{deg4-34_4_points}
 If $P$ belongs to a $(-1)$-curve $A$ and there exist $(-1)$-curves and $(-2)$-curves   which form the following dual graph:
\begin{figure}[h!]
    \centering
 \includegraphics[width=12cm]{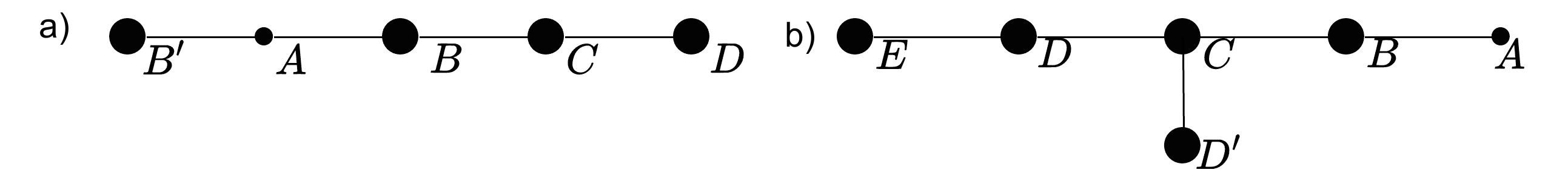}
    \caption{Dual graph: $(-K_S)^2=4$ and $\delta_P(S)=\frac{3}{4}$ with $\tau(A)=4$}
\end{figure}
 \par   Then $\tau(A)=4$ and the Zariski Decomposition of the divisor $-K_S-vA$ is given by:
{\allowdisplaybreaks \begin{align*}
&{\text{\bf a). }}&P(v)=-K_S-vA-\frac{v}{4}(2B'+ 3B+2C+D)\text{ if }v\in[0,4]. \\&& N(v)=\frac{v}{4}(2B'+ 3B+2C+D)\text{ if }v\in[0,4].\\
&{\text{\bf b). }}&P(v)=-K_S-vA-\frac{v}{4}( 5B+6C+4D+2E+3D') \text{ if }v\in[0,4].\\&&N(v)=\frac{v}{4}(5B+6C+4D+2E+3D')\text{ if }v\in[0,4].
\end{align*}}
Moreover, 
$$(P(v))^2=\frac{(v - 4)^2}{4 }
 P(v)\cdot A=1-\frac{v}{4}\text{ if }v\in[0,4].$$
In this case $\delta_P(S)=\frac{3}{4}\text{ if }P\in A$.
\end{lemma}
\begin{proof}
The Zariski Decomposition in part a). follows from $-K_S-vA\sim_{\DR} (4-v)A+2B'+ 3B+2C+D$. A similar statement holds in other parts.
We have
$S_{S}(A)=\frac{4}{3}$.
Thus, $\delta_P(S)\le \frac{3}{4}$ for $P\in A$. Moreover for $P\in A\backslash B_1$ we have
$h(v)\le  \frac{ (4 - v) (3 v + 4)}{32}\text{ if }v\in[0,4].$
So 
$S(W_{\bullet,\bullet}^{A};P)\le 1\le \frac{4}{3} $.
We get that $\delta_P(S)= \frac{3}{4}$ for $P\in A\backslash B_1$.
\end{proof}
\begin{lemma}\label{deg4-23_1_2_3_points}
 If $P$ belongs to a $(-2)$-curve $A$ and there exist $(-1)$-curves and $(-2)$-curves   which form the following dual graph:
\begin{figure}[h!]
    \centering
\includegraphics[width=13.5cm]{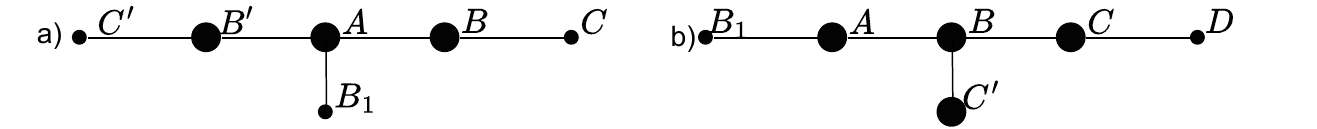}
    \caption{Dual graph: $(-K_S)^2=4$ and $\delta_P(S)=\frac{2}{3}$}
\end{figure}
 \par    Then $\tau(A)=3$ and the Zariski Decomposition of the divisor $-K_S-vA$ is given by:
{\allowdisplaybreaks \begin{align*}
&{\text{\bf a). }}&P(v)=\begin{cases}
-K_S-vA-\frac{v}{2}(B+B')\text{ if }v\in[0,1],\\
-K_S-vA-\frac{v}{2}(B+B')-(v-1)B_1\text{ if }v\in[1,2],\\
-K_S-vA-(v-1)(B+B'+B_1)-(v-2)(C+C')\text{ if }v\in[2,3].
\end{cases}\\&&N(v)=\begin{cases}
\frac{v}{2}(B+B')\text{ if }v\in[0,1],\\
\frac{v}{2}(B+B')+(v-1)B_1\text{ if }v\in[1,2],\\
(v-1)(B+B'+B_1)+(v-2)(C+C')\text{ if }v\in[2,3].
\end{cases}\\
&{\text{\bf b). }}&P(v)=\begin{cases}
-K_S-vA-\frac{v}{2}(2B+C+C')\text{ if }v\in[0,1],\\
-K_S-vA-\frac{v}{2}(2B+C+C')-(v-1)B_1\text{ if }v\in[1,2],\\
-K_S-vA-(v-1)(2B+C'+B_1)-(2v-3)C-(2v-4)D\text{ if }v\in[2,3].
\end{cases}\\&&N(v)=\begin{cases}
\frac{v}{2}(2B+C+C')\text{ if }v\in[0,1],\\
\frac{v}{2}(2B+C+C')+(v-1)B_1\text{ if }v\in[1,2],\\
(v-1)(2B+C'+B_1)+(2v-3)C+(2v-4)D\text{ if }v\in[2,3].
\end{cases}
\end{align*}}
Moreover, 
$$(P(v))^2=\begin{cases} (2-v)(2+v) \text{ if }v\in[0,1],\\
5 - 2v \text{ if }v\in[1,2],\\
(3 - v)^2\text{ if }v\in[2,3].
\end{cases}P(v)\cdot  A=\begin{cases} v\text{ if } v \in[0,1],\\
1\text{ if } v \in[1,2],\\
3-v \text{ if }v\in[2,3].
\end{cases}$$
In this case $\delta_P(S)=\frac{2}{3}\text{ if }P\in A\backslash B$.
\end{lemma}
\begin{proof}
The Zariski Decomposition in part a). follows from $-K_S-vA\sim_{\DR} (3-v)A+C'+2B'+2B+C'+2B_1$.
We have
$S_S(A)=\frac{3}{2}$.
Thus, $\delta_P(S)\le \frac{2}{3}$ for $P\in  A$. Note that fo $P\in A\backslash B$ we have:
$$h(v)\le\begin{cases} 
\frac{v^2}{2}\text{ if } v \in[0,1],\\
v+\frac{1}{2}\text{ if } v \in[1,2],\\
\frac{ (3 - v) (v + 1)}{2} \text{ if }v\in[2,3].
\end{cases}$$
So $S(W_{\bullet,\bullet}^{A};P)\le 1\le\frac{3}{2}$.
Thus, $\delta_P(S)=\frac{2}{3}$ if $P\in A\backslash B$.
\end{proof}
\begin{lemma}\label{deg4-2437_32_2_3_points}
 If $P$ belongs to a $(-2)$-curve $A$ and there exist $(-1)$-curves and $(-2)$-curves   which form the following dual graph:
\begin{figure}[h!]
    \centering
\includegraphics[width=5.5cm]{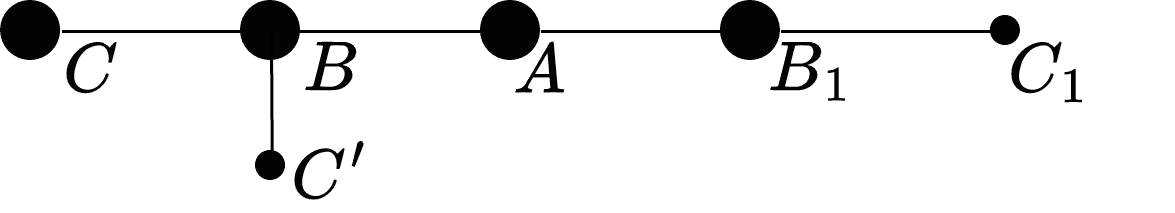}
    \caption{Dual graph: $(-K_S)^2=4$ and $\delta_P(S)=\frac{24}{37}$}
\end{figure}
 \par   Then $\tau(A)=3$ and the Zariski Decomposition of the divisor $-K_S-vA$ is given by:
{\allowdisplaybreaks \begin{align*}
&&P(v)=\begin{cases}
-K_S-vA-\frac{v}{2}B_1-\frac{v}{3}(2B+C)\text{ if }v\in\big[0,\frac{3}{2}\big],\\
-K_S-vA-\frac{v}{2}B_1-(v-1)(2B+C)-(2v-3)C'\text{ if }v\in\big[\frac{3}{2},2\big],\\
-K_S-vA-(v-1)(2B+C+B_1)-(2v-3)C'-(v-2)C_1\text{ if }v\in[2,3].
\end{cases}\\&&N(v)=\begin{cases}
\frac{v}{2}B_1+\frac{v}{3}(2B+C)\text{ if }v\in\big[0,\frac{3}{2}\big],\\
\frac{v}{2}B_1+(v-1)(2B+C)+(2v-3)C'\text{ if }v\in\big[\frac{3}{2},2\big],\\
(v-1)(2B+C+B_1)+(2v-3)C'+(v-2)C_1\text{ if }v\in[2,3].
\end{cases}
\end{align*}}
Moreover, 
$$(P(v))^2=\begin{cases} 4-\frac{5v^2}{6}  \text{ if }v\in\big[0,\frac{3}{2}\big],\\
7-4v+\frac{v^2}{2}  \text{ if }v\in\big[\frac{3}{2},2\big],\\
(3 - v)^2\text{ if }v\in[2,3].
\end{cases}P(v)\cdot  A\begin{cases} \frac{5v}{6} \text{ if } v \in\big[0,\frac{3}{2}\big],\\
2-\frac{v}{2} \text{ if } v \in\big[\frac{3}{2},2\big],\\
3-v \text{ if }v\in[2,3].
\end{cases}$$
In this case $\delta_P(S)=\frac{24}{37}\text{ if }P\in A\backslash B$.
\end{lemma}

\begin{proof}
The Zariski Decomposition follows from $-K_S-vA\sim_{\DR} (3-v)A+2B_1+C_1+4B+2C+3C'$.
We have
$S_S(A)=\frac{37}{24}$.
Thus, $\delta_P(S)\le \frac{24}{37}$ for $P\in  A$. Note that for $P\in A\backslash B$ we have: 
$$h(v)\le \begin{cases} 
\frac{55v^2}{72}\text{ if } v \in\big[0,\frac{3}{2}\big],\\
\frac{(4-v)(4+v)}{8}\text{ if } v \in\big[\frac{3}{2},2\big],\\
\frac{ (3 - v) (v + 1)}{2} \text{ if }v\in[2,3].
\end{cases}$$
So $S(W_{\bullet,\bullet}^{A};P)\le \frac{13}{12}\le\frac{37}{24}$.
Thus, $\delta_P(S)=\frac{24}{37}$ if $P\in A\backslash B$.
\end{proof}
\begin{lemma}\label{deg4-914_53_3_points}
 If $P$ belongs to a $(-2)$-curve $A$ and there exist $(-1)$-curves and $(-2)$-curves   which form the following dual graph:
 \begin{figure}[h!]
    \centering
\includegraphics[width=6cm]{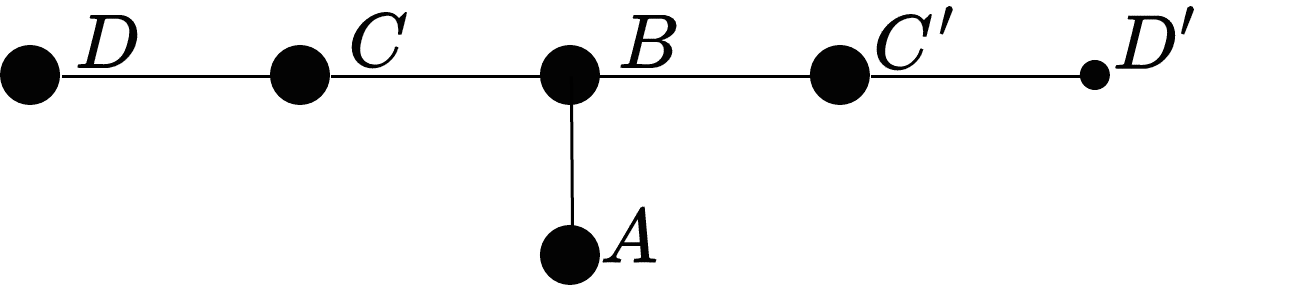}
    \caption{Dual graph: $(-K_S)^2=4$ and $\delta_P(S)=\frac{9}{14}$}
\end{figure}
 \par     Then $\tau(A)=3$ and the Zariski Decomposition of the divisor $-K_S-vA$ is given by:
{\allowdisplaybreaks \begin{align*}
&&P(v)=\begin{cases}
-K_S-v A-\frac{v}{5}(6B+4C+3C'+2D)\text{ if }v\in\big[0,\frac{5}{3}\big],\\
-K_S-vA-(v-1)(3B+2C+D)-(3v-4)C'-(3v-5)D'\text{ if }v\in\big[\frac{5}{3},3\big].
\end{cases}\\&&N(v)=\begin{cases}\frac{v}{5}(6B+4C+3C'+2D)\text{ if }v\in\big[0,\frac{5}{3}\big],\\
(v-1)(3B+2C+D)+(3v-4)C'+(3v-5)D'\text{ if }v\in\big[\frac{5}{3},3\big].
\end{cases}
\end{align*}}
Moreover, 
$$(P(v))^2=\begin{cases}4-\frac{4v^2}{5}   \text{ if }v\in\big[0,\frac{5}{3}\big],\\
(3-v)^2\text{ if }v\in\big[\frac{5}{3},3\big].
\end{cases}P(v)\cdot  A=\begin{cases}\frac{4v}{5}\text{ if }v\in\big[0,\frac{5}{3}\big],\\
3-v\text{ if }v\in\big[\frac{5}{3},3\big].
\end{cases}$$
In this case $\delta_P(S)=\frac{9}{14}\text{ if }P\in A\backslash B$.
\end{lemma}
\begin{proof}
The Zariski Decomposition follows from $-K_S-vA\sim_{\DR} (3-v)A+4D'+5C'+6B+4C+2D$.
We have $S_S(A)=\frac{14}{9}$.
Thus, $\delta_P(S)\le \frac{9}{14}$ for $P\in A$. Note that if $P\in A\backslash B$ then 
    $$h(v)=\begin{cases}\frac{8v^2}{25} \text{ if }v\in\big[0,\frac{5}{3}\big],\\
\frac{(3-v)^2}{2}\text{ if }v\in\big[\frac{5}{3},3\big].
\end{cases}$$
So $S(W_{\bullet,\bullet}^{A};P)\le\frac{14}{9}$.
Thus, $\delta_P(S)=\frac{9}{14}$ if $P\in A\backslash B$.
\end{proof}
\begin{lemma}\label{deg4-611_1_3_4_points}
 If $P$ belongs to a $(-2)$-curve $A$ and there exist $(-1)$-curves and $(-2)$-curves   which form the following dual graph:
\begin{figure}[h!]
    \centering
\includegraphics[width=6cm]{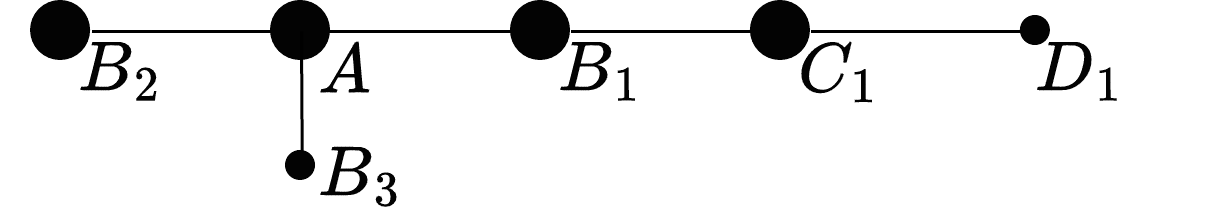}
    \caption{Dual graph: $(-K_S)^2=4$ and $\delta_P(S)=\frac{6}{11}$}
\end{figure}
 \par   Then $\tau(A)=4$ and the Zariski Decomposition of the divisor $-K_S-vA$ is given by:
{\allowdisplaybreaks \begin{align*}
&&P(v)=\begin{cases}
-K_S-vA-\frac{v}{3}(2B_1+C_1)-\frac{v}{2}B_2\text{ if }v\in[0,1],\\
-K_S-vA-\frac{v}{3}(2B_1+C_1)-\frac{v}{2}B_2-(v-1)B_3\text{ if }v\in[1,3],\\
-K_S-vA-(v-1)(B_1+B_3)-(v-2)C_1-(v-3)D_1-\frac{v}{2}B_2\text{ if }v\in[3,4].
\end{cases}\\&&N(v)=
\begin{cases}
\frac{v}{3}(2B_1+C_1)+\frac{v}{2}B_2\text{ if }v\in[0,1],\\
\frac{v}{3}(2B_1+C_1)+\frac{v}{2}B_2+(v-1)B_3\text{ if }v\in[1,3],\\
(v-1)(B_1+B_3)+(v-2)C_1+(v-3)D_1+\frac{v}{2}B_2\text{ if }v\in[3,4].
\end{cases}
\end{align*}}
Moreover, 
$$(P(v))^2=\begin{cases} 4-\frac{5v^2}{6}  \text{ if }v\in[0,1],\\
5-2v+\frac{v^2}{6}  \text{ if }v\in[1,3],\\
\frac{(4-v)^2}{2}\text{ if }v\in[3,4].
\end{cases}P(v)\cdot  A=\begin{cases} \frac{5v}{6} \text{ if } v \in[0,1],\\
1-\frac{v}{6} \text{ if } v \in[1,3],\\
2-\frac{v}{2}  \text{ if }v\in[3,4].
\end{cases}$$
In this case $\delta_P(S)=\frac{6}{11}\text{ if }P\in A$.
\end{lemma}

\begin{proof}
The Zariski Decomposition follows from $-K_S-vA\sim_{\DR} (4-v)A+3B_1+2C_1+D_1+2B_2+3B_3$.
We have $S_S(A)=\frac{11}{6}$.
Thus, $\delta_P(S)\le \frac{6}{11}$ for $P\in  A$. Note that if $P\in A\backslash (B_1\cup B_3)$ or $P\in A\cap B_1$ or $P\in A\cap B_3$ we have:
$$h(v)\le\begin{cases} 
\frac{55v^2}{72}\text{ if } v \in[0,1],\\
\frac{(6 - v) (5 v + 6)}{72}\text{ if } v \in[1,3],\\
\frac{ (4-v) (v + 4)}{8} \text{ if }v\in[3,4].
\end{cases}
\text{ or }
h(v)\le \begin{cases} \frac{ 65v^2}{72}\text{ if } v \in[0,1],\\
\frac{ (6 - v) (7 v + 6)}{72}\text{ if } v \in[1,3],\\
\frac{ 3 (4 - v) v}{8} \text{ if }v\in[3,4].
\end{cases}$$$$
\text{ or }
h(v)\le \begin{cases} \frac{25v^2}{72}\text{ if } v \in[0,1],\\
\frac{(6 - v) (11 v - 6)}{72}\text{ if } v \in[1,3],\\
\frac{ 3 (4 - v) v}{8} \text{ if }v\in[3,4].
\end{cases}
$$
So $S(W_{\bullet,\bullet}^{A};P)\le \frac{11}{9}\le\frac{11}{6}$ or $S(W_{\bullet,\bullet}^{A};P)\le \frac{37}{24}\le\frac{11}{6}$ $S(W_{\bullet,\bullet}^{A};P)\le \frac{29}{24}\le\frac{11}{6}$.
Thus, $\delta_P(S)=\frac{6}{11}$ if $P\in A$.
\end{proof}
\begin{lemma}\label{deg4-nearDpoints}
 If $P$ belongs to a $(-1)$-curve $A$ and there exist $(-1)$-curves and $(-2)$-curves   which form the following dual graph:
\begin{figure}[h!]
    \centering
 \includegraphics[width=11cm]{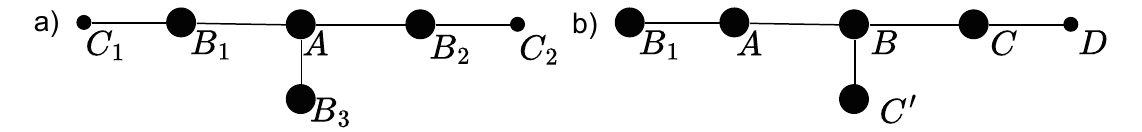}
    \caption{Dual graph: $(-K_S)^2=4$ and $\delta_P(S)=\frac{1}{2}$ with $\tau(A)=4$}
\end{figure}
 \par   Then $\tau(A)=4$ and the Zariski Decomposition of the divisor $-K_S-vA$ is given by:
{\allowdisplaybreaks \begin{align*}
&{\text{\bf a). }}&P(v)=\begin{cases}-K_S-vA-\frac{v}{2}(B_1+B_2+B_3)\text{ if }v\in[0,2],\\
-K_S-vA-\frac{v}{2}B_1-(v-1)(B_2+B_3)-(v-2)(C_2+C_3)\text{ if }v\in[2,4].\\
\end{cases}\\
&&N(v)=\begin{cases}\frac{v}{2}(B_1+B_2+B_3)\text{ if }v\in[0,2],\\
\frac{v}{2}B_1+(v-1)(B_2+B_3)+(v-2)(C_2+C_3)\text{ if }v\in[2,4].
\end{cases}\\
&{\text{\bf b). }}&P(v)=\begin{cases}-K_S-vA-\frac{v}{2} (B_1+2B+C+C')\text{ if }v\in[0,2],\\
-K_S-vA-\frac{v}{2}E_1-(v-1)(2B+C')-(2v-3)C-(2v-4)D\text{ if }v\in[2,4].
\end{cases}\\
&&N(v)=\begin{cases}\frac{v}{2}(B_1+2B+C+C')\text{ if }v\in[0,2],\\
\frac{v}{2}B_1+(v-1)(2B+C')+(2v-3)C+(2v-4)D\text{ if }v\in[2,4].
\end{cases}
\end{align*}}
Moreover, 
$$(P(v))^2=\begin{cases} 4-\frac{v^2}{2}  \text{ if }v\in[0,2],\\
\frac{(4-v)^2}{2}\text{ if }v\in[2,4].
\end{cases}P(v)\cdot  A=\begin{cases}\frac{v}{2}\text{ if } v  \in[0,2]\\
2-\frac{v}{2} \text{ if }v\in[2,4].
\end{cases}$$
In this case $\delta_P(S)=\frac{1}{2}$ if $P\in A\backslash B$.
\end{lemma}
\begin{proof}
The Zariski Decomposition in part a). follows from $-K_S-vA\sim_{\DR} (4-v)A+3B_1+2C_1+3B_2+2C_2+2B_3$. A similar statement holds in other parts.
Thus,
$S_S(A)=2$
Thus, $\delta_P(S)\le \frac{1}{2}$ for $P\in  A$. Note that we have that 
$$h(v)\le \begin{cases}
\frac{3v^2}{8}\text{ if }v\in[0,2],\\
\frac{ (4 - v) (v + 4)}{8} \text{ if }v\in[2,4].
\end{cases}\text{ or }
h(v)=\begin{cases}\frac{3v^2}{8}\text{ if }v\in[0,2],\\
\frac{ 3 (4 - v) v}{8} \text{ if }v\in[2,4].
\end{cases}$$
So $ S(W_{\bullet,\bullet}^{A};P)\le \frac{3}{2}\le 2$
or $ S(W_{\bullet,\bullet}^{A};P)\le \frac{4}{3}\le 2$. Thus, $\delta_P(S)=\frac{1}{2}$.
\end{proof}
\begin{lemma}\label{deg4-12_1_5_points}
 If $P$ belongs to a $(-2)$-curve $A$ and there exist $(-1)$-curves and $(-2)$-curves   which form the following dual graph:
 \begin{figure}[h!]
    \centering
\includegraphics[width=6.5cm]{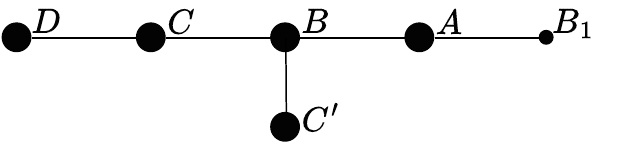}
    \caption{Dual graph: $(-K_S)^2=4$ and $\delta_P(S)=\frac{1}{2}$ with $\tau(A)=5$}
\end{figure}
 \par   Then $\tau(A)=5$ and the Zariski Decomposition of the divisor $-K_S-vA$ is given by:
 {\allowdisplaybreaks \begin{align*}
&&P(v)=\begin{cases}
-K_S-v A-\frac{v}{5}(6B+4C+3C'+2D)\text{ if }v\in[0,1],\\
-K_S-vA-\frac{v}{5}(6B+4C+3C'+2D)-(v-1)B_1\text{ if }v\in[1,5].
\end{cases}\\&&N(v)=\begin{cases}\frac{v}{5}(6B+4C+3C'+2D)\text{ if }v\in[0,1],\\
\frac{v}{5}(6B+4C+3C'+2D)+(v-1)B_1\text{ if }v\in[1,5].
\end{cases}
\end{align*}}
Moreover, 
$$(P(v))^2=\begin{cases}4-\frac{4v^2}{5}  \text{ if }v\in[0,1],\\
\frac{(5-v)^2}{5}\text{ if }v\in[1,5].
\end{cases}P(v)\cdot  A=\begin{cases}
\frac{4v}{5}\text{ if }v\in[0,1],\\
1-\frac{v}{5}\text{ if }v\in[1,5].
\end{cases}$$
In this case $\delta_P(S)=\frac{1}{2}\text{ if }P\in A\backslash B$.
\end{lemma}
\begin{proof}
The Zariski Decomposition follows from $-K_S-vA\sim_{\DR} (5-v)A+3C'+6B+4C+2D+4B_1$.
We have
$S_S(A)=2$.
Thus, $\delta_P(S)\le \frac{1}{2}$ for $P\in A$. Note that for $P\in A\backslash B$ we have:
$$h(v)=\begin{cases}
\frac{8v^2}{25} \text{ if }v\in[0,1],\\
\frac{ (5 - v) (9 v - 5)}{50} \text{ if }v\in[1,5].
\end{cases}$$
So $S(W_{\bullet,\bullet}^{A};P)\le\frac{4}{3}\le 2$.
Thus, $\delta_P(S)=\frac{1}{2}$ if $P\in  A\backslash B$.
\end{proof}
\begin{lemma}\label{deg4-38_2_6_points}
 If $P$ belongs to a $(-2)$-curve $A$ and there exist $(-1)$-curves and $(-2)$-curves   which form the following dual graph:
\begin{figure}[h!]
    \centering
\includegraphics[width=6.5cm]{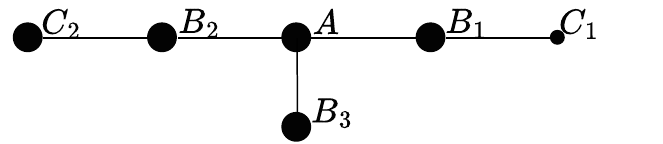}
    \caption{Dual graph: $(-K_S)^2=4$ and $\delta_P(S)=\frac{3}{8}$}
\end{figure}
 \par     Then $\tau(A)=6$ and the Zariski Decomposition of the divisor $-K_S-vA$ is given by:
 {\allowdisplaybreaks \begin{align*}
&&P(v)=\begin{cases}
-K_S-v A-\frac{v}{6}(4C_2+2B_2+3B_3+3B_1)\text{ if }v\in[0,2],\\
-K_S-vA-\frac{v}{6}(4C_2+2B_2+3B_3)-(v-1)B_1-(v-2)C_1\text{ if }v\in[2,6].
\end{cases}\\&&N(v)=\begin{cases}\frac{v}{6}(4C_2+2B_2+3B_3+3B_1)\text{ if }v\in[0,2],\\
\frac{v}{6}(4C_2+2B_2+3B_3)+(v-1)B_1+(v-2)C_1\text{ if }v\in[2,6].
\end{cases}
\end{align*}}
Moreover, 
$$(P(v))^2=\begin{cases} 
4-\frac{v^2}{3}  \text{ if }v\in[0,2],\\
\frac{(6-v)^2}{6}\text{ if }v\in[2,6].
\end{cases}P(v)\cdot  A=\begin{cases}\frac{v}{3}\text{ if }v\in[0,2],\\
1-\frac{v}{6} \text{ if }v\in[2,6].
\end{cases}$$
In this case $\delta_P(S)=\frac{3}{8}\text{ if }P\in A$.
\end{lemma}
\begin{proof}
The Zariski Decomposition follows from $-K_S-vA\sim_{\DR} (6-v)A+5B_1+4C_1+4B_2+2C_2+3B_3$.
We have
$S_S(A)=\frac{8}{3}$.
Thus, $\delta_P(S)\le \frac{3}{8}$ for $P\in A$. Note that
$$h(v)\le\begin{cases}
\frac{5v^2}{18} \text{ if }v\in[0,2],\\
\frac{ (6 - v) (7 v + 6)}{72}\text{ if }v\in[2,6].
\end{cases}
\text{ or }
h(v)\le \begin{cases}
\frac{2v^2}{9} \text{ if }v\in[0,2],\\
\frac{(6 - v) (5 v + 6)}{72}\text{ if }v\in[2,6].
\end{cases}$$
So 
$S(W_{\bullet,\bullet}^{A};P)\le \frac{8}{3}$
or
$S(W_{\bullet,\bullet}^{A};P)\le\frac{14}{9}\le \frac{8}{3}$.
Thus, $\delta_P(S)=\frac{3}{8}$ if $P\in  A$.
\end{proof}

\subsection{Finding $\delta$-invariants for degree $4$}

Let $X$ be a singular del Pezzo surface of degree $4$ with and $S$ be a minimal resolution of $X$. Then there are several possible cases:
\begin{itemize}
    \item[I.] $X$ has an $\DA_1$ singularity and contains $12$ lines. In this case, we let $E$ be the exceptional divisor, $L_{i}$, $L_{i,j}$ for $i\in\{1,2,3,4\}$, $j\in \{1,2\}$ be the lines on $S$, 
    \item[II.] $X$ has two $\DA_1$ singularities and contains $9$ lines. In this case, we let $E_i$ be for $i\in\{1,2\}$ be the exceptional divisors, $L_{12},$ $L_{i,j}$, $L_{12,ij}$ for $i,j\in\{1,2\}$  be the lines on $S$,
    \item[III.] $X$ has two $\DA_1$ singularities and contains $8$ lines. In this case, we let $E_i$ be for $i\in\{1,2\}$ be the exceptional divisors,  $L_{i,j}$, for $i\in\{1,2\}$, $j\in\{1,2,3,4\}$   be the lines on $S$,
     \item[IV.] $X$ has three $\DA_1$ singularities and contains $6$ lines. In this case, we let $E_1$, $E_1'$ and $E_2$ be the exceptional divisors, $L_{12}$, $L_{12}'$, $L_{1,i}$ and $L_{1,i}'$ for $i\in\{1,2\}$   be the lines on $S$, 
     \item[V.] $X$ has four $\DA_1$ singularities and contains $4$ lines. In this case, we let $E_i$ for $i\in\{1,2,3,4\}$ be the exceptional divisors, $L_{12}$, $L_{23}$, $L_{34}$, $L_{14}$  be the lines on $S$,
     \item[VI.] $X$ has an $\DA_2$ singularity and contains $8$ lines. In this case, we let $E$ and $E'$ be the exceptional divisors, $L_{i}$, $L_{i,1}$ and $L_{i,1}'$ for $i\in\{1,2\}$  be the lines on $S$,
     \item[VII.] $X$ has $\DA_2$ and $\DA_1$ singularities and contains $6$ lines. In this case, we let $E_1$, $E_1'$ and $E_2$ be the exceptional divisors, $L_{1,i}$, $L_{12,i}$ for $i\in\{1,2\}$, $L_2$ and $L_{12}$  be the lines on $S$,
     \item[VIII.] $X$ has $\DA_2$ and two $\DA_1$ singularities and contains $4$ lines. In this case, we let $E_i$ and $E_i'$ for $i\in\{1,2\}$ be the exceptional divisors, $L_{2}$, $L_2'$, $L_{12}$ and $L_{12}'$   be the lines on $S$,
     \item[IX.] $X$ has an $\DA_3$ singularity and contains $5$ lines. In this case, we let $E_1$, $E_1'$ and $E_2$ be the exceptional divisors, $L_2$, $L_{1}$, $L_1'$, $L_{1,1}$ and $L_{1,1}'$  be the lines on $S$,
    \item[X.] $X$ has an $\DA_3$ singularity and contains $4$ lines. In this case, we let
     $E_1$, $E_1'$ and $E_2$ be the exceptional divisors,  $L_{1,i}$ and $L_{1,i}'$ for $i\in\{1,2\}$  be the lines on $S$,
      \item[XI.] $X$ has  $\DA_3$ and $\DA_1$ singularities and contains $3$ lines. In this case, we let $E_1$, $E_1'$, $E_2$ and $E_3$ be the exceptional divisors, $L_{13}$, $L_{1,1}$ and $L_{1,2}$  be the lines on $S$,
      \item[XII.] $X$ has  $\DA_3$ and two $\DA_1$ singularities and contains $2$ lines. In this case, we let $E_1$, $E_1'$, $E_2$, $E_3$ and $E_3'$ be the exceptional divisors,  $L_{13}$ and $L_{13}'$ be the lines on $S$,
      \item[XIII.] $X$ has an $\DA_4$ singularity and contains $3$ lines. In this case, we let $E_i$ for $i\in\{1,2,3,4\}$ be the exceptional divisors, $L_2$, $L_{4}$ and $L_{4,1}$   be the lines on $S$,
      \item[XIV.] $X$ has an $\mathbb{D}_4$ singularity and contains $2$ lines. In this case, we let $E$, $E_1$, $E_1'$ and $E_2$ be the exceptional divisors,  $L_{1}$ and $L_{1}'$   be the lines on $S$,
      \item[XV.] $X$ has an $\mathbb{D}_5$ singularity and contains $1$ line. In this case, we let
      $E$ and $E_i$ for $i\in\{1,2,3,4\}$ be the exceptional divisors,  $L_{4}$ be a line on $S$.
\end{itemize}
such that the dual graph of the $(-1)$-curves  and $(-2)$-curves on $S$ is given on the picture below:
 \begin{figure}[h!]
    \centering
\hspace*{0cm}\includegraphics[width=16.5cm]{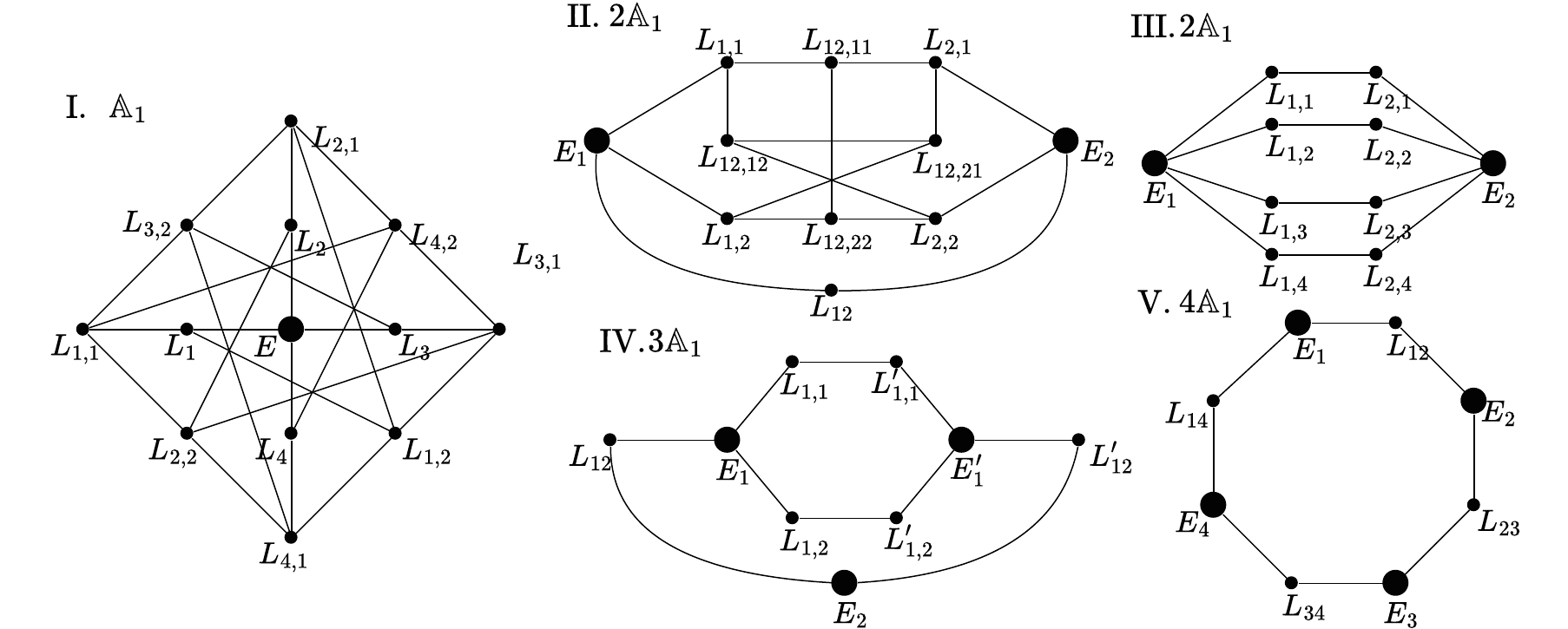}
\end{figure}
 \begin{figure}[h!]
    \centering
\hspace*{0cm}\includegraphics[width=16.5cm]{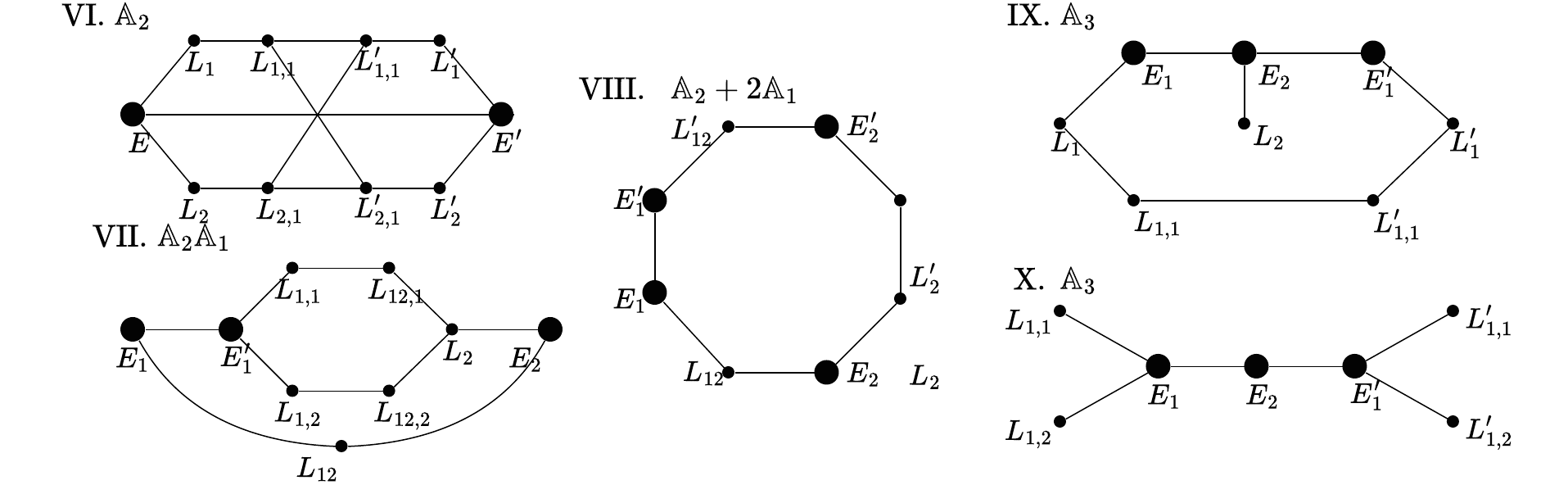}
\end{figure}
\begin{figure}[h!]
    \centering
\hspace*{0cm}\includegraphics[width=16cm]{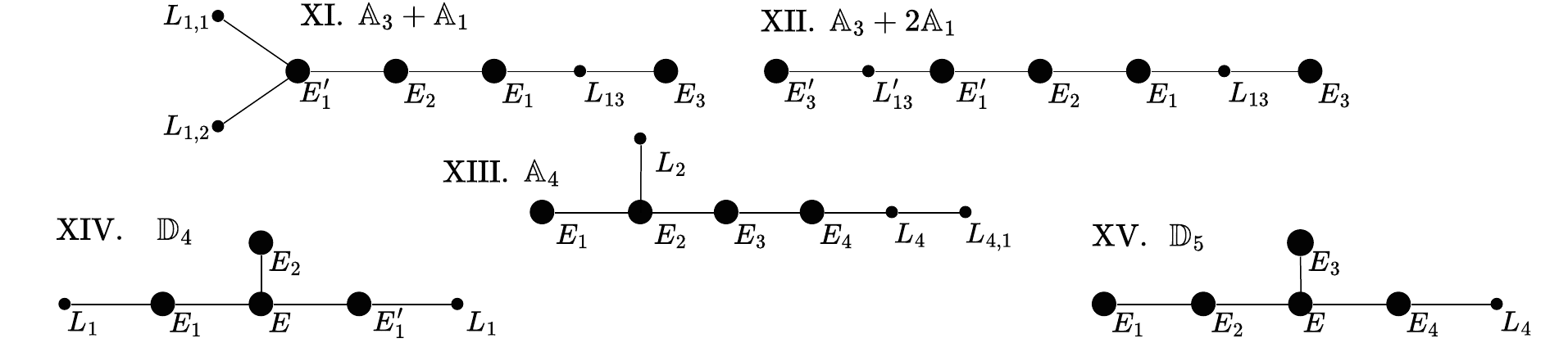}
\end{figure}

 One has:
{\hspace*{0cm}
\begin{itemize}
    \item[I.] $\delta(X)=1$ since depending on the position of point $P\in S$ we have 
    \begin{table}[h!]
 
\begin{tabular}{ | c || c | c | c | c |c |}
 \hline 
 $P$ & $E$ & $\big(\bigcup_{i\in\{1,2,3,4\}}L_{i}\big)\backslash E$ & 2 curves in $\mathbf{L}_1$ &  1 curve in $\mathbf{L}_1\backslash \big(\bigcup_{i\in\{1,2,3,4\}}L_{i}\big)$ & o/w \\\hline
$\delta_P(S)$ & $1$ & $\frac{6}{5}$ &$\frac{4}{3}$ & $\frac{18}{13}$ & $\frac{3}{2}$\\\hline
    \end{tabular}
\caption{Local $\delta$-invariants: $(-K_S)^2=4$ and  $\DA_1$ singularity}
\end{table}
\\ where $\mathbf{L}_1:= \big\{L_{i,j}\big| i\in\{1,2,3,4\}, j\in\{1,2\}\big\}$.
\item[II.]  $\delta(X)=1$ since depending on the position of point $P\in S$ we have
\begin{table}[h!]
 
\begin{tabular}{ | c || c | c | c | c |c |}
 \hline 
 $P$ & $ \mathbf{E}_2\cup L_{12}$ & $\big(\bigcup_{i,j\in\{1,2\}} L_{i,j} \big)\backslash \mathbf{E}_2$ &  $\mathbf{L}_2$ &   $\mathbf{L}_{12} \backslash \Big(\mathbf{L}_2\cup \bigcup_{i,j\in\{1,2\}} L_{i,j}\Big)$ & o/w \\\hline
$\delta_P(S)$ & $1$ & $\frac{6}{5}$ &$\frac{4}{3}$ & $\frac{18}{13}$ & $\frac{3}{2}$\\\hline
    \end{tabular}
\caption{Local $\delta$-invariants: $(-K_S)^2=4$ and  $2\DA_1$ singularities (9 lines)}
\end{table}
\\ where $\mathbf{E}_2:=E_1\cup E_2$,  $\mathbf{L}_2:= \big(L_{12, 11}\cap L_{12, 22}\big)\cup \big(L_{12, 12}\cap L_{12, 21}\big)$,  $\mathbf{L}_{12}:=\bigcup_{i,j\in\{1,2\}}L_{12,ij}$.
 \item[III.]  $\delta(X)=1$ since depending on the position of point $P\in S$ we have
 \begin{table}[h!]
 
\begin{tabular}{ | c || c | c | c | }
 \hline 
 $P$ & $E_1\cup E_2$ & $\big(\bigcup_{i\in\{1,2\},j\in\{1,2,3,4\}}L_{i,j}\big)\backslash ( E_1\cup E_2)$  & o/w \\\hline
$\delta_P(S)$ & $1$ & $\frac{6}{5}$  & $\frac{3}{2}$\\\hline
    \end{tabular}
\caption{Local $\delta$-invariants: $(-K_S)^2=4$ and  $2\DA_1$ singularities (8 lines)}
\end{table}
\newpage
 \item[IV.]  $\delta(X)=1$ since depending on the position of point $P\in S$ we have
  \begin{table}[h!]
 
\begin{tabular}{ | c || c | c | c | }
 \hline 
 $P$ & $E_1\cup E_1' \cup E_2\cup L_{12}\cup L_{12}'$ & $\big(\bigcup_{i\in \{1,2\}} L_{1,i}\cup L_{1,i}'\big)\backslash (E_1\cup E_1')$  & o/w \\\hline
$\delta_P(S)$ & $1$ & $\frac{6}{5}$  & $\frac{3}{2}$\\\hline
    \end{tabular}
\caption{Local $\delta$-invariants: $(-K_S)^2=4$ and  $3\DA_1$ singularities}
\end{table}
  \item[V.]  $\delta(X)=1$ since depending on the position of point $P\in S$ we have
  \begin{table}[h!]
\begin{tabular}{ | c || c | c |  }
 \hline 
 $P$ & $\big(\bigcup_{i\in\{1,2,3,4\}}E_i\big)\cup L_{12}\cup L_{23}\cup L_{34}\cup L_{14}$   & o/w \\\hline
$\delta_P(S)$ & $1$   & $\frac{3}{2}$\\\hline
    \end{tabular}
\caption{Local $\delta$-invariants: $(-K_S)^2=4$ and  $4\DA_1$ singularities}
\end{table}
\item[VI.]  $\delta(X)=\frac{6}{7}$ since depending on the position of point $P\in S$ we have
\begin{table}[h!]
 
\begin{tabular}{ | c || c | c | c | c |c |}
 \hline 
 $P$ & $\mathbf{E}_6$ & $\big(\bigcup_{i\in\{1,2\}} L_i\cup L_i'\big)\backslash \mathbf{E}_6$ &  $\mathbf{L}_6$ &   $\big(\bigcup_{i\in\{1,2\}} L_{i,1}\cup L_{i,1}'\big)\backslash \Big( \mathbf{L}_6 \cup \mathbf{E}_6\Big)$ & o/w \\\hline
$\delta_P(S)$ & $\frac{6}{7}$  & $\frac{8}{7}$ &$\frac{4}{3}$ & $\frac{18}{13}$ & $\frac{3}{2}$\\\hline
      \end{tabular} \\
       where $\mathbf{E}_6:=E\cup E'$, $\mathbf{L}_6:=  \big(L_{1,1}\cap L_{1,1}'\big)\cup \big(L_{1,1}\cap L_{2,1}'\big)\cup \big(L_{2,1}\cap L_{1,1}'\big)\cup \big(L_{2,1}\cap L_{2,1}'\big)$.
\caption{Local $\delta$-invariants: $(-K_S)^2=4$ and  $\DA_2$ singularity}
\end{table}
     \item[VII.] $\delta(X)=\frac{6}{7}$ since depending on the position of point $P\in S$ we have
 \begin{table}[h!]
 
\begin{tabular}{ | c || c | c | c | c |c | c |}
 \hline 
 $P$ & $E_1\cup E_1'\cup L_{12}$ & $E_2 \backslash L_{12}$ &  $\mathbf{L}_7\backslash E_1'$ &   $L_{2}\backslash E_2$ & $(L_{12,1}\cup L_{12,2})\backslash (\mathbf{L}_7\cup L_{2}) $ & o/w \\\hline
$\delta_P(S)$ & $\frac{6}{7}$ & $1$ & $\frac{8}{7}$ & $\frac{6}{5}$ & $\frac{18}{13}$ & $\frac{3}{2}$\\\hline
    \end{tabular}\\
     where $\mathbf{L}_7:=L_{1,1}\cup L_{1,2}$.
\caption{Local $\delta$-invariants: $(-K_S)^2=4$ and  $\DA_2\DA_1$ singularities}
\end{table}
\item[VIII.] $\delta(X)=\frac{6}{7}$ since depending on the position of point $P\in S$ we have
\begin{table}[h!]
 
\begin{tabular}{ | c || c | c | c | c |}
 \hline 
 $P$ & $E_1\cup E_1'\cup L_{12}\cup L_{12}'$ & $(E_2\cup E_2')\backslash (L_{12}\cup L_{12}')$ &  $(L_2\cup L_2')\backslash (E_2\cup E_2')$  & o/w \\\hline
$\delta_P(S)$ & $\frac{6}{7}$  & $1$ &$\frac{6}{5}$ &  $\frac{3}{2}$\\\hline
    \end{tabular}
\caption{Local $\delta$-invariants: $(-K_S)^2=4$ and  $\DA_22\DA_1$ singularities}
\end{table}
\item[IX.] $\delta(X)=\frac{3}{4}$ since depending on the position of point $P\in S$ we have
\begin{table}[h!]
 
\begin{tabular}{ | c || c | c | c | c | c | c | c |}
 \hline 
 $P$ & $E_2$ & $\mathbf{E}_9 \backslash E_2$ &  $(L_1\cup L_{1}') \backslash  \mathbf{E}_9$ & $L_{2}\backslash E_2$ & $\mathbf{L}_9$ & $(L_{1,1}\cup L_{1,1}')\backslash\mathbf{L}_9$ & o/w \\\hline
$\delta_P(S)$ & $\frac{2}{3}$  & $\frac{24}{29}$ & $\frac{12}{11}$ & $1$ &$\frac{4}{3}$ & $\frac{18}{13}$ &  $\frac{3}{2}$\\\hline
    \end{tabular}\\
     where $\mathbf{E}_9:=E_1\cup E_1'$, $\mathbf{L}_9:=L_{1,1}\cap L_{1,1}'$.
\caption{Local $\delta$-invariants: $(-K_S)^2=4$ and  $\DA_3$ singularity (5 lines)}
\end{table}
\newpage
\item[X.] $\delta(X)=\frac{3}{4}$ since depending on the position of point $P\in S$ we have
\begin{table}[h!]
\begin{tabular}{ | c || c | c | c |}
 \hline 
 $P$ & $E_1\cup E_1'\cup E_2$ & $\big(\bigcup_{i\in\{1,2\}} L_{1,i}\cup L_{1,i}' \big) \backslash  ( E_1\cup E_1')$  & o/w \\\hline
$\delta_P(S)$ & $\frac{3}{4}$  &$\frac{9}{8}$ &  $\frac{3}{2}$\\\hline
    \end{tabular}
\caption{Local $\delta$-invariants: $(-K_S)^2=4$ and  $\DA_3$ singularity (4 lines)}
\end{table}
\item[XI.]  $\delta(X)=\frac{3}{4}$ since depending on the position of point $P\in S$ we have
\begin{table}[h!]
 
\begin{tabular}{ | c || c | c | c | c |}
 \hline 
 $P$ & $E_1\cup E_1'\cup E_2 \cup L_{13}$ & $E_3 \backslash  L_{13}$  & $(L_{1,1}\cup L_{1,2})\backslash E_1'$ & o/w \\\hline
$\delta_P(S)$ & $\frac{3}{4}$ & $1$ &$\frac{9}{8}$ &  $\frac{3}{2}$\\\hline
    \end{tabular}
\caption{Local $\delta$-invariants: $(-K_S)^2=4$ and  $\DA_3\DA_1$ singularities}
\end{table}
  \item[XII.] $\delta(X)=\frac{3}{4}$ since
  \begin{table}[h!]
 
\begin{tabular}{ | c || c | c | c |}
 \hline 
 $P$ & $ E_1\cup E_1'\cup E_2 \cup L_{13} \cup L_{13}'$ & $(E_3\cup E_3')\backslash  ( L_{13} \cup L_{13}')$   & o/w \\\hline
$\delta_P(S)$ & $\frac{3}{4}$ & $1$ & $\frac{3}{2}$\\\hline
    \end{tabular}
\caption{Local $\delta$-invariants: $(-K_S)^2=4$ and  $\DA_32\DA_1$ singularities}
\end{table}
 \item[XIII.] $\delta(X)=\frac{6}{11}$ since depending on the position of point $P\in S$ we have
 \begin{table}[h!]
 
\begin{tabular}{ | c || c | c | c |c | c | c | c | c |}
 \hline 
 $P$ & $E_2$ & $E_3\backslash E_2$ & $E_4\backslash E_3$ & $L_2\backslash E_2$ & $E_1\backslash E_2$  & $L_4\backslash E_4$ & $L_{4,1}\backslash L_4$ & o/w \\\hline
$\delta_P(S)$ & $\frac{6}{11}$ & $\frac{24}{37}$ & $\frac{4}{5}$ & $\frac{24}{29}$ & $\frac{9}{11}$ & $\frac{24}{23}$ & $\frac{18}{13}$ & $\frac{3}{2}$\\\hline
    \end{tabular}
\caption{Local $\delta$-invariants: $(-K_S)^2=4$ and  $\DA_4$ singularity}
\end{table}
\item[XIV.] $\delta(X)=\frac{1}{2}$ since depending on the position of point $P\in S$ we have
 \begin{table}[h!]
 
\begin{tabular}{ | c || c | c | c | c | c |}
 \hline 
 $P$ & $E$ & $ (E_1\cup E_1')\backslash E$ & $E_2\backslash E$ & $(L_1\cup L_1')\backslash (E_1\cup E_1')$  & o/w \\\hline
$\delta_P(S)$ & $\frac{1}{2}$ & $\frac{2}{3}$ & $\frac{3}{4}$ & $1$ & $\frac{3}{2}$\\\hline
    \end{tabular}
\caption{Local $\delta$-invariants: $(-K_S)^2=4$ and  $\mathbb{D}_4$ singularity}
\end{table}
\item[XV.]  $\delta(X)=\frac{3}{8}$ since depending on the position of point $P\in S$ we have
\begin{table}[h!]
 
\begin{tabular}{ | c || c | c | c | c | c | }
 \hline 
 $P$ & $E_3$ & $(E_2\cup E_4)\backslash E_3$ & $ E\backslash  E_3$ & $(E_1\cup L_4)\backslash  (E_2\cup E_4)$   & o/w \\\hline
$\delta_P(S)$ & $\frac{3}{8}$ & $\frac{1}{2}$ & $\frac{9}{14}$ & $\frac{3}{4}$ & $\frac{3}{2}$\\\hline
    \end{tabular}
\caption{Local $\delta$-invariants: $(-K_S)^2=4$ and  $\mathbb{D}_5$ singularity}
\end{table}
\end{itemize}}
\begin{proof}
We prove each case separately using lemmas from the previous section.
    \begin{itemize}
        \item[I.] If $P\in E$, the assertion follows from Lemma \ref{deg4-A1points} [a).]. 
 If $P\in\big(\bigcup_{i\in\{1,2,3,4\}}L_{i}\big)\backslash E$, the assertion follows from Lemma \ref{deg4-nearA1points} [a).]. 
   If $P$ is the intersection of two $(-1)$-curves in  $\big\{L_{i,j}\big| i\in\{1,2,3,4\}, j\in\{1,2\}\big\}$, the assertion follows from Lemma \ref{deg4-twolines}. 
  If $P$  belongs to exactly one curve in  $\big\{L_{i,j}\big| i\in\{1,2,3,4\}, j\in\{1,2\}\big\} \backslash \big(\bigcup_{i\in\{1,2,3,4\}}L_{i}\big)$, the assertion follows from Lemma \ref{deg4-conic}.
  If $P$ is a general point, the assertion follows from Lemma \ref{deg4-generalpoint}.
   \item[II.] If $P\in E_1\cup E_2$, the assertion follows from Lemma \ref{deg4-A1points} [b).]. 
   If $P\in L_{12}\backslash (E_1\cup E_2)$, the assertion follows from Lemma \ref{deg4-near2A1points}.
 If $P\in \big(\bigcup_{i,j\in\{1,2\}} L_{i,j} \big)\backslash (E_1\cup E_2)$, the assertion follows from Lemma \ref{deg4-nearA1points} [a).]. 
  If $P\in \big(L_{12, 11}\cap L_{12, 22}\big)\cup \big(L_{12, 12}\cap L_{12, 21}\big) $, the assertion follows from Lemma \ref{deg4-twolines}. 
 If $ P\in \big(\bigcup_{i,j\in\{1,2\}}L_{12,ij}) \backslash \Big(\big(L_{12, 11}\cap L_{12, 22}\big)\cup \big(L_{12, 12}\cap L_{12, 21}\big)\cup \bigcup_{i,j\in\{1,2\}} L_{i,j}\Big)$, the assertion follows from Lemma \ref{deg4-conic}.
  If $P$ is a general point, the assertion follows from Lemma \ref{deg4-generalpoint}.
   \item[III.] If $P\in  E_1\cup E_2$, the assertion follows from Lemma \ref{deg4-A1points} [a).]. 
 If $P\in \big(\bigcup_{i\in\{1,2\},j\in\{1,2,3,4\}}L_{i,j}\big)\backslash ( E_1\cup E_2)$, the assertion follows from Lemma \ref{deg4-nearA1points} [b).]. 
   If $P$ is a general point, the assertion follows from Lemma \ref{deg4-generalpoint}.
   \item[IV.] If $P\in E_2$, the assertion follows from Lemma \ref{deg4-A1points} [c).].
 If $P\in E_1\cup E_1'$, the assertion follows from Lemma \ref{deg4-A1points} [b).].
  If $P\in (L_{12}\cup L_{12}')\backslash (E_1\cup E_1'\cup E_2)$, the assertion follows from Lemma \ref{deg4-near2A1points}.
 If $P\in \big(\bigcup_{i\in \{1,2\}} L_{1,i}\cup L_{1,i}'\big)\backslash (E_1\cup E_1')$, the assertion follows from Lemma \ref{deg4-nearA1points} [b).]. 
  If $P$ is a general point, the assertion follows from Lemma \ref{deg4-generalpoint}.
  \item[V.] If $P\in \bigcup_{i\in\{1,2,3,4\}}E_i$, the assertion follows from Lemma \ref{deg4-A1points} [c).].
 If $P\in (L_{12}\cup L_{23}\cup L_{34}\cup L_{14})\backslash \big(\bigcup_{i\in\{1,2,3,4\}}E_i\big)$, the assertion follows from Lemma \ref{deg4-near2A1points} [a).].
   If $P$ is a general point, the assertion follows from Lemma \ref{deg4-generalpoint}.
   \item[VI.] If $P\in E\cup E'$, the assertion follows from Lemma \ref{deg4-A2points} [a).].
 If $P\in \big(\bigcup_{i\in\{1,2\}} L_i\cup L_i'\big)\backslash (E\cup E')$, the assertion follows from Lemma \ref{deg4-nearA2points}.
  If $P\in \big(L_{1,1}\cap L_{1,1}'\big)\cup \big(L_{1,1}\cap L_{2,1}'\big)\cup \big(L_{2,1}\cap L_{1,1}'\big)\cup \big(L_{2,1}\cap L_{2,1}'\big)$, the assertion follows from Lemma \ref{deg4-twolines}.
 If $ P\in \big(\bigcup_{i\in\{1,2\}} L_{i,1}\cup L_{i,1}'\big)\backslash \Big( \big(L_{1,1}\cap L_{1,1}'\big)\cup \big(L_{1,1}\cap L_{2,1}'\big)\cup \big(L_{2,1}\cap L_{1,1}'\big)\cup \big(L_{2,1}\cap L_{2,1}'\big)\cup E\cup E'\Big)$, the assertion follows from Lemma \ref{deg4-conic}.
  If $P$ is a general point, the assertion follows from Lemma \ref{deg4-generalpoint}.
  \item[VII.] If $P\in E_1'$, the assertion follows from Lemma \ref{deg4-A2points} [a).].
 If $P\in E_1$, the assertion follows from Lemma \ref{deg4-A2points} [b).].
  If $P\in L_{12}\backslash E_1$, the assertion follows from Lemma \ref{deg4-67-2_3_points}.
 If $P\in E_2\backslash L_{12}$, the assertion follows from Lemma \ref{deg4-A1points} [d).].
 If $P\in (L_{1,1}\cup L_{1,2})\backslash E_1$, the assertion follows from Lemma \ref{deg4-nearA2points}.
 If $P\in L_{2}\backslash E_2$, the assertion follows from Lemma \ref{deg4-nearA1points} [a).]. 
 If $ P\in (L_{12,1}\cup L_{12,2})\backslash (L_{1,1}\cup L_{1,2}\cup L_{2})$, the assertion follows from Lemma \ref{deg4-conic}.
  If $P$ is a general point, the assertion follows from Lemma \ref{deg4-generalpoint}.
  \item[VIII.] If $P\in (L_{12}\cup L_{12}')\backslash (E_1\cup E_1')$, the assertion follows from Lemma \ref{deg4-67-2_3_points}.
  If $P\in E_1\cup E_1'$, the assertion follows from Lemma \ref{deg4-A2points} [b).].
  If $P\in (E_2\cup E_2')\backslash (L_{12}\cup L_{12}')$, the assertion follows from Lemma \ref{deg4-A1points} [d).].
 If $P\in  (L_2\cup L_2')\backslash (E_2\cup E_2')$, the assertion follows from Lemma \ref{deg4-nearA1points} [b).]. 
  If $P$ is a general point, the assertion follows from Lemma \ref{deg4-generalpoint}.
  \item[IX.] If $P\in E_2$, the assertion follows from Lemma \ref{deg4-centerA3points}.
 If $P\in (E_1\cup E_1')\backslash E_2$, the assertion follows from Lemma \ref{deg4-2429_1_32_2_points}.
  If $P\in (L_1\cup L_1')\backslash (E_1\cup E_1')$, the assertion follows from Lemma \ref{deg4-1211_1_2_points}.
 If $P\in L_{2}\backslash E_2$, the assertion follows from Lemma \ref{deg4-near2A1points} [b).].
 If $P= L_{1,1}\cap L_{1,1}'$, the assertion follows from Lemma \ref{deg4-twolines}.
 If $ P\in (L_{1,1}\cup L_{1,1}')\backslash (L_{1,1}\cap L_{1,1}')$, the assertion follows from Lemma \ref{deg4-conic}.
  If $P$ is a general point, the assertion follows from Lemma \ref{deg4-generalpoint}.
  \item[X.] If $P\in E_2$, the assertion follows from Lemma \ref{deg4-centerA3points}. 
 If $P\in (E_1\cup E_1')\backslash E_2$, the assertion follows from Lemma \ref{deg4-pointA3bound} [a).]. 
  If $P\in \big(\bigcup_{i\in\{1,2\}} L_{1,i}\cup L_{1,i}' \big) \backslash  ( E_1\cup E_1') $, the assertion follows from Lemma \ref{deg4-nearA3points}. 
  If $P$ is a general point, the assertion follows from Lemma \ref{deg4-generalpoint}.
  \item[XI.] If $P\in E_2$, the assertion follows from Lemma  \ref{deg4-23_1_2_3_points} [a).]. 
 If $P\in E_1'\backslash E_2$, the assertion follows from Lemma \ref{deg4-pointA3bound} [a).]. 
  If $P\in E_1\backslash E_2$, the assertion follows from Lemma \ref{deg4-pointA3bound} [b).]. 
 If $P\in L_{13}\backslash E_1$, the assertion follows from Lemma \ref{deg4-nearA1A3points} [a).]. 
 If $P\in E_3 \backslash  L_{13}$, the assertion follows from Lemma \ref{deg4-A1points} [e).].
 If $P\in (L_{1,1}\cup L_{1,2})\backslash E_1'$, the assertion follows from Lemma \ref{deg4-nearA3points}. 
  If $P$ is a general point, the assertion follows from Lemma \ref{deg4-generalpoint}.
    \item[XII.] If $P\in E_2$, the assertion follows from Lemma \ref{deg4-centerA3points}. 
 If $P\in (E_1\cup E_1')\backslash E_2$, the assertion follows from Lemma \ref{deg4-pointA3bound}. 
  If $P\in (L_{13} \cup L_{13}')\backslash (E_1\cup E_1')$, the assertion follows from Lemma \ref{deg4-nearA1A3points} [a).]. 
 If $P\in (E_3\cup E_3')\backslash  ( L_{13} \cup L_{13}')$, the assertion follows from Lemma \ref{deg4-A1points} [e).].
  If $P$ is a general point, the assertion follows from Lemma \ref{deg4-generalpoint}.
   \item[XIII.]  If $P\in E_2$, the assertion follows from Lemma \ref{deg4-611_1_3_4_points}. 
 If $P\in E_3\backslash E_2$, the assertion follows from Lemma \ref{deg4-2437_32_2_3_points}.
  If $P\in E_1\backslash E_2$, the assertion follows from Lemma \ref{deg4-911_43_2_points}.
 If $P\in E_4\backslash E_3$, the assertion follows from Lemma \ref{deg4-45_1_2_points}.
 If $P\in L_2\backslash E_2$, the assertion follows from Lemma \ref{deg4-2429_52_3_points}.
 If $P\in L_4\backslash E_4$, the assertion follows from Lemma \ref{deg4-2423_1_52_points}.
 If $ P\in L_{4,1}\backslash L_4$, the assertion follows from Lemma \ref{deg4-conic}.
  If $P$ is a general point, the assertion follows from Lemma \ref{deg4-generalpoint}.
  \item[XIV.] If $ P\in E$, the assertion follows from Lemma \ref{deg4-nearDpoints} [a).].
  If $ P\in (E_1\cup E_1')\backslash E$, the assertion follows from Lemma \ref{deg4-23_1_2_3_points} [b).].
  If $P\in E_2\backslash E$, the assertion follows from Lemma \ref{deg4-centerA3points} [b).].
 If $P\in (L_1\cup L_1')\backslash (E_1\cup E_1')$, the assertion follows from Lemma \ref{deg4-near2A1points} [c).].
  If $P$ is a general point, the assertion follows from Lemma \ref{deg4-generalpoint}.
  \item[XV.] If $ P\in E_3$, the assertion follows from Lemma \ref{deg4-38_2_6_points}.
  If $ P\in E_2\backslash E_3$, the assertion follows from Lemma \ref{deg4-nearDpoints} [b).].
  If $ P\in E_4\backslash E_3$, the assertion follows from Lemma \ref{deg4-12_1_5_points}.
 If $P\in E_2\backslash E_3$, the assertion follows from Lemma \ref{deg4-centerA3points} [c).].
 If $P\in L_{4}\backslash E_4$, the assertion follows from Lemma \ref{deg4-nearA1A3points} [b).]. 
 If $P\in E\backslash E_3$, the assertion follows from Lemma \ref{deg4-914_53_3_points}.
  If $P$ is a general point, the assertion follows from Lemma \ref{deg4-generalpoint}.
    \end{itemize}
\end{proof}

 \printbibliography
  
\end{document}